\documentclass[preprint,3p,sort&compress]{elsarticle}
\makeatletter
\def\ps@pprintTitle{%
 \let\@oddhead\@empty
 \let\@evenhead\@empty
 \def\@oddfoot{\centerline{\thepage}}%
 \let\@evenfoot\@oddfoot}
\makeatother

\journal{{\tt arXiv.org}}

\usepackage{graphicx}
\usepackage[dvipsnames]{xcolor}
\usepackage{amsmath,amsthm,amssymb}
\usepackage{bbm}
\usepackage{upgreek}
\usepackage{mathtools}
\usepackage{verbatim}
\usepackage{relsize}

\usepackage{subcaption}
\usepackage{booktabs}

\usepackage[linesnumbered,ruled]{algorithm2e}
\usepackage{multirow}

\usepackage[hypertexnames=false]{hyperref}
\usepackage[sort&compress,capitalize,noabbrev]{cleveref}

\usepackage{pgfplots}
\pgfplotsset{compat=newest}

\newcommand{\externaltikz}[2]{\includegraphics{Externals/#1}} 

\usepackage[symbol]{footmisc}

\title{\texorpdfstring{Distributed model order reduction of a model for microtubule-based cell polarization using HAPOD\footnote[1]{This is an early draft, numerical experiments are missing and there are probably many errors or inconsistencies.}}{
Distributed model order reduction of a model for microtubule-based cell polarization using HAPOD}}
\author[tl]{Tobias Leibner}
\author[mm]{Maja Matis}
\author[mo]{Mario Ohlberger}
\author[sr]{Stephan Rave}
\address[tl,mo.sr]{Faculty of Mathematics and Computer Science, University of Münster, Einsteinstrasse 62, 48149 Münster, {\tt tobias.leibner@uni-muenster.de, mario.ohlberger@uni-muenster.de, stephan.rave@uni-muenster.de}}
\address[mm]{Institute of Cell Biology, Center for Molecular Biology of Inflammation, University of Münster, Von-Esmarch-Strasse 56, 48149 Münster, {\tt matism@uni-muenster.de}}
\date{}

\def\quand{\quad \mbox{and} \quad}

\theoremstyle{plain}
\newtheorem{theorem}{Theorem}[section]

\newtheorem{remark}[theorem]{Remark}

\numberwithin{equation}{section}            
\numberwithin{figure}{section}
\numberwithin{table}{section}

\DeclareMathOperator*{\argmin}{argmin}

\DeclareMathOperator{\diverg}{div}
\DeclareMathOperator{\Vdiverg}{\mathbf{div}}

\newcommand*{\landauO}{O}

\newcommand*{\outernormal}{\Vn}

\newcommand{\pluseq}{\mathrel{+}=}
\newcommand{\multeq}{\mathrel{*}=}

\providecommand\given{}
\newcommand\SetSymbol[1][]{%
  \nonscript\:#1\vert{}
  \allowbreak{}
  \nonscript\:
  \mathopen{}}
\DeclarePairedDelimiterX\Set[1]\{\}{%
\renewcommand\given{\,\SetSymbol[\delimsize]\,}
\,#1\,
}

\newcommand*{\R}{\mathbb{R}}

\newcommand*{\Reals}[1][]{\R^{#1}}

\newcommand*{\Lp}[1]{L^{#1}}

\newcommand*{\LpError}[2][h]{E_{#1}^{#2}}
\newcommand*{\LpMeanErrorAbs}[1]{\overline{E}_{\mathrm{abs}}^{#1}}
\newcommand*{\LpMeanErrorRel}[1]{\overline{E}_{\mathrm{rel}}^{#1}}

\newcommand*{\Ltwo}{\Lp{2}}
\newcommand*{\Ltwozero}{\Lp{2}_0}
\newcommand*{\Hspace}[1]{H^{#1}}
\newcommand*{\Hone}{\Hspace{1}}
\newcommand*{\Honezero}{\Hspace{1}_0}
\newcommand*{\Hminusone}{\Hspace{-1}}

\newcommand*{\abs}[1]{\left\vert #1 \right\vert}
\newcommand*{\norm}[2][]{\left\lVert #2 \right\rVert_{#1}}

\newcommand*{\generalnorm}[1]{\left\lVert #1 \right\rVert}
\newcommand*{\hapodhilbertspace}{W}
\newcommand*{\hapodhilbertspacepfield}{W_{\hspace{-0.06cm}\pfieldvec}}
\newcommand*{\hapodhilbertspaceofield}{W_{\hspace{-0.05cm}\ofieldvec}}
\newcommand*{\hapodhilbertspacestokes}{W_{\hspace{-0.05cm}\stokesvec}}

\newcommand*{\Winnerproduct}[2]{{({#1},{#2})}_{\hapodhilbertspace}}
\newcommand*{\Winnerproductpfield}[2]{{({#1},{#2})}_{\hapodhilbertspacepfield}}
\newcommand*{\Winnerproductofield}[2]{{({#1},{#2})}_{\hapodhilbertspaceofield}}
\newcommand*{\Winnerproductstokes}[2]{{({#1},{#2})}_{\hapodhilbertspacestokes}}
\newcommand*{\Wnorm}[1]{\generalnorm{#1}_{\hapodhilbertspace}}

\newcommand*{\innerproductmatrix}{\MW}
\newcommand*{\innerproductmatrixpfield}{\MW_{\hspace{-0.06cm}\pfieldvec}}
\newcommand*{\innerproductmatrixofield}{\MW_{\hspace{-0.06cm}\ofieldvec}}
\newcommand*{\innerproductmatrixstokes}{\MW_{\hspace{-0.06cm}\stokesvec}}

\newcommand*{\Vzero}{\mathbf{0}} 
\newcommand*{\Mzero}{\mathbf{0}} 
\newcommand*{\massmatrix}{\MM}

\newcommand*{\ellipticmatrix}{\ME}

\newcommand*{\eyematrix}{\MI}

\newcommand*{\timeindex}{k}
\newcommand*{\numtimesteps}{K}
\newcommand*{\dimension}{d}
\newcommand*{\dimindex}{k}
\newcommand*{\feindex}{i}
\newcommand*{\feindexalt}{j}
\newcommand*{\newtonindex}{r}
\newcommand*{\paramindex}{m}

\newcommand*{\chunkindex}{a} 

\newcommand*{\dt}{\partial_\timevar}
\newcommand*{\gradient}{\nabla}
\newcommand*{\Vgradient}{\boldsymbol{\nabla}}

\newcommand*{\laplace}{\Delta}
\newcommand*{\Vlaplace}{\boldsymbol{\Delta}}
\newcommand*{\partialderiv}[2]{\frac{\partial #1}{\partial #2}}
\newcommand*{\variationalderiv}[2]{\frac{\delta #1}{\delta #2}}
\newcommand*{\normalderiv}[1]{\frac{\partial #1}{\partial n}}
\newcommand*{\normalderivvec}[1]{\frac{\partial #1}{\partial \outernormal}}
\newcommand*{\jacobianop}{D}
\newcommand*{\doubledotproduct}[2]{#1\mathbin{:}#2}

\newcommand*{\x}{x}

\newcommand*{\domain}{X}
\newcommand*{\domainboundary}{\partial \domain}
\newcommand*{\spatialvar}{\Vx}

\newcommand*{\gridwidth}{h}

\newcommand*{\timevar}{t}
\newcommand*{\timestep}{\Delta\timevar}

\newcommand*{\tf}{t_f} 

\makeatletter
\@tfor\next:=abcdefghijklmnopqrstuvwxyzPQ\do{%
  \def\command@factory#1{%
    \expandafter\def\csname V#1\endcsname{\mathbf{#1}}
  }
  \expandafter\command@factory\next{}
}
\makeatother

\makeatletter
\def\greekvectors#1{%
\@for\next:=#1\do{%
\def\X##1;{%
\expandafter\def\csname V##1\endcsname{\boldsymbol{\csname##1\endcsname}}
}
\expandafter\X\next;
}
}
\makeatother

\greekvectors{alpha,beta,delta,iota,gamma,lambda,theta,mu,nu,eta,chi,phi,psi,Gamma,Chi,varsigma,varphi}

\makeatletter
\@tfor\next:=ABCDEFGHIJKLMOPQRSTUVWXYZ\do{%
  \def\command@factory#1{%
    \expandafter\def\csname M#1\endcsname{\mathbf{#1}}
  }
  \expandafter\command@factory\next{}
}
\makeatother

\makeatletter
\def\greekmatrices#1{%
\@for\next:=#1\do{%
\def\X##1;{%
\expandafter\def\csname M##1\endcsname{\boldsymbol{\csname##1\endcsname}}
}
\expandafter\X\next;
}
}
\makeatother

\greekmatrices{Gamma,Chi,Omega,Sigma,Psi,Phi,Theta}

\newcommand*{\param}{\Veta}

\newcommand*{\fedimension}{n}
\newcommand*{\snapshotcount}{s}

\newcommand*{\reddimension}{N}

\newcommand*{\paramspace}{\mathcal{P}}
\newcommand*{\trainparams}{\mathcal{P}_{train}}
\newcommand*{\ntrain}{n_{train}}
\newcommand*{\subspace}{V}

\newcommand{\Vredp}{V_{\pfieldvec}}
\newcommand{\Vredo}{V_{\ofieldvec}}
\newcommand{\Vreds}{V_{\stokesvec}}

\newcommand{\Bredp}{\mathbf{V}_{\pfieldvec}}
\newcommand{\Bredo}{\mathbf{V}_{\ofieldvec}}
\newcommand{\Breds}{\mathbf{V}_{\stokesvec}}
\newcommand{\Credp}{\mathbf{C}_{\pfieldvec}}
\newcommand{\Credo}{\mathbf{C}_{\ofieldvec}}
\newcommand{\Creds}{\mathbf{C}_{\stokesvec}}
\newcommand{\redpfieldvec}{\overline{\pfieldvec}}
\newcommand{\redofieldvec}{\overline{\ofieldvec}}
\newcommand{\redstokesvec}{\overline{\stokesvec}}
\newcommand{\redpfieldvecdimension}{{N_{\pfieldvec}}}
\newcommand{\redofieldvecdimension}{{N_{\ofieldvec}}}
\newcommand{\redstokesvecdimension}{{N_{\stokesvec}}}
\newcommand{\colpfieldvecdimension}{{M_{\pfieldvec}}}
\newcommand{\colofieldvecdimension}{{M_{\ofieldvec}}}
\newcommand{\colstokesvecdimension}{{M_{\stokesvec}}}
\newcommand{\pfielddofmatrix}{\mathbf{Z}_{\pfieldvec}}
\newcommand{\ofielddofmatrix}{\mathbf{Z}_{\ofieldvec}}
\newcommand{\stokesdofmatrix}{\mathbf{Z}_{\stokesvec}}
\newcommand{\pfielddof}[1]{{d_{\pfieldvec,#1}}}
\newcommand{\ofielddof}[1]{{d_{\ofieldvec,#1}}}
\newcommand{\stokesdof}[1]{{d_{\stokesvec,#1}}}

\newcommand*{\fevec}{\Vy}

\newcommand*{\fespace}[1][p]{V_{\gridwidth}^{#1}}
\newcommand*{\fespacezero}[1][p]{V_{\gridwidth,0}^{#1}}

\newcommand*{\snapshot}{\Vv}

\newcommand*{\snapshotmatrix}{\MS}

\newcommand*{\snapshotset}{\mathcal{S}}

\newcommand*{\snapshotsetpfieldi}[1]{\snapshotset_{\pfieldvec,#1}}
\newcommand*{\snapshotsetofieldi}[1]{\snapshotset_{\ofieldvec,#1}}
\newcommand*{\snapshotsetstokesi}[1]{\snapshotset_{\stokesvec,#1}}
\newcommand*{\snapshotsetany}{\snapshotset_{\anyvec}}
\newcommand*{\snapshotsetanyi}[1]{\snapshotset_{\anyvec,#1}}
\newcommand*{\residualset}{\mathcal{R}}

\newcommand*{\residualsetpfieldi}[1]{\residualset_{\pfieldvec,#1}}
\newcommand*{\residualsetofieldi}[1]{\residualset_{\ofieldvec,#1}}
\newcommand*{\residualsetstokesi}[1]{\residualset_{\stokesvec,#1}}
\newcommand*{\residualsetany}{\residualset_{\anyvec}}
\newcommand*{\residualsetanyi}[1]{\residualset_{\anyvec,#1}}

\newcommand*{\newtonitspfield}{n_{\mathrm{it},\pfieldvec}}
\newcommand*{\newtonitsofield}{n_{\mathrm{it},\ofieldvec}}
\newcommand*{\newtonitsstokes}{n_{\mathrm{it},\stokesvec}}
\newcommand*{\newtonitsany}{n_{\mathrm{it},\anyvec}}

\newcommand*{\orthprojection}[1][]{P_{#1}}

\newcommand*{\orthprojectionmatrix}{\MP}

\newcommand*{\chunksize}{l}
\newcommand*{\chunksizelast}{l_{last}}
\newcommand*{\numchunks}{n_{\mathrm{chunks}}}

\newcommand*{\hapodtree}{\mathcal{T}}
\newcommand{\T}{\mathcal{T}}
\newcommand*{\hapodtreedepth}{L_{\hapodtree}}

\newcommand{\Snap}{\mathcal{S}}
\newcommand{\AccSnap}{\widetilde{\mathcal{S}}}

\newcommand{\POD}{\operatorname{POD}}

\newcommand{\HAPOD}{\operatorname{HAPOD}}

\newcommand{\children}[2][\T]{\mathcal{C}_{#1}(#2)} 
\newcommand{\childmap}[1][\T]{\mathcal{C}_{#1}}

\newcommand*{\intvar}[1]{\,\mathrm{d}#1}
\newcommand*{\bulkintegral}[1]{\int\limits_{\domain}\,{ #1 }\,\intvar{x}}

\newcommand*{\boundaryintegral}[1]{\int_{\partial\domain}\,{ #1 }\,\intvar{x}}
\newcommand*{\boundaryint}{\int_{\partial\domain}\,}

\newcommand*{\signeddistance}{r}

\newcommand*{\polaritynumber}{\text{Pa}}
\newcommand*{\ofielddoublewellparam}{c_1}
\newcommand*{\bendingcapillarynumber}{\text{Be}}
\newcommand*{\capillarynumber}{\text{Ca}}
\newcommand*{\fracrotationaldynamicviscosity}{\kappa}
\newcommand*{\pfieldmobility}{\gamma}
\newcommand*{\pfieldepsilon}{\epsilon}

\newcommand*{\vorticitytensor}{\boldsymbol{\Omega}}
\newcommand*{\deformationtensor}{\MD}
\newcommand*{\shapefactor}{\xi}
\newcommand*{\reynoldsnumber}{\text{Re}}

\newcommand*{\activeforcenumber}{\text{Fa}}

\newcommand*{\freeenergy}{E}
\newcommand*{\filamentenergy}{E_d}
\newcommand*{\surfaceenergy}{E_S}
\newcommand*{\kineticenergy}{E_{kin}}

\newcommand*{\doublewellpotential}{W}
\newcommand*{\doublewellpotentialprime}{\doublewellpotential^{\prime}}
\newcommand*{\doublewellpotentialtwoprime}{\doublewellpotential^{\prime\prime}}
\newcommand*{\doublewellpotentialthreeprime}{\doublewellpotential^{\prime\prime\prime}}

\newcommand*{\pfield}{\phi}
\newcommand*{\pfieldinitial}{\pfield_0}
\newcommand*{\pfieldphinat}{\pfield^{\natural}}
\newcommand*{\pfieldcellindicator}{\tilde{\pfield}}
\newcommand*{\Vphinat}{\Vphi^{\natural}}
\newcommand*{\pfieldmu}{\mu}
\newcommand*{\pfieldvec}{\Vvarphi}
\newcommand*{\anyvec}{\Valpha}

\newcommand*{\firstpfieldtest}{\chi}
\newcommand*{\secondpfieldtest}{\chi^{\natural}}
\newcommand*{\thirdpfieldtest}{\nu}
\newcommand*{\pfieldfebasisfunc}{\varphi}
\newcommand*{\pfieldpolynomialorder}{p}

\newcommand*{\pfieldfedimension}{n_{\gridwidth}^{\pfieldpolynomialorder}}
\newcommand*{\pfieldvecdimension}{{n_{\pfieldvec}}}
\newcommand*{\pfieldfecoeff}{\underline{\pfield}}
\newcommand*{\pfieldphinatfecoeff}{\underline{\pfieldphinat}}
\newcommand*{\pfieldmufecoeff}{\underline{\pfieldmu}}
\newcommand*{\pfieldfecoeffvec}{\underline{\Vphi}}

\newcommand*{\pfieldphinatfecoeffvec}{\underline{\Vphinat}}
\newcommand*{\pfieldphinatfecoeffvectime}[1]{\underline{\Vphi}^{\natural,#1}}

\newcommand*{\pfieldmufecoeffvec}{\underline{\Vmu}}

\newcommand*{\pfieldmassmatrix}{\massmatrix_{pf}}
\newcommand*{\pfieldmassmatrixentry}[2]{M_{#1#2}}
\newcommand*{\pfieldellipticmatrix}{\ellipticmatrix}
\newcommand*{\pfieldellipticmatrixentry}[2]{K_{#1#2}}

\newcommand*{\pfieldBmatrix}{\MB}
\newcommand*{\pfieldBmatrixentry}[2]{B_{#1#2}}

\newcommand*{\pfieldfe}{\pfield_{\gridwidth}}
\newcommand*{\pfieldphinatfe}{\pfieldphinat_{\gridwidth}}
\newcommand*{\pfieldphinatfetime}[1]{\pfield_{\gridwidth}^{\natural,#1}}
\newcommand*{\pfieldmufe}{\pfieldmu_{\gridwidth}}
\newcommand*{\pfieldferhsvec}{\Va}
\newcommand*{\pfieldferhsvecentry}[1]{a_{#1}}
\newcommand*{\pfieldfenonlinf}{\Vf}
\newcommand*{\pfieldfenonlinfentry}[1]{f_{#1}}
\newcommand*{\pfieldfenonling}{\Vg}
\newcommand*{\pfieldfenonlingentry}[1]{g_{#1}}

\newcommand*{\pfieldjacobianmatrix}{\MJ_{pf}}
\newcommand*{\pfieldjacobianlinearpreconditioner}{\hat{\MJ}_{pf}}
\newcommand*{\pfieldresidual}{\Vr_{\pfieldvec}}
\newcommand*{\anyresidual}{\Vr_{\anyvec}}

\newcommand*{\ofield}{\Vd}
\newcommand*{\ofieldvec}{\Vo}

\newcommand*{\ofieldinitial}{\ofield_0}
\newcommand*{\ofieldPnat}{\ofield^{\natural}}
\newcommand*{\firstofieldtest}{\Vq}
\newcommand*{\secondofieldtest}{\firstofieldtest^{\natural}}
\newcommand*{\ofieldfe}{\ofield_{\gridwidth}}
\newcommand*{\ofieldfebasisfunc}{\Vvarphi}

\newcommand*{\ofieldfedimension}{{(\pfieldfedimension)}^{\dimension}}
\newcommand*{\ofieldvecdimension}{{n_{\ofieldvec}}}
\newcommand*{\ofieldfecoeff}{\underline{d}}
\newcommand*{\ofieldPnatfecoeff}{\underline{d}^{\natural}}
\newcommand*{\ofieldfecoeffvec}{\underline{\ofield}}
\newcommand*{\ofieldPnatfecoeffvec}{\underline{\ofieldPnat}}
\newcommand*{\ofieldPnatfecoeffvectime}[1]{\underline{\Vd}^{\natural,#1}}

\newcommand*{\ofieldmassmatrix}{\massmatrix_{of}}
\newcommand*{\stokesmassmatrix}{\massmatrix_{st}}
\newcommand*{\ofieldmassmatrixentry}[2]{M_{#1#2}}
\newcommand*{\ofieldBmatrix}{\MB}
\newcommand*{\ofieldBmatrixentry}[2]{B_{#1#2}}
\newcommand*{\ofieldCmatrix}{\MC}
\newcommand*{\ofieldGmatrix}{\MG}
\newcommand*{\ofieldCmatrixentry}[2]{C_{#1#2}}
\newcommand*{\ofieldellipticmatrix}{\ellipticmatrix}
\newcommand*{\ofieldellipticmatrixentry}[2]{K_{#1#2}}
\newcommand*{\ofieldPnatfe}{\ofieldPnat_{\gridwidth}}
\newcommand*{\ofieldPnatfetime}[1]{\ofield_{\gridwidth}^{\natural,#1}}
\newcommand*{\ofieldfenonlinf}{\Vf}
\newcommand*{\ofieldfenonlinfentry}[1]{f_{#1}}
\newcommand*{\ofieldjacobianmatrix}{\MJ_{of}}

\newcommand*{\ofieldresidual}{\Vr_{\ofieldvec}}

\newcommand*{\stokesvelocity}{\Vu}
\newcommand*{\stokesvec}{\Vs}
\newcommand*{\stokespressurepolynomialorder}{p}

\newcommand*{\stokesvelocityfedimension}{{(n_{h,0}^{\stokespressurepolynomialorder+1})}^{\dimension}}
\newcommand*{\stokesvecdimension}{{n_{\stokesvec}}}
\newcommand*{\stokesvelocityfebasisfunc}{\Vpsi}
\newcommand*{\stokesvelocityfe}{\stokesvelocity_{\gridwidth}}
\newcommand*{\stokesvelocityfecoeff}{\underline{u}}
\newcommand*{\stokesvelocityfecoeffvec}{\underline{\Vu}}
\newcommand*{\stokesvelocityinitial}{\stokesvelocity_0}

\newcommand*{\stokesferhsvec}{\Va}
\newcommand*{\stokesferhsvecentry}[1]{{a_{#1}}}
\newcommand*{\stokespressure}{p}
\newcommand*{\stokespressurefe}{\stokespressure_{\gridwidth}}
\newcommand*{\stokespressurefecoeff}{\underline{\stokespressure}}
\newcommand*{\stokespressurefecoeffvec}{\underline{\Vp}}

\newcommand*{\viscousstress}{\boldsymbol{\sigma}_{viscous}}
\newcommand*{\ericksenstress}{\boldsymbol{\sigma}_{ericksen}}
\newcommand*{\activestress}{\boldsymbol{\sigma}_{active}}
\newcommand*{\distortionstress}{\boldsymbol{\sigma}_{dist}}
\newcommand*{\firststokestest}{\Vv}
\newcommand*{\secondstokestest}{q}
\newcommand*{\stokesAmatrix}{\MA}
\newcommand*{\stokesAmatrixentry}[2]{A_{#1#2}}
\newcommand*{\stokesBmatrix}{\MB}
\newcommand*{\stokesBmatrixentry}[2]{B_{#1#2}}
\newcommand*{\stokesresidual}{\Vr_{\stokesvec}}

\definecolor{greenyellow}   {cmyk}{0.15, 0, 0.69, 0}
\definecolor{yellow}        {cmyk}{0, 0, 1, 0}
\definecolor{goldenrod}     {cmyk}{0, 0.10, 0.84, 0}
\definecolor{dandelion}     {cmyk}{0, 0.29, 0.84, 0}
\definecolor{apricot}       {cmyk}{0, 0.32, 0.52, 0}
\definecolor{peach}         {cmyk}{0, 0.50, 0.70, 0}
\definecolor{melon}         {cmyk}{0, 0.46, 0.50, 0}
\definecolor{yelloworange}  {cmyk}{0, 0.42, 1, 0}
\definecolor{orange}        {cmyk}{0, 0.61, 0.87, 0}
\definecolor{burntorange}   {cmyk}{0, 0.51, 1, 0}
\definecolor{bittersweet}   {cmyk}{0, 0.75, 1, 0.24}
\definecolor{redorange}     {cmyk}{0, 0.77, 0.87, 0}
\definecolor{mahogany}      {cmyk}{0, 0.85, 0.87, 0.35}
\definecolor{maroon}        {cmyk}{0, 0.87, 0.68, 0.32}
\definecolor{brickred}      {cmyk}{0, 0.89, 0.94, 0.28}
\definecolor{red}           {cmyk}{0, 1, 1, 0}
\definecolor{orangered}     {cmyk}{0, 1, 0.50, 0}
\definecolor{rubinered}     {cmyk}{0, 1, 0.13, 0}
\definecolor{wildstrawberry}{cmyk}{0, 0.96, 0.39, 0}
\definecolor{salmon}        {cmyk}{0, 0.53, 0.38, 0}
\definecolor{carnationpink} {cmyk}{0, 0.63, 0, 0}
\definecolor{magenta}       {cmyk}{0, 1, 0, 0}
\definecolor{violetred}     {cmyk}{0, 0.81, 0, 0}
\definecolor{rhodamine}     {cmyk}{0, 0.82, 0, 0}
\definecolor{mulberry}      {cmyk}{0.34, 0.90, 0, 0.02}
\definecolor{redviolet}     {cmyk}{0.07, 0.90, 0, 0.34}
\definecolor{fuchsia}       {cmyk}{0.47, 0.91, 0, 0.08}
\definecolor{lavender}      {cmyk}{0, 0.48, 0, 0}
\definecolor{thistle}       {cmyk}{0.12, 0.59, 0, 0}
\definecolor{orchid}        {cmyk}{0.32, 0.64, 0, 0}
\definecolor{darkorchid}    {cmyk}{0.40, 0.80, 0.20, 0}
\definecolor{purple}        {cmyk}{0.45, 0.86, 0, 0}
\definecolor{plum}          {cmyk}{0.50, 1, 0, 0}
\definecolor{violet}        {cmyk}{0.79, 0.88, 0, 0}
\definecolor{royalpurple}   {cmyk}{0.75, 0.90, 0, 0}
\definecolor{blueviolet}    {cmyk}{0.86, 0.91, 0, 0.04}
\definecolor{periwinkle}    {cmyk}{0.57, 0.55, 0, 0}
\definecolor{cadetblue}     {cmyk}{0.62, 0.57, 0.23, 0}
\definecolor{cornflowerblue}{cmyk}{0.65, 0.13, 0, 0}
\definecolor{midnightblue}  {cmyk}{0.98, 0.13, 0, 0.43}
\definecolor{navyblue}      {cmyk}{0.94, 0.54, 0, 0}
\definecolor{royalblue}     {cmyk}{1, 0.50, 0, 0}
\definecolor{blue}          {cmyk}{1, 1, 0, 0}
\definecolor{cerulean}      {cmyk}{0.94, 0.11, 0, 0}
\definecolor{cyan}          {cmyk}{1, 0, 0, 0}
\definecolor{processblue}   {cmyk}{0.96, 0, 0, 0}
\definecolor{skyblue}       {cmyk}{0.62, 0, 0.12, 0}
\definecolor{turquoise}     {cmyk}{0.85, 0, 0.20, 0}
\definecolor{tealblue}      {cmyk}{0.86, 0, 0.34, 0.02}
\definecolor{aquamarine}    {cmyk}{0.82, 0, 0.30, 0}
\definecolor{bluegreen}     {cmyk}{0.85, 0, 0.33, 0}
\definecolor{emerald}       {cmyk}{1, 0, 0.50, 0}
\definecolor{junglegreen}   {cmyk}{0.99, 0, 0.52, 0}
\definecolor{seagreen}      {cmyk}{0.69, 0, 0.50, 0}
\definecolor{green}         {cmyk}{1, 0, 1, 0}
\definecolor{forestgreen}   {cmyk}{0.91, 0, 0.88, 0.12}
\definecolor{pinegreen}     {cmyk}{0.92, 0, 0.59, 0.25}
\definecolor{limegreen}     {cmyk}{0.50, 0, 1, 0}
\definecolor{yellowgreen}   {cmyk}{0.44, 0, 0.74, 0}
\definecolor{springgreen}   {cmyk}{0.26, 0, 0.76, 0}
\definecolor{olivegreen}    {cmyk}{0.64, 0, 0.95, 0.40}
\definecolor{rawsienna}     {cmyk}{0, 0.72, 1, 0.45}
\definecolor{sepia}         {cmyk}{0, 0.83, 1, 0.70}
\definecolor{brown}         {cmyk}{0, 0.81, 1, 0.60}
\definecolor{tan}           {cmyk}{0.14, 0.42, 0.56, 0}
\definecolor{gray}          {cmyk}{0, 0, 0, 0.50}
\definecolor{black}         {cmyk}{0, 0, 0, 1}
\definecolor{white}         {cmyk}{0, 0, 0, 0}

\pgfplotscreateplotcyclelist{color parula}{%
  royalblue!70,every mark/.append style={solid,line width = 0pt,fill=royalblue!60!black},mark=ball\\%
  goldenrod!50!yelloworange,every mark/.append style={solid,fill=goldenrod!30!black},mark=*\\%
  green,every mark/.append style={black,fill=yellowgreen},mark=diamond*\\%
  bluegreen,every mark/.append style={fill=black},mark=triangle*\\%
  royalblue!70,densely dashed,every mark/.append style={solid,line width = 0pt,fill=royalblue!60!black},mark=ball\\%
  goldenrod!50!yelloworange,densely dashed,every mark/.append style={solid,black,fill=goldenrod},mark=diamond*\\%
  yellowgreen,densely dashed,every mark/.append style={solid,fill=yellowgreen!60!royalblue},mark=square*, mark size= 1pt\\%
  bluegreen,densely dashed,mark size = 1.5pt,every mark/.append style={solid,fill=white},mark=otimes*\\%
  royalblue,densely dashed,mark size = 1.5pt,every mark/.append style={solid,fill=royalblue!30!black},mark=pentagon*\\%
  goldenrod!40!yelloworange,every mark/.append style={solid,fill=goldenrod!30!black},mark=|\\%
}

\pgfplotscreateplotcyclelist{color parula pairwise}{%
royalblue!70,solid,every mark/.append style={solid,line width = 0pt,fill=royalblue!60!black},mark=ball\\%
royalblue!70,densely dashed,every mark/.append style={solid,line width = 0pt,fill=royalblue!60!black},mark=*\\%
goldenrod!50!yelloworange,solid,every mark/.append
style={solid,fill=goldenrod!30!black},mark=diamond*\\%
goldenrod!50!yelloworange,densely dashed, every mark/.append
style={solid,fill=goldenrod!30!black},mark=triangle*\\%
green,solid,every mark/.append style={solid, black,fill=yellowgreen},mark=square*\\%
green,densely dashed,every mark/.append style={solid, black, fill=yellowgreen},mark=pentagon*\\%
bluegreen,every mark/.append style={fill=black},mark=triangle*\\%
blueviolet,densely dashed,every mark/.append style={solid,fill=blueviolet},mark=square*, mark size= 1pt\\%
bluegreen,densely dashed,mark size = 1.5pt,every mark/.append style={solid,fill=white},mark=otimes*\\%
royalblue,densely dashed,mark size = 1.5pt,every mark/.append style={solid,fill=royalblue!30!black},mark=pentagon*\\%
goldenrod!40!yelloworange,every mark/.append style={solid,fill=goldenrod!30!black},mark=|\\%
}

\pgfplotsset{
  discard if not/.style 2 args={
      x filter/.code={
          \edef\tempa{\thisrow{#1}}
          \edef\tempb{#2}
          \ifx\tempa\tempb{}
          \else
            
          \fi
        }
    }
}

\usepgfplotslibrary{groupplots}
\newlength{\figurewidth}
\newlength{\figureheight}
\newlength{\figurewidthfortwo}
\newlength{\figureheightforthree}
\pgfplotsset{%
  compat=1.11,
  colormap={parula}{%
      rgb(0pt)=(0.2081,0.1663,0.5292); 
      rgb(1pt)=(0.208355,0.16778,0.532238); 
      rgb(2pt)=(0.208611,0.169261,0.535275); 
      rgb(3pt)=(0.208866,0.170741,0.538313); 
      rgb(4pt)=(0.209121,0.172222,0.54135); 
      rgb(5pt)=(0.209376,0.173702,0.544388); 
      rgb(6pt)=(0.209632,0.175183,0.547425); 
      rgb(7pt)=(0.209887,0.176663,0.550463); 
      rgb(8pt)=(0.210134,0.178144,0.553505); 
      rgb(9pt)=(0.210338,0.179624,0.556568); 
      rgb(10pt)=(0.210542,0.181105,0.559631); 
      rgb(11pt)=(0.210746,0.182585,0.562694); 
      rgb(12pt)=(0.210944,0.184066,0.565763); 
      rgb(13pt)=(0.211123,0.185546,0.568852); 
      rgb(14pt)=(0.211302,0.187027,0.57194); 
      rgb(15pt)=(0.21148,0.188507,0.575029); 
      rgb(16pt)=(0.211642,0.189996,0.578117); 
      rgb(17pt)=(0.21177,0.191502,0.581206); 
      rgb(18pt)=(0.211897,0.193008,0.584295); 
      rgb(19pt)=(0.212025,0.194514,0.587383); 
      rgb(20pt)=(0.212132,0.19602,0.590472); 
      rgb(21pt)=(0.212208,0.197526,0.59356); 
      rgb(22pt)=(0.212285,0.199032,0.596649); 
      rgb(23pt)=(0.212361,0.200538,0.599738); 
      rgb(24pt)=(0.212413,0.202044,0.602839); 
      rgb(25pt)=(0.212438,0.20355,0.605953); 
      rgb(26pt)=(0.212464,0.205056,0.609067); 
      rgb(27pt)=(0.212489,0.206562,0.612181); 
      rgb(28pt)=(0.212471,0.208083,0.61531); 
      rgb(29pt)=(0.21242,0.209614,0.61845); 
      rgb(30pt)=(0.212368,0.211146,0.621589); 
      rgb(31pt)=(0.212317,0.212677,0.624729); 
      rgb(32pt)=(0.212216,0.214209,0.627868); 
      rgb(33pt)=(0.212088,0.215741,0.631008); 
      rgb(34pt)=(0.211961,0.217272,0.634148); 
      rgb(35pt)=(0.211833,0.218804,0.637287); 
      rgb(36pt)=(0.211668,0.220354,0.640446); 
      rgb(37pt)=(0.211489,0.221911,0.643611); 
      rgb(38pt)=(0.21131,0.223468,0.646776); 
      rgb(39pt)=(0.211132,0.225025,0.649941); 
      rgb(40pt)=(0.210848,0.226603,0.653107); 
      rgb(41pt)=(0.210541,0.228186,0.656272); 
      rgb(42pt)=(0.210235,0.229768,0.659437); 
      rgb(43pt)=(0.209929,0.231351,0.662602); 
      rgb(44pt)=(0.209553,0.232934,0.665767); 
      rgb(45pt)=(0.20917,0.234516,0.668932); 
      rgb(46pt)=(0.208787,0.236099,0.672098); 
      rgb(47pt)=(0.208405,0.237681,0.675263); 
      rgb(48pt)=(0.20787,0.239289,0.678453); 
      rgb(49pt)=(0.207334,0.240897,0.681644); 
      rgb(50pt)=(0.206798,0.242505,0.684835); 
      rgb(51pt)=(0.206255,0.244114,0.688025); 
      rgb(52pt)=(0.205617,0.245722,0.691216); 
      rgb(53pt)=(0.204979,0.24733,0.694407); 
      rgb(54pt)=(0.204341,0.248938,0.697597); 
      rgb(55pt)=(0.203675,0.250554,0.700792); 
      rgb(56pt)=(0.202858,0.252213,0.704008); 
      rgb(57pt)=(0.202041,0.253872,0.707224); 
      rgb(58pt)=(0.201225,0.255531,0.710441); 
      rgb(59pt)=(0.200372,0.257184,0.713657); 
      rgb(60pt)=(0.199402,0.258818,0.716873); 
      rgb(61pt)=(0.198432,0.260452,0.720089); 
      rgb(62pt)=(0.197462,0.262085,0.723305); 
      rgb(63pt)=(0.196419,0.263735,0.726522); 
      rgb(64pt)=(0.195219,0.26542,0.729738); 
      rgb(65pt)=(0.19402,0.267105,0.732954); 
      rgb(66pt)=(0.19282,0.268789,0.73617); 
      rgb(67pt)=(0.191549,0.270474,0.739386); 
      rgb(68pt)=(0.19017,0.272159,0.742603); 
      rgb(69pt)=(0.188792,0.273843,0.745819); 
      rgb(70pt)=(0.187414,0.275528,0.749035); 
      rgb(71pt)=(0.1859,0.277237,0.752264); 
      rgb(72pt)=(0.184241,0.278973,0.755505); 
      rgb(73pt)=(0.182581,0.280709,0.758747); 
      rgb(74pt)=(0.180922,0.282444,0.761989); 
      rgb(75pt)=(0.179133,0.284209,0.765245); 
      rgb(76pt)=(0.177244,0.285996,0.768512); 
      rgb(77pt)=(0.175356,0.287783,0.77178); 
      rgb(78pt)=(0.173467,0.289569,0.775047); 
      rgb(79pt)=(0.171363,0.291406,0.778314); 
      rgb(80pt)=(0.169142,0.293269,0.781581); 
      rgb(81pt)=(0.166922,0.295132,0.784849); 
      rgb(82pt)=(0.164701,0.296996,0.788116); 
      rgb(83pt)=(0.162238,0.298934,0.791365); 
      rgb(84pt)=(0.159686,0.300899,0.794606); 
      rgb(85pt)=(0.157133,0.302865,0.797848); 
      rgb(86pt)=(0.15458,0.30483,0.80109); 
      rgb(87pt)=(0.151738,0.306858,0.804352); 
      rgb(88pt)=(0.148828,0.3089,0.80762); 
      rgb(89pt)=(0.145918,0.310942,0.810887); 
      rgb(90pt)=(0.143008,0.312984,0.814154); 
      rgb(91pt)=(0.139687,0.31514,0.81733); 
      rgb(92pt)=(0.136318,0.31731,0.820495); 
      rgb(93pt)=(0.132949,0.319479,0.82366); 
      rgb(94pt)=(0.129579,0.321649,0.826826); 
      rgb(95pt)=(0.125811,0.323918,0.829841); 
      rgb(96pt)=(0.122033,0.32619,0.832853); 
      rgb(97pt)=(0.118256,0.328462,0.835865); 
      rgb(98pt)=(0.114458,0.330737,0.838862); 
      rgb(99pt)=(0.110349,0.333059,0.841619); 
      rgb(100pt)=(0.106239,0.335382,0.844376); 
      rgb(101pt)=(0.102129,0.337705,0.847132); 
      rgb(102pt)=(0.0979874,0.340021,0.849835); 
      rgb(103pt)=(0.093648,0.342292,0.852209); 
      rgb(104pt)=(0.0893087,0.344564,0.854583); 
      rgb(105pt)=(0.0849694,0.346836,0.856957); 
      rgb(106pt)=(0.08063,0.349091,0.859234); 
      rgb(107pt)=(0.0762907,0.351286,0.861174); 
      rgb(108pt)=(0.0719514,0.353481,0.863114); 
      rgb(109pt)=(0.067612,0.355676,0.865053); 
      rgb(110pt)=(0.0633195,0.357817,0.866853); 
      rgb(111pt)=(0.0591333,0.359833,0.868333); 
      rgb(112pt)=(0.0549471,0.36185,0.869814); 
      rgb(113pt)=(0.050761,0.363866,0.871294); 
      rgb(114pt)=(0.0466838,0.365823,0.872626); 
      rgb(115pt)=(0.0427784,0.367687,0.873724); 
      rgb(116pt)=(0.038873,0.36955,0.874821); 
      rgb(117pt)=(0.0349676,0.371414,0.875919); 
      rgb(118pt)=(0.0315066,0.373217,0.876872); 
      rgb(119pt)=(0.0285456,0.374953,0.877664); 
      rgb(120pt)=(0.0255847,0.376688,0.878455); 
      rgb(121pt)=(0.0226237,0.378424,0.879246); 
      rgb(122pt)=(0.0202132,0.380061,0.879868); 
      rgb(123pt)=(0.0182477,0.381618,0.880353); 
      rgb(124pt)=(0.0162823,0.383175,0.880838); 
      rgb(125pt)=(0.0143168,0.384732,0.881323); 
      rgb(126pt)=(0.0127892,0.386241,0.881695); 
      rgb(127pt)=(0.0115129,0.387721,0.882001); 
      rgb(128pt)=(0.0102366,0.389202,0.882307); 
      rgb(129pt)=(0.00896036,0.390682,0.882614); 
      rgb(130pt)=(0.00812372,0.392089,0.88281); 
      rgb(131pt)=(0.00746006,0.393468,0.882963); 
      rgb(132pt)=(0.0067964,0.394846,0.883116); 
      rgb(133pt)=(0.00613273,0.396224,0.883269); 
      rgb(134pt)=(0.00581622,0.397562,0.88332); 
      rgb(135pt)=(0.00558649,0.398889,0.883346); 
      rgb(136pt)=(0.00535676,0.400217,0.883371); 
      rgb(137pt)=(0.00512703,0.401544,0.883397); 
      rgb(138pt)=(0.00516757,0.402804,0.883332); 
      rgb(139pt)=(0.00524414,0.404054,0.883256); 
      rgb(140pt)=(0.00532072,0.405305,0.883179); 
      rgb(141pt)=(0.0053973,0.406556,0.883103); 
      rgb(142pt)=(0.00572012,0.407757,0.882952); 
      rgb(143pt)=(0.00605195,0.408957,0.882799); 
      rgb(144pt)=(0.00638378,0.410157,0.882646); 
      rgb(145pt)=(0.00672643,0.411355,0.882489); 
      rgb(146pt)=(0.00728799,0.412529,0.882259); 
      rgb(147pt)=(0.00784955,0.413704,0.88203); 
      rgb(148pt)=(0.00841111,0.414878,0.8818); 
      rgb(149pt)=(0.00898919,0.416045,0.881564); 
      rgb(150pt)=(0.00967838,0.417168,0.881283); 
      rgb(151pt)=(0.0103676,0.418292,0.881002); 
      rgb(152pt)=(0.0110568,0.419415,0.880721); 
      rgb(153pt)=(0.011773,0.420532,0.880435); 
      rgb(154pt)=(0.0125898,0.42163,0.880129); 
      rgb(155pt)=(0.0134066,0.422728,0.879823); 
      rgb(156pt)=(0.0142234,0.423825,0.879516); 
      rgb(157pt)=(0.0150703,0.424915,0.879195); 
      rgb(158pt)=(0.0159892,0.425987,0.878838); 
      rgb(159pt)=(0.0169081,0.427059,0.87848); 
      rgb(160pt)=(0.017827,0.428132,0.878123); 
      rgb(161pt)=(0.0187748,0.429194,0.877746); 
      rgb(162pt)=(0.0197703,0.430241,0.877338); 
      rgb(163pt)=(0.0207658,0.431287,0.876929); 
      rgb(164pt)=(0.0217613,0.432334,0.876521); 
      rgb(165pt)=(0.0227802,0.43338,0.876113); 
      rgb(166pt)=(0.0238267,0.434427,0.875704); 
      rgb(167pt)=(0.0248733,0.435473,0.875296); 
      rgb(168pt)=(0.0259198,0.43652,0.874887); 
      rgb(169pt)=(0.0269802,0.437553,0.874451); 
      rgb(170pt)=(0.0280523,0.438574,0.873992); 
      rgb(171pt)=(0.0291243,0.439595,0.873532); 
      rgb(172pt)=(0.0301964,0.440616,0.873073); 
      rgb(173pt)=(0.0312844,0.441621,0.872614); 
      rgb(174pt)=(0.032382,0.442616,0.872154); 
      rgb(175pt)=(0.0334796,0.443612,0.871695); 
      rgb(176pt)=(0.0345772,0.444607,0.871235); 
      rgb(177pt)=(0.0357108,0.445603,0.870758); 
      rgb(178pt)=(0.0368595,0.446598,0.870273); 
      rgb(179pt)=(0.0380081,0.447594,0.869788); 
      rgb(180pt)=(0.0391568,0.448589,0.869303); 
      rgb(181pt)=(0.0402652,0.449565,0.868798); 
      rgb(182pt)=(0.0413628,0.450535,0.868287); 
      rgb(183pt)=(0.0424604,0.451505,0.867777); 
      rgb(184pt)=(0.043558,0.452474,0.867266); 
      rgb(185pt)=(0.0445889,0.453444,0.866756); 
      rgb(186pt)=(0.0456099,0.454414,0.866245); 
      rgb(187pt)=(0.0466309,0.455384,0.865735); 
      rgb(188pt)=(0.047652,0.456354,0.865224); 
      rgb(189pt)=(0.0486,0.457324,0.864714); 
      rgb(190pt)=(0.0495444,0.458294,0.864203); 
      rgb(191pt)=(0.0504889,0.459264,0.863692); 
      rgb(192pt)=(0.0514315,0.460234,0.863181); 
      rgb(193pt)=(0.0523249,0.461204,0.862645); 
      rgb(194pt)=(0.0532183,0.462174,0.862109); 
      rgb(195pt)=(0.0541117,0.463144,0.861573); 
      rgb(196pt)=(0.0549991,0.464111,0.861034); 
      rgb(197pt)=(0.0558414,0.465056,0.860472); 
      rgb(198pt)=(0.0566838,0.466,0.859911); 
      rgb(199pt)=(0.0575261,0.466944,0.859349); 
      rgb(200pt)=(0.0583532,0.467889,0.858793); 
      rgb(201pt)=(0.0591189,0.468833,0.858257); 
      rgb(202pt)=(0.0598847,0.469778,0.857721); 
      rgb(203pt)=(0.0606505,0.470722,0.857185); 
      rgb(204pt)=(0.0614018,0.471667,0.856641); 
      rgb(205pt)=(0.0621165,0.472611,0.85608); 
      rgb(206pt)=(0.0628312,0.473556,0.855518); 
      rgb(207pt)=(0.0635459,0.4745,0.854957); 
      rgb(208pt)=(0.064242,0.475444,0.854405); 
      rgb(209pt)=(0.0649057,0.476389,0.853868); 
      rgb(210pt)=(0.0655694,0.477333,0.853332); 
      rgb(211pt)=(0.066233,0.478278,0.852796); 
      rgb(212pt)=(0.0668625,0.479222,0.852249); 
      rgb(213pt)=(0.0674495,0.480167,0.851687); 
      rgb(214pt)=(0.0680366,0.481111,0.851126); 
      rgb(215pt)=(0.0686237,0.482056,0.850564); 
      rgb(216pt)=(0.0691838,0.483,0.850003); 
      rgb(217pt)=(0.0697198,0.483944,0.849441); 
      rgb(218pt)=(0.0702559,0.484889,0.84888); 
      rgb(219pt)=(0.0707919,0.485833,0.848318); 
      rgb(220pt)=(0.0712967,0.486778,0.847772); 
      rgb(221pt)=(0.0717817,0.487722,0.847236); 
      rgb(222pt)=(0.0722667,0.488667,0.8467); 
      rgb(223pt)=(0.0727517,0.489611,0.846164); 
      rgb(224pt)=(0.0732012,0.490573,0.845628); 
      rgb(225pt)=(0.0736351,0.491543,0.845092); 
      rgb(226pt)=(0.0740691,0.492513,0.844556); 
      rgb(227pt)=(0.074503,0.493483,0.84402); 
      rgb(228pt)=(0.0748973,0.494433,0.843484); 
      rgb(229pt)=(0.0752802,0.495378,0.842948); 
      rgb(230pt)=(0.0756631,0.496322,0.842412); 
      rgb(231pt)=(0.0760459,0.497267,0.841876); 
      rgb(232pt)=(0.0763631,0.498233,0.841362); 
      rgb(233pt)=(0.0766694,0.499203,0.840851); 
      rgb(234pt)=(0.0769757,0.500173,0.840341); 
      rgb(235pt)=(0.077282,0.501143,0.83983); 
      rgb(236pt)=(0.0775162,0.502137,0.83932); 
      rgb(237pt)=(0.0777459,0.503132,0.838809); 
      rgb(238pt)=(0.0779757,0.504128,0.838298); 
      rgb(239pt)=(0.0782042,0.505123,0.837789); 
      rgb(240pt)=(0.0783829,0.506093,0.837304); 
      rgb(241pt)=(0.0785616,0.507063,0.836819); 
      rgb(242pt)=(0.0787402,0.508033,0.836334); 
      rgb(243pt)=(0.0789135,0.509008,0.835851); 
      rgb(244pt)=(0.0790411,0.510029,0.835392); 
      rgb(245pt)=(0.0791688,0.51105,0.834932); 
      rgb(246pt)=(0.0792964,0.512071,0.834473); 
      rgb(247pt)=(0.0794048,0.513092,0.834018); 
      rgb(248pt)=(0.0794303,0.514113,0.833584); 
      rgb(249pt)=(0.0794559,0.515134,0.83315); 
      rgb(250pt)=(0.0794814,0.516155,0.832717); 
      rgb(251pt)=(0.0794862,0.517183,0.832289); 
      rgb(252pt)=(0.0794351,0.51823,0.831881); 
      rgb(253pt)=(0.0793841,0.519276,0.831473); 
      rgb(254pt)=(0.079333,0.520323,0.831064); 
      rgb(255pt)=(0.079255,0.521369,0.830665); 
      rgb(256pt)=(0.0791273,0.522416,0.830282); 
      rgb(257pt)=(0.0789997,0.523462,0.829899); 
      rgb(258pt)=(0.0788721,0.524509,0.829516); 
      rgb(259pt)=(0.0786889,0.525589,0.829156); 
      rgb(260pt)=(0.0784336,0.526712,0.828824); 
      rgb(261pt)=(0.0781784,0.527835,0.828492); 
      rgb(262pt)=(0.0779231,0.528958,0.82816); 
      rgb(263pt)=(0.077615,0.530081,0.827868); 
      rgb(264pt)=(0.0772577,0.531205,0.827613); 
      rgb(265pt)=(0.0769003,0.532328,0.827357); 
      rgb(266pt)=(0.0765429,0.533451,0.827102); 
      rgb(267pt)=(0.0761243,0.534589,0.826862); 
      rgb(268pt)=(0.0756649,0.535738,0.826632); 
      rgb(269pt)=(0.0752054,0.536886,0.826403); 
      rgb(270pt)=(0.0747459,0.538035,0.826173); 
      rgb(271pt)=(0.0742168,0.539219,0.825961); 
      rgb(272pt)=(0.0736553,0.540418,0.825756); 
      rgb(273pt)=(0.0730937,0.541618,0.825552); 
      rgb(274pt)=(0.0725321,0.542818,0.825348); 
      rgb(275pt)=(0.0718925,0.544037,0.825183); 
      rgb(276pt)=(0.0712288,0.545262,0.82503); 
      rgb(277pt)=(0.0705652,0.546487,0.824877); 
      rgb(278pt)=(0.0699015,0.547713,0.824723); 
      rgb(279pt)=(0.0691514,0.548938,0.824614); 
      rgb(280pt)=(0.0683856,0.550163,0.824511); 
      rgb(281pt)=(0.0676198,0.551388,0.824409); 
      rgb(282pt)=(0.0668541,0.552614,0.824307); 
      rgb(283pt)=(0.0660408,0.553886,0.824205); 
      rgb(284pt)=(0.065224,0.555162,0.824103); 
      rgb(285pt)=(0.0644072,0.556439,0.824001); 
      rgb(286pt)=(0.0635892,0.557715,0.823899); 
      rgb(287pt)=(0.0626703,0.558991,0.823848); 
      rgb(288pt)=(0.0617514,0.560268,0.823797); 
      rgb(289pt)=(0.0608324,0.561544,0.823746); 
      rgb(290pt)=(0.0599087,0.56282,0.823693); 
      rgb(291pt)=(0.0589387,0.564096,0.823616); 
      rgb(292pt)=(0.0579688,0.565373,0.82354); 
      rgb(293pt)=(0.0569988,0.566649,0.823463); 
      rgb(294pt)=(0.0560243,0.567925,0.823386); 
      rgb(295pt)=(0.0550288,0.569202,0.82331); 
      rgb(296pt)=(0.0540333,0.570478,0.823233); 
      rgb(297pt)=(0.0530378,0.571754,0.823157); 
      rgb(298pt)=(0.0520423,0.57303,0.82308); 
      rgb(299pt)=(0.0510468,0.574307,0.823004); 
      rgb(300pt)=(0.0500514,0.575583,0.822927); 
      rgb(301pt)=(0.0490559,0.576859,0.82285); 
      rgb(302pt)=(0.0480604,0.578127,0.822756); 
      rgb(303pt)=(0.0470649,0.579377,0.822629); 
      rgb(304pt)=(0.0460694,0.580628,0.822501); 
      rgb(305pt)=(0.0450739,0.581879,0.822374); 
      rgb(306pt)=(0.0441,0.583119,0.822235); 
      rgb(307pt)=(0.0431556,0.584344,0.822082); 
      rgb(308pt)=(0.0422111,0.585569,0.821929); 
      rgb(309pt)=(0.0412667,0.586795,0.821776); 
      rgb(310pt)=(0.0403351,0.58802,0.821597); 
      rgb(311pt)=(0.0394162,0.589245,0.821392); 
      rgb(312pt)=(0.0384973,0.59047,0.821188); 
      rgb(313pt)=(0.0375784,0.591695,0.820984); 
      rgb(314pt)=(0.0367495,0.592891,0.820735); 
      rgb(315pt)=(0.0359838,0.594065,0.820454); 
      rgb(316pt)=(0.035218,0.595239,0.820173); 
      rgb(317pt)=(0.0344523,0.596413,0.819892); 
      rgb(318pt)=(0.0337721,0.597553,0.819595); 
      rgb(319pt)=(0.0331339,0.598676,0.819288); 
      rgb(320pt)=(0.0324958,0.599799,0.818982); 
      rgb(321pt)=(0.0318577,0.600923,0.818676); 
      rgb(322pt)=(0.0312964,0.602026,0.818312); 
      rgb(323pt)=(0.0307604,0.603124,0.817929); 
      rgb(324pt)=(0.0302243,0.604222,0.817546); 
      rgb(325pt)=(0.0296883,0.605319,0.817163); 
      rgb(326pt)=(0.0292375,0.606395,0.816738); 
      rgb(327pt)=(0.0288036,0.607468,0.816304); 
      rgb(328pt)=(0.0283697,0.60854,0.81587); 
      rgb(329pt)=(0.0279357,0.609612,0.815436); 
      rgb(330pt)=(0.0275721,0.610637,0.814955); 
      rgb(331pt)=(0.0272147,0.611658,0.81447); 
      rgb(332pt)=(0.0268574,0.612679,0.813985); 
      rgb(333pt)=(0.0265,0.6137,0.8135); 
      rgb(334pt)=(0.0262447,0.614695,0.812964); 
      rgb(335pt)=(0.0259895,0.615691,0.812428); 
      rgb(336pt)=(0.0257342,0.616686,0.811892); 
      rgb(337pt)=(0.0254853,0.61768,0.811352); 
      rgb(338pt)=(0.0253066,0.61865,0.810765); 
      rgb(339pt)=(0.0251279,0.61962,0.810177); 
      rgb(340pt)=(0.0249492,0.62059,0.80959); 
      rgb(341pt)=(0.024779,0.621551,0.808995); 
      rgb(342pt)=(0.0246514,0.62247,0.808357); 
      rgb(343pt)=(0.0245237,0.623389,0.807719); 
      rgb(344pt)=(0.0243961,0.624308,0.80708); 
      rgb(345pt)=(0.0242748,0.625221,0.80643); 
      rgb(346pt)=(0.0241727,0.626114,0.805741); 
      rgb(347pt)=(0.0240706,0.627008,0.805051); 
      rgb(348pt)=(0.0239685,0.627901,0.804362); 
      rgb(349pt)=(0.0238832,0.628786,0.803656); 
      rgb(350pt)=(0.0238321,0.629654,0.802916); 
      rgb(351pt)=(0.0237811,0.630522,0.802176); 
      rgb(352pt)=(0.02373,0.631389,0.801435); 
      rgb(353pt)=(0.023679,0.632247,0.800685); 
      rgb(354pt)=(0.0236279,0.633089,0.799919); 
      rgb(355pt)=(0.0235769,0.633932,0.799153); 
      rgb(356pt)=(0.0235258,0.634774,0.798387); 
      rgb(357pt)=(0.0234748,0.635604,0.797596); 
      rgb(358pt)=(0.0234237,0.63642,0.79678); 
      rgb(359pt)=(0.0233727,0.637237,0.795963); 
      rgb(360pt)=(0.0233216,0.638054,0.795146); 
      rgb(361pt)=(0.0232706,0.638856,0.794329); 
      rgb(362pt)=(0.0232195,0.639647,0.793512); 
      rgb(363pt)=(0.0231685,0.640439,0.792695); 
      rgb(364pt)=(0.0231174,0.64123,0.791879); 
      rgb(365pt)=(0.0230832,0.642005,0.791011); 
      rgb(366pt)=(0.0230577,0.64277,0.790118); 
      rgb(367pt)=(0.0230321,0.643536,0.789225); 
      rgb(368pt)=(0.0230066,0.644302,0.788331); 
      rgb(369pt)=(0.0229811,0.645049,0.787438); 
      rgb(370pt)=(0.0229556,0.645789,0.786544); 
      rgb(371pt)=(0.02293,0.646529,0.785651); 
      rgb(372pt)=(0.0229045,0.647269,0.784758); 
      rgb(373pt)=(0.022858,0.64801,0.783843); 
      rgb(374pt)=(0.0228069,0.64875,0.782924); 
      rgb(375pt)=(0.0227559,0.64949,0.782005); 
      rgb(376pt)=(0.0227048,0.65023,0.781086); 
      rgb(377pt)=(0.0227,0.650947,0.780144); 
      rgb(378pt)=(0.0227,0.651662,0.7792); 
      rgb(379pt)=(0.0227,0.652377,0.778256); 
      rgb(380pt)=(0.0227,0.653092,0.777311); 
      rgb(381pt)=(0.0228261,0.653781,0.776341); 
      rgb(382pt)=(0.0229538,0.65447,0.775371); 
      rgb(383pt)=(0.0230814,0.655159,0.774402); 
      rgb(384pt)=(0.0232108,0.655849,0.77343); 
      rgb(385pt)=(0.023364,0.656538,0.772434); 
      rgb(386pt)=(0.0235171,0.657227,0.771439); 
      rgb(387pt)=(0.0236703,0.657916,0.770443); 
      rgb(388pt)=(0.0238312,0.658602,0.769444); 
      rgb(389pt)=(0.0240354,0.659265,0.768423); 
      rgb(390pt)=(0.0242396,0.659929,0.767402); 
      rgb(391pt)=(0.0244438,0.660592,0.766381); 
      rgb(392pt)=(0.0247021,0.661256,0.765354); 
      rgb(393pt)=(0.025136,0.66192,0.764307); 
      rgb(394pt)=(0.02557,0.662583,0.763261); 
      rgb(395pt)=(0.0260039,0.663247,0.762214); 
      rgb(396pt)=(0.0264541,0.663911,0.761168); 
      rgb(397pt)=(0.026939,0.664574,0.760121); 
      rgb(398pt)=(0.027424,0.665238,0.759074); 
      rgb(399pt)=(0.027909,0.665902,0.758028); 
      rgb(400pt)=(0.028445,0.666555,0.756971); 
      rgb(401pt)=(0.0290577,0.667193,0.755899); 
      rgb(402pt)=(0.0296703,0.667832,0.754827); 
      rgb(403pt)=(0.0302829,0.66847,0.753755); 
      rgb(404pt)=(0.030994,0.669095,0.752683); 
      rgb(405pt)=(0.0318108,0.669708,0.751611); 
      rgb(406pt)=(0.0326276,0.670321,0.750539); 
      rgb(407pt)=(0.0334444,0.670933,0.749467); 
      rgb(408pt)=(0.0343045,0.67156,0.748366); 
      rgb(409pt)=(0.0351979,0.672198,0.747243); 
      rgb(410pt)=(0.0360913,0.672837,0.74612); 
      rgb(411pt)=(0.0369847,0.673475,0.744996); 
      rgb(412pt)=(0.0380432,0.674096,0.743873); 
      rgb(413pt)=(0.0391919,0.674709,0.74275); 
      rgb(414pt)=(0.0403405,0.675322,0.741627); 
      rgb(415pt)=(0.0414892,0.675934,0.740504); 
      rgb(416pt)=(0.0427123,0.676528,0.739381); 
      rgb(417pt)=(0.0439631,0.677115,0.738258); 
      rgb(418pt)=(0.0452138,0.677702,0.737135); 
      rgb(419pt)=(0.0464646,0.678289,0.736011); 
      rgb(420pt)=(0.0477153,0.678897,0.734868); 
      rgb(421pt)=(0.0489661,0.67951,0.733719); 
      rgb(422pt)=(0.0502168,0.680123,0.73257); 
      rgb(423pt)=(0.0514676,0.680735,0.731422); 
      rgb(424pt)=(0.0529237,0.681325,0.73025); 
      rgb(425pt)=(0.0544042,0.681912,0.729076); 
      rgb(426pt)=(0.0558847,0.682499,0.727902); 
      rgb(427pt)=(0.0573652,0.683086,0.726728); 
      rgb(428pt)=(0.0587709,0.683673,0.725553); 
      rgb(429pt)=(0.0601748,0.68426,0.724379); 
      rgb(430pt)=(0.0615787,0.684847,0.723205); 
      rgb(431pt)=(0.0629946,0.685435,0.722028); 
      rgb(432pt)=(0.0646027,0.686022,0.720803); 
      rgb(433pt)=(0.0662108,0.686609,0.719577); 
      rgb(434pt)=(0.0678189,0.687196,0.718352); 
      rgb(435pt)=(0.069427,0.687779,0.717131); 
      rgb(436pt)=(0.0710351,0.688341,0.715931); 
      rgb(437pt)=(0.0726432,0.688902,0.714731); 
      rgb(438pt)=(0.0742514,0.689464,0.713532); 
      rgb(439pt)=(0.0758709,0.690026,0.712326); 
      rgb(440pt)=(0.07753,0.690587,0.711101); 
      rgb(441pt)=(0.0791892,0.691149,0.709876); 
      rgb(442pt)=(0.0808483,0.69171,0.70865); 
      rgb(443pt)=(0.0825387,0.692272,0.707417); 
      rgb(444pt)=(0.0843,0.692833,0.706167); 
      rgb(445pt)=(0.0860613,0.693395,0.704916); 
      rgb(446pt)=(0.0878225,0.693956,0.703665); 
      rgb(447pt)=(0.089564,0.694518,0.702405); 
      rgb(448pt)=(0.0912742,0.69508,0.701128); 
      rgb(449pt)=(0.0929844,0.695641,0.699852); 
      rgb(450pt)=(0.0946946,0.696203,0.698576); 
      rgb(451pt)=(0.0965009,0.696752,0.697299); 
      rgb(452pt)=(0.0984153,0.697288,0.696023); 
      rgb(453pt)=(0.10033,0.697824,0.694747); 
      rgb(454pt)=(0.102244,0.69836,0.693471); 
      rgb(455pt)=(0.10413,0.698896,0.69218); 
      rgb(456pt)=(0.105994,0.699432,0.690878); 
      rgb(457pt)=(0.107857,0.699968,0.689577); 
      rgb(458pt)=(0.10972,0.700505,0.688275); 
      rgb(459pt)=(0.111632,0.701041,0.686973); 
      rgb(460pt)=(0.113572,0.701577,0.685671); 
      rgb(461pt)=(0.115512,0.702113,0.684369); 
      rgb(462pt)=(0.117452,0.702649,0.683068); 
      rgb(463pt)=(0.119429,0.703185,0.681747); 
      rgb(464pt)=(0.12142,0.703721,0.68042); 
      rgb(465pt)=(0.123411,0.704257,0.679093); 
      rgb(466pt)=(0.125402,0.704793,0.677765); 
      rgb(467pt)=(0.127372,0.705308,0.676438); 
      rgb(468pt)=(0.129338,0.705819,0.675111); 
      rgb(469pt)=(0.131303,0.706329,0.673783); 
      rgb(470pt)=(0.133269,0.70684,0.672456); 
      rgb(471pt)=(0.135369,0.70735,0.671084); 
      rgb(472pt)=(0.137488,0.707861,0.669705); 
      rgb(473pt)=(0.139607,0.708371,0.668327); 
      rgb(474pt)=(0.141725,0.708882,0.666949); 
      rgb(475pt)=(0.143795,0.709392,0.665595); 
      rgb(476pt)=(0.145862,0.709903,0.664242); 
      rgb(477pt)=(0.14793,0.710414,0.662889); 
      rgb(478pt)=(0.150003,0.710924,0.661534); 
      rgb(479pt)=(0.152198,0.711435,0.66013); 
      rgb(480pt)=(0.154394,0.711945,0.658726); 
      rgb(481pt)=(0.156589,0.712456,0.657322); 
      rgb(482pt)=(0.158784,0.712963,0.655922); 
      rgb(483pt)=(0.160979,0.713448,0.654543); 
      rgb(484pt)=(0.163174,0.713933,0.653165); 
      rgb(485pt)=(0.16537,0.714418,0.651786); 
      rgb(486pt)=(0.16757,0.714908,0.650397); 
      rgb(487pt)=(0.169791,0.715419,0.648968); 
      rgb(488pt)=(0.172012,0.715929,0.647538); 
      rgb(489pt)=(0.174232,0.71644,0.646109); 
      rgb(490pt)=(0.176483,0.716935,0.64468); 
      rgb(491pt)=(0.178806,0.717395,0.64325); 
      rgb(492pt)=(0.181129,0.717854,0.641821); 
      rgb(493pt)=(0.183452,0.718314,0.640391); 
      rgb(494pt)=(0.185755,0.718783,0.638952); 
      rgb(495pt)=(0.188027,0.719268,0.637497); 
      rgb(496pt)=(0.190299,0.719753,0.636042); 
      rgb(497pt)=(0.192571,0.720238,0.634587); 
      rgb(498pt)=(0.194913,0.720711,0.633132); 
      rgb(499pt)=(0.197338,0.72117,0.631677); 
      rgb(500pt)=(0.199762,0.72163,0.630223); 
      rgb(501pt)=(0.202187,0.722089,0.628768); 
      rgb(502pt)=(0.204612,0.722549,0.627299); 
      rgb(503pt)=(0.207037,0.723008,0.625818); 
      rgb(504pt)=(0.209462,0.723468,0.624338); 
      rgb(505pt)=(0.211887,0.723927,0.622857); 
      rgb(506pt)=(0.214328,0.724386,0.621377); 
      rgb(507pt)=(0.216778,0.724846,0.619896); 
      rgb(508pt)=(0.219229,0.725305,0.618416); 
      rgb(509pt)=(0.221679,0.725765,0.616935); 
      rgb(510pt)=(0.224202,0.726188,0.615455); 
      rgb(511pt)=(0.226754,0.726597,0.613974); 
      rgb(512pt)=(0.229307,0.727005,0.612494); 
      rgb(513pt)=(0.231859,0.727414,0.611014); 
      rgb(514pt)=(0.234392,0.727842,0.609513); 
      rgb(515pt)=(0.236919,0.728276,0.608007); 
      rgb(516pt)=(0.239446,0.72871,0.606501); 
      rgb(517pt)=(0.241973,0.729144,0.604995); 
      rgb(518pt)=(0.244611,0.729556,0.603467); 
      rgb(519pt)=(0.247266,0.729964,0.601935); 
      rgb(520pt)=(0.24992,0.730372,0.600404); 
      rgb(521pt)=(0.252575,0.730781,0.598872); 
      rgb(522pt)=(0.25523,0.731189,0.597365); 
      rgb(523pt)=(0.257884,0.731598,0.595859); 
      rgb(524pt)=(0.260539,0.732006,0.594353); 
      rgb(525pt)=(0.263194,0.732414,0.592846); 
      rgb(526pt)=(0.265848,0.732796,0.591314); 
      rgb(527pt)=(0.268503,0.733179,0.589783); 
      rgb(528pt)=(0.271158,0.733562,0.588251); 
      rgb(529pt)=(0.27383,0.733945,0.58672); 
      rgb(530pt)=(0.276638,0.734328,0.585188); 
      rgb(531pt)=(0.279446,0.734711,0.583657); 
      rgb(532pt)=(0.282254,0.735094,0.582125); 
      rgb(533pt)=(0.285051,0.735471,0.580594); 
      rgb(534pt)=(0.287808,0.735829,0.579062); 
      rgb(535pt)=(0.290565,0.736186,0.577531); 
      rgb(536pt)=(0.293322,0.736544,0.575999); 
      rgb(537pt)=(0.2961,0.736894,0.574468); 
      rgb(538pt)=(0.298933,0.737226,0.572936); 
      rgb(539pt)=(0.301767,0.737557,0.571405); 
      rgb(540pt)=(0.3046,0.737889,0.569873); 
      rgb(541pt)=(0.307452,0.738221,0.568351); 
      rgb(542pt)=(0.310336,0.738553,0.566845); 
      rgb(543pt)=(0.313221,0.738885,0.565339); 
      rgb(544pt)=(0.316105,0.739217,0.563833); 
      rgb(545pt)=(0.318978,0.739537,0.562315); 
      rgb(546pt)=(0.321837,0.739843,0.560784); 
      rgb(547pt)=(0.324696,0.74015,0.559252); 
      rgb(548pt)=(0.327555,0.740456,0.557721); 
      rgb(549pt)=(0.330468,0.740749,0.556216); 
      rgb(550pt)=(0.333429,0.741029,0.554736); 
      rgb(551pt)=(0.336389,0.74131,0.553255); 
      rgb(552pt)=(0.33935,0.741591,0.551775); 
      rgb(553pt)=(0.342296,0.741872,0.550279); 
      rgb(554pt)=(0.345231,0.742153,0.548773); 
      rgb(555pt)=(0.348167,0.742433,0.547267); 
      rgb(556pt)=(0.351102,0.742714,0.545761); 
      rgb(557pt)=(0.354038,0.742977,0.54429); 
      rgb(558pt)=(0.356973,0.743232,0.542835); 
      rgb(559pt)=(0.359908,0.743488,0.54138); 
      rgb(560pt)=(0.362844,0.743743,0.539925); 
      rgb(561pt)=(0.365839,0.743959,0.53847); 
      rgb(562pt)=(0.368851,0.744163,0.537015); 
      rgb(563pt)=(0.371863,0.744367,0.53556); 
      rgb(564pt)=(0.374875,0.744571,0.534105); 
      rgb(565pt)=(0.377843,0.744775,0.532672); 
      rgb(566pt)=(0.380804,0.74498,0.531243); 
      rgb(567pt)=(0.383765,0.745184,0.529814); 
      rgb(568pt)=(0.386726,0.745388,0.528384); 
      rgb(569pt)=(0.389711,0.745568,0.527003); 
      rgb(570pt)=(0.392697,0.745747,0.525624); 
      rgb(571pt)=(0.395684,0.745926,0.524246); 
      rgb(572pt)=(0.39867,0.746104,0.522868); 
      rgb(573pt)=(0.401657,0.746257,0.521489); 
      rgb(574pt)=(0.404643,0.74641,0.520111); 
      rgb(575pt)=(0.40763,0.746563,0.518732); 
      rgb(576pt)=(0.410611,0.746716,0.517359); 
      rgb(577pt)=(0.413546,0.746869,0.516032); 
      rgb(578pt)=(0.416482,0.747023,0.514705); 
      rgb(579pt)=(0.419417,0.747176,0.513377); 
      rgb(580pt)=(0.422357,0.747319,0.512055); 
      rgb(581pt)=(0.425318,0.747421,0.510753); 
      rgb(582pt)=(0.428279,0.747523,0.509451); 
      rgb(583pt)=(0.43124,0.747626,0.50815); 
      rgb(584pt)=(0.43418,0.747735,0.506848); 
      rgb(585pt)=(0.437065,0.747862,0.505546); 
      rgb(586pt)=(0.439949,0.74799,0.504244); 
      rgb(587pt)=(0.442834,0.748117,0.502942); 
      rgb(588pt)=(0.445727,0.748227,0.501659); 
      rgb(589pt)=(0.448637,0.748304,0.500408); 
      rgb(590pt)=(0.451547,0.74838,0.499157); 
      rgb(591pt)=(0.454457,0.748457,0.497906); 
      rgb(592pt)=(0.457333,0.748522,0.496667); 
      rgb(593pt)=(0.460167,0.748573,0.495441); 
      rgb(594pt)=(0.463,0.748624,0.494216); 
      rgb(595pt)=(0.465833,0.748675,0.492991); 
      rgb(596pt)=(0.468667,0.748726,0.491779); 
      rgb(597pt)=(0.4715,0.748777,0.490579); 
      rgb(598pt)=(0.474333,0.748829,0.48938); 
      rgb(599pt)=(0.477167,0.74888,0.48818); 
      rgb(600pt)=(0.479969,0.748931,0.486995); 
      rgb(601pt)=(0.482752,0.748982,0.485821); 
      rgb(602pt)=(0.485534,0.749033,0.484647); 
      rgb(603pt)=(0.488316,0.749084,0.483473); 
      rgb(604pt)=(0.491081,0.7491,0.482316); 
      rgb(605pt)=(0.493838,0.7491,0.481168); 
      rgb(606pt)=(0.496595,0.7491,0.480019); 
      rgb(607pt)=(0.499351,0.7491,0.47887); 
      rgb(608pt)=(0.502069,0.74912,0.477722); 
      rgb(609pt)=(0.504775,0.749145,0.476573); 
      rgb(610pt)=(0.50748,0.749171,0.475424); 
      rgb(611pt)=(0.510186,0.749196,0.474276); 
      rgb(612pt)=(0.512892,0.7492,0.47317); 
      rgb(613pt)=(0.515598,0.7492,0.472073); 
      rgb(614pt)=(0.518303,0.7492,0.470975); 
      rgb(615pt)=(0.521009,0.7492,0.469877); 
      rgb(616pt)=(0.523644,0.749176,0.46878); 
      rgb(617pt)=(0.526273,0.749151,0.467682); 
      rgb(618pt)=(0.528902,0.749125,0.466585); 
      rgb(619pt)=(0.531531,0.7491,0.465487); 
      rgb(620pt)=(0.53416,0.749074,0.464415); 
      rgb(621pt)=(0.536789,0.749049,0.463343); 
      rgb(622pt)=(0.539418,0.749023,0.462271); 
      rgb(623pt)=(0.542043,0.748998,0.461199); 
      rgb(624pt)=(0.544621,0.748972,0.460127); 
      rgb(625pt)=(0.547199,0.748947,0.459055); 
      rgb(626pt)=(0.549777,0.748921,0.457983); 
      rgb(627pt)=(0.55235,0.748891,0.45692); 
      rgb(628pt)=(0.554903,0.74884,0.455899); 
      rgb(629pt)=(0.557456,0.748789,0.454878); 
      rgb(630pt)=(0.560008,0.748738,0.453857); 
      rgb(631pt)=(0.562554,0.748687,0.452829); 
      rgb(632pt)=(0.565081,0.748636,0.451783); 
      rgb(633pt)=(0.567608,0.748585,0.450736); 
      rgb(634pt)=(0.570135,0.748534,0.449689); 
      rgb(635pt)=(0.572653,0.748474,0.44866); 
      rgb(636pt)=(0.575155,0.748397,0.447665); 
      rgb(637pt)=(0.577656,0.748321,0.446669); 
      rgb(638pt)=(0.580158,0.748244,0.445674); 
      rgb(639pt)=(0.582649,0.748168,0.444678); 
      rgb(640pt)=(0.585125,0.748091,0.443683); 
      rgb(641pt)=(0.587601,0.748014,0.442687); 
      rgb(642pt)=(0.590077,0.747938,0.441692); 
      rgb(643pt)=(0.59254,0.747861,0.440709); 
      rgb(644pt)=(0.59499,0.747785,0.439739); 
      rgb(645pt)=(0.597441,0.747708,0.438769); 
      rgb(646pt)=(0.599891,0.747632,0.437799); 
      rgb(647pt)=(0.602311,0.747555,0.436814); 
      rgb(648pt)=(0.604711,0.747478,0.435819); 
      rgb(649pt)=(0.60711,0.747402,0.434823); 
      rgb(650pt)=(0.60951,0.747325,0.433828); 
      rgb(651pt)=(0.611909,0.747232,0.432867); 
      rgb(652pt)=(0.614308,0.747129,0.431922); 
      rgb(653pt)=(0.616708,0.747027,0.430978); 
      rgb(654pt)=(0.619107,0.746925,0.430033); 
      rgb(655pt)=(0.621487,0.746823,0.429089); 
      rgb(656pt)=(0.623861,0.746721,0.428144); 
      rgb(657pt)=(0.626235,0.746619,0.4272); 
      rgb(658pt)=(0.628609,0.746517,0.426256); 
      rgb(659pt)=(0.630962,0.746393,0.425311); 
      rgb(660pt)=(0.63331,0.746266,0.424367); 
      rgb(661pt)=(0.635658,0.746138,0.423422); 
      rgb(662pt)=(0.638007,0.746011,0.422478); 
      rgb(663pt)=(0.640332,0.745906,0.421557); 
      rgb(664pt)=(0.642654,0.745804,0.420638); 
      rgb(665pt)=(0.644977,0.745702,0.419719); 
      rgb(666pt)=(0.6473,0.7456,0.4188); 
      rgb(667pt)=(0.649623,0.745472,0.417881); 
      rgb(668pt)=(0.651946,0.745345,0.416962); 
      rgb(669pt)=(0.654268,0.745217,0.416043); 
      rgb(670pt)=(0.656587,0.745089,0.415124); 
      rgb(671pt)=(0.658859,0.744962,0.414205); 
      rgb(672pt)=(0.661131,0.744834,0.413286); 
      rgb(673pt)=(0.663402,0.744707,0.412368); 
      rgb(674pt)=(0.665674,0.744579,0.411453); 
      rgb(675pt)=(0.667946,0.744451,0.410559); 
      rgb(676pt)=(0.670218,0.744324,0.409666); 
      rgb(677pt)=(0.672489,0.744196,0.408773); 
      rgb(678pt)=(0.674755,0.744062,0.407879); 
      rgb(679pt)=(0.677001,0.743909,0.406986); 
      rgb(680pt)=(0.679247,0.743756,0.406092); 
      rgb(681pt)=(0.681494,0.743603,0.405199); 
      rgb(682pt)=(0.68374,0.743458,0.404306); 
      rgb(683pt)=(0.685986,0.74333,0.403412); 
      rgb(684pt)=(0.688232,0.743203,0.402519); 
      rgb(685pt)=(0.690479,0.743075,0.401626); 
      rgb(686pt)=(0.692704,0.742937,0.400732); 
      rgb(687pt)=(0.694899,0.742784,0.399839); 
      rgb(688pt)=(0.697094,0.742631,0.398945); 
      rgb(689pt)=(0.699289,0.742477,0.398052); 
      rgb(690pt)=(0.701497,0.742324,0.397171); 
      rgb(691pt)=(0.703718,0.742171,0.396303); 
      rgb(692pt)=(0.705939,0.742018,0.395435); 
      rgb(693pt)=(0.708159,0.741865,0.394568); 
      rgb(694pt)=(0.710351,0.741712,0.3937); 
      rgb(695pt)=(0.71252,0.741559,0.392832); 
      rgb(696pt)=(0.71469,0.741405,0.391964); 
      rgb(697pt)=(0.71686,0.741252,0.391096); 
      rgb(698pt)=(0.719029,0.741082,0.390228); 
      rgb(699pt)=(0.721199,0.740904,0.38936); 
      rgb(700pt)=(0.723369,0.740725,0.388492); 
      rgb(701pt)=(0.725538,0.740546,0.387625); 
      rgb(702pt)=(0.727708,0.740386,0.386757); 
      rgb(703pt)=(0.729878,0.740233,0.385889); 
      rgb(704pt)=(0.732047,0.74008,0.385021); 
      rgb(705pt)=(0.734217,0.739927,0.384153); 
      rgb(706pt)=(0.736366,0.739753,0.383285); 
      rgb(707pt)=(0.73851,0.739574,0.382417); 
      rgb(708pt)=(0.740654,0.739395,0.38155); 
      rgb(709pt)=(0.742798,0.739217,0.380682); 
      rgb(710pt)=(0.744919,0.739038,0.379837); 
      rgb(711pt)=(0.747038,0.738859,0.378995); 
      rgb(712pt)=(0.749156,0.738681,0.378152); 
      rgb(713pt)=(0.751275,0.738502,0.37731); 
      rgb(714pt)=(0.753394,0.738323,0.376442); 
      rgb(715pt)=(0.755512,0.738145,0.375574); 
      rgb(716pt)=(0.757631,0.737966,0.374707); 
      rgb(717pt)=(0.75975,0.737789,0.373841); 
      rgb(718pt)=(0.761868,0.737636,0.372998); 
      rgb(719pt)=(0.763987,0.737483,0.372156); 
      rgb(720pt)=(0.766105,0.73733,0.371314); 
      rgb(721pt)=(0.76822,0.737169,0.370471); 
      rgb(722pt)=(0.770313,0.736965,0.369629); 
      rgb(723pt)=(0.772406,0.73676,0.368786); 
      rgb(724pt)=(0.774499,0.736556,0.367944); 
      rgb(725pt)=(0.776592,0.736358,0.367096); 
      rgb(726pt)=(0.778686,0.736179,0.366228); 
      rgb(727pt)=(0.780779,0.736001,0.36536); 
      rgb(728pt)=(0.782872,0.735822,0.364492); 
      rgb(729pt)=(0.784957,0.735643,0.363632); 
      rgb(730pt)=(0.787024,0.735465,0.36279); 
      rgb(731pt)=(0.789092,0.735286,0.361948); 
      rgb(732pt)=(0.791159,0.735107,0.361105); 
      rgb(733pt)=(0.793227,0.734929,0.360263); 
      rgb(734pt)=(0.795295,0.73475,0.359421); 
      rgb(735pt)=(0.797362,0.734571,0.358578); 
      rgb(736pt)=(0.79943,0.734392,0.357736); 
      rgb(737pt)=(0.801485,0.734214,0.356881); 
      rgb(738pt)=(0.803527,0.734035,0.356014); 
      rgb(739pt)=(0.805569,0.733856,0.355146); 
      rgb(740pt)=(0.807611,0.733678,0.354278); 
      rgb(741pt)=(0.809668,0.733499,0.353424); 
      rgb(742pt)=(0.811735,0.73332,0.352582); 
      rgb(743pt)=(0.813803,0.733142,0.35174); 
      rgb(744pt)=(0.81587,0.732963,0.350897); 
      rgb(745pt)=(0.817921,0.732784,0.350038); 
      rgb(746pt)=(0.819963,0.732606,0.349171); 
      rgb(747pt)=(0.822005,0.732427,0.348303); 
      rgb(748pt)=(0.824047,0.732248,0.347435); 
      rgb(749pt)=(0.826071,0.73207,0.346567); 
      rgb(750pt)=(0.828087,0.731891,0.345699); 
      rgb(751pt)=(0.830104,0.731712,0.344831); 
      rgb(752pt)=(0.83212,0.731534,0.343963); 
      rgb(753pt)=(0.834158,0.731355,0.343095); 
      rgb(754pt)=(0.8362,0.731176,0.342228); 
      rgb(755pt)=(0.838242,0.730998,0.34136); 
      rgb(756pt)=(0.840284,0.730819,0.340492); 
      rgb(757pt)=(0.842303,0.73064,0.339624); 
      rgb(758pt)=(0.84432,0.730462,0.338756); 
      rgb(759pt)=(0.846336,0.730283,0.337888); 
      rgb(760pt)=(0.848353,0.730104,0.33702); 
      rgb(761pt)=(0.850369,0.729926,0.336153); 
      rgb(762pt)=(0.852386,0.729747,0.335285); 
      rgb(763pt)=(0.854402,0.729568,0.334417); 
      rgb(764pt)=(0.856419,0.729391,0.333546); 
      rgb(765pt)=(0.858435,0.729238,0.332627); 
      rgb(766pt)=(0.860452,0.729085,0.331708); 
      rgb(767pt)=(0.862468,0.728932,0.330789); 
      rgb(768pt)=(0.864481,0.728778,0.329874); 
      rgb(769pt)=(0.866472,0.728625,0.32898); 
      rgb(770pt)=(0.868463,0.728472,0.328087); 
      rgb(771pt)=(0.870454,0.728319,0.327194); 
      rgb(772pt)=(0.872445,0.728166,0.326295); 
      rgb(773pt)=(0.874436,0.728013,0.325376); 
      rgb(774pt)=(0.876427,0.727859,0.324457); 
      rgb(775pt)=(0.878418,0.727706,0.323538); 
      rgb(776pt)=(0.880417,0.727561,0.322619); 
      rgb(777pt)=(0.882433,0.727433,0.3217); 
      rgb(778pt)=(0.88445,0.727306,0.320781); 
      rgb(779pt)=(0.886466,0.727178,0.319862); 
      rgb(780pt)=(0.888463,0.72705,0.318933); 
      rgb(781pt)=(0.890429,0.726923,0.317989); 
      rgb(782pt)=(0.892394,0.726795,0.317044); 
      rgb(783pt)=(0.894359,0.726668,0.3161); 
      rgb(784pt)=(0.896337,0.726552,0.315132); 
      rgb(785pt)=(0.898328,0.72645,0.314136); 
      rgb(786pt)=(0.900319,0.726348,0.313141); 
      rgb(787pt)=(0.90231,0.726246,0.312145); 
      rgb(788pt)=(0.904301,0.726158,0.31115); 
      rgb(789pt)=(0.906292,0.726081,0.310154); 
      rgb(790pt)=(0.908283,0.726005,0.309159); 
      rgb(791pt)=(0.910274,0.725928,0.308163); 
      rgb(792pt)=(0.912249,0.725851,0.307151); 
      rgb(793pt)=(0.914214,0.725775,0.30613); 
      rgb(794pt)=(0.91618,0.725698,0.305109); 
      rgb(795pt)=(0.918145,0.725622,0.304088); 
      rgb(796pt)=(0.920111,0.7256,0.303031); 
      rgb(797pt)=(0.922076,0.7256,0.301959); 
      rgb(798pt)=(0.924041,0.7256,0.300886); 
      rgb(799pt)=(0.926007,0.7256,0.299814); 
      rgb(800pt)=(0.927972,0.7256,0.298722); 
      rgb(801pt)=(0.929938,0.7256,0.297624); 
      rgb(802pt)=(0.931903,0.7256,0.296527); 
      rgb(803pt)=(0.933869,0.7256,0.295429); 
      rgb(804pt)=(0.935812,0.725668,0.294264); 
      rgb(805pt)=(0.937752,0.725744,0.29309); 
      rgb(806pt)=(0.939692,0.725821,0.291916); 
      rgb(807pt)=(0.941632,0.725897,0.290741); 
      rgb(808pt)=(0.943571,0.726023,0.289518); 
      rgb(809pt)=(0.945511,0.726151,0.288293); 
      rgb(810pt)=(0.947451,0.726278,0.287068); 
      rgb(811pt)=(0.949389,0.726411,0.285839); 
      rgb(812pt)=(0.951278,0.726641,0.284537); 
      rgb(813pt)=(0.953167,0.72687,0.283235); 
      rgb(814pt)=(0.955056,0.7271,0.281933); 
      rgb(815pt)=(0.956938,0.72734,0.280622); 
      rgb(816pt)=(0.958776,0.727646,0.279243); 
      rgb(817pt)=(0.960614,0.727952,0.277865); 
      rgb(818pt)=(0.962451,0.728259,0.276486); 
      rgb(819pt)=(0.964273,0.728597,0.275086); 
      rgb(820pt)=(0.966034,0.729057,0.273606); 
      rgb(821pt)=(0.967795,0.729516,0.272126); 
      rgb(822pt)=(0.969557,0.729976,0.270645); 
      rgb(823pt)=(0.971288,0.730473,0.269135); 
      rgb(824pt)=(0.972947,0.73106,0.267552); 
      rgb(825pt)=(0.974606,0.731647,0.265969); 
      rgb(826pt)=(0.976265,0.732234,0.264387); 
      rgb(827pt)=(0.977857,0.732879,0.262785); 
      rgb(828pt)=(0.979338,0.733619,0.261151); 
      rgb(829pt)=(0.980818,0.734359,0.259518); 
      rgb(830pt)=(0.982299,0.735099,0.257884); 
      rgb(831pt)=(0.983697,0.73591,0.256227); 
      rgb(832pt)=(0.984999,0.736803,0.254542); 
      rgb(833pt)=(0.986301,0.737697,0.252858); 
      rgb(834pt)=(0.987603,0.73859,0.251173); 
      rgb(835pt)=(0.988753,0.739566,0.249474); 
      rgb(836pt)=(0.989774,0.740613,0.247764); 
      rgb(837pt)=(0.990795,0.741659,0.246054); 
      rgb(838pt)=(0.991816,0.742706,0.244344); 
      rgb(839pt)=(0.992677,0.743816,0.242681); 
      rgb(840pt)=(0.993443,0.744965,0.241048); 
      rgb(841pt)=(0.994209,0.746114,0.239414); 
      rgb(842pt)=(0.994975,0.747262,0.23778); 
      rgb(843pt)=(0.995578,0.748465,0.236165); 
      rgb(844pt)=(0.996114,0.74969,0.234557); 
      rgb(845pt)=(0.99665,0.750915,0.232949); 
      rgb(846pt)=(0.997186,0.752141,0.231341); 
      rgb(847pt)=(0.997562,0.753386,0.229813); 
      rgb(848pt)=(0.997893,0.754637,0.228307); 
      rgb(849pt)=(0.998225,0.755887,0.226801); 
      rgb(850pt)=(0.998557,0.757138,0.225295); 
      rgb(851pt)=(0.998711,0.758433,0.223856); 
      rgb(852pt)=(0.998839,0.759735,0.222426); 
      rgb(853pt)=(0.998966,0.761037,0.220997); 
      rgb(854pt)=(0.999094,0.762339,0.219567); 
      rgb(855pt)=(0.999076,0.763641,0.218186); 
      rgb(856pt)=(0.99905,0.764942,0.216808); 
      rgb(857pt)=(0.999025,0.766244,0.21543); 
      rgb(858pt)=(0.998995,0.767546,0.214054); 
      rgb(859pt)=(0.998868,0.768848,0.212752); 
      rgb(860pt)=(0.99874,0.77015,0.21145); 
      rgb(861pt)=(0.998613,0.771451,0.210149); 
      rgb(862pt)=(0.998473,0.772756,0.208856); 
      rgb(863pt)=(0.998243,0.774083,0.207631); 
      rgb(864pt)=(0.998014,0.775411,0.206405); 
      rgb(865pt)=(0.997784,0.776738,0.20518); 
      rgb(866pt)=(0.997539,0.77806,0.20396); 
      rgb(867pt)=(0.997232,0.779362,0.20276); 
      rgb(868pt)=(0.996926,0.780664,0.201561); 
      rgb(869pt)=(0.99662,0.781966,0.200361); 
      rgb(870pt)=(0.996299,0.783268,0.199168); 
      rgb(871pt)=(0.995942,0.784569,0.197994); 
      rgb(872pt)=(0.995584,0.785871,0.19682); 
      rgb(873pt)=(0.995227,0.787173,0.195646); 
      rgb(874pt)=(0.994842,0.788475,0.19449); 
      rgb(875pt)=(0.994408,0.789777,0.193367); 
      rgb(876pt)=(0.993974,0.791078,0.192244); 
      rgb(877pt)=(0.99354,0.79238,0.191121); 
      rgb(878pt)=(0.993083,0.793671,0.190021); 
      rgb(879pt)=(0.992598,0.794947,0.188949); 
      rgb(880pt)=(0.992113,0.796223,0.187877); 
      rgb(881pt)=(0.991628,0.797499,0.186805); 
      rgb(882pt)=(0.99113,0.798789,0.185732); 
      rgb(883pt)=(0.990619,0.800091,0.18466); 
      rgb(884pt)=(0.990109,0.801393,0.183588); 
      rgb(885pt)=(0.989598,0.802695,0.182516); 
      rgb(886pt)=(0.989072,0.803996,0.18146); 
      rgb(887pt)=(0.988536,0.805298,0.180413); 
      rgb(888pt)=(0.988,0.8066,0.179367); 
      rgb(889pt)=(0.987464,0.807902,0.17832); 
      rgb(890pt)=(0.98691,0.809186,0.177291); 
      rgb(891pt)=(0.986349,0.810462,0.17627); 
      rgb(892pt)=(0.985787,0.811738,0.175249); 
      rgb(893pt)=(0.985226,0.813015,0.174228); 
      rgb(894pt)=(0.984644,0.814311,0.173207); 
      rgb(895pt)=(0.984057,0.815613,0.172186); 
      rgb(896pt)=(0.98347,0.816914,0.171165); 
      rgb(897pt)=(0.982883,0.818216,0.170144); 
      rgb(898pt)=(0.982296,0.819518,0.169145); 
      rgb(899pt)=(0.981709,0.82082,0.16815); 
      rgb(900pt)=(0.981122,0.822122,0.167154); 
      rgb(901pt)=(0.980535,0.823423,0.166159); 
      rgb(902pt)=(0.979947,0.824725,0.165163); 
      rgb(903pt)=(0.97936,0.826027,0.164168); 
      rgb(904pt)=(0.978773,0.827329,0.163172); 
      rgb(905pt)=(0.978186,0.828631,0.162177); 
      rgb(906pt)=(0.977599,0.829932,0.161181); 
      rgb(907pt)=(0.977012,0.831234,0.160186); 
      rgb(908pt)=(0.976425,0.832536,0.15919); 
      rgb(909pt)=(0.975838,0.833841,0.158195); 
      rgb(910pt)=(0.975251,0.835168,0.157199); 
      rgb(911pt)=(0.974664,0.836495,0.156204); 
      rgb(912pt)=(0.974077,0.837823,0.155208); 
      rgb(913pt)=(0.973489,0.83915,0.154213); 
      rgb(914pt)=(0.972902,0.840477,0.153217); 
      rgb(915pt)=(0.972315,0.841805,0.152222); 
      rgb(916pt)=(0.971728,0.843132,0.151226); 
      rgb(917pt)=(0.971155,0.844466,0.150224); 
      rgb(918pt)=(0.970619,0.845819,0.149203); 
      rgb(919pt)=(0.970083,0.847172,0.148182); 
      rgb(920pt)=(0.969547,0.848525,0.147161); 
      rgb(921pt)=(0.96902,0.849886,0.14614); 
      rgb(922pt)=(0.968509,0.851265,0.145119); 
      rgb(923pt)=(0.967999,0.852643,0.144098); 
      rgb(924pt)=(0.967488,0.854022,0.143077); 
      rgb(925pt)=(0.967,0.855411,0.142056); 
      rgb(926pt)=(0.966541,0.856815,0.141035); 
      rgb(927pt)=(0.966081,0.858219,0.140014); 
      rgb(928pt)=(0.965622,0.859623,0.138992); 
      rgb(929pt)=(0.965189,0.86104,0.137945); 
      rgb(930pt)=(0.96478,0.862469,0.136873); 
      rgb(931pt)=(0.964372,0.863899,0.135801); 
      rgb(932pt)=(0.963963,0.865328,0.134729); 
      rgb(933pt)=(0.96357,0.866773,0.133657); 
      rgb(934pt)=(0.963187,0.868228,0.132585); 
      rgb(935pt)=(0.962805,0.869683,0.131513); 
      rgb(936pt)=(0.962422,0.871138,0.130441); 
      rgb(937pt)=(0.962091,0.87261,0.129351); 
      rgb(938pt)=(0.961785,0.874091,0.128253); 
      rgb(939pt)=(0.961478,0.875571,0.127156); 
      rgb(940pt)=(0.961172,0.877052,0.126058); 
      rgb(941pt)=(0.960885,0.878571,0.124961); 
      rgb(942pt)=(0.960605,0.880103,0.123863); 
      rgb(943pt)=(0.960324,0.881634,0.122765); 
      rgb(944pt)=(0.960043,0.883166,0.121668); 
      rgb(945pt)=(0.959849,0.884719,0.120549); 
      rgb(946pt)=(0.95967,0.886276,0.119426); 
      rgb(947pt)=(0.959491,0.887833,0.118302); 
      rgb(948pt)=(0.959313,0.88939,0.117179); 
      rgb(949pt)=(0.959181,0.890995,0.116032); 
      rgb(950pt)=(0.959054,0.892603,0.114884); 
      rgb(951pt)=(0.958926,0.894211,0.113735); 
      rgb(952pt)=(0.958799,0.895819,0.112587); 
      rgb(953pt)=(0.958748,0.897453,0.111464); 
      rgb(954pt)=(0.958697,0.899086,0.110341); 
      rgb(955pt)=(0.958646,0.90072,0.109217); 
      rgb(956pt)=(0.958602,0.902359,0.108089); 
      rgb(957pt)=(0.958628,0.904043,0.106915); 
      rgb(958pt)=(0.958653,0.905728,0.105741); 
      rgb(959pt)=(0.958679,0.907413,0.104567); 
      rgb(960pt)=(0.958718,0.909102,0.103393); 
      rgb(961pt)=(0.95882,0.910812,0.102219); 
      rgb(962pt)=(0.958922,0.912522,0.101044); 
      rgb(963pt)=(0.959024,0.914232,0.0998703); 
      rgb(964pt)=(0.959153,0.915962,0.0986961); 
      rgb(965pt)=(0.959357,0.917749,0.0975219); 
      rgb(966pt)=(0.959561,0.919536,0.0963477); 
      rgb(967pt)=(0.959765,0.921323,0.0951736); 
      rgb(968pt)=(0.959996,0.923118,0.0939907); 
      rgb(969pt)=(0.960277,0.924931,0.092791); 
      rgb(970pt)=(0.960557,0.926743,0.0915913); 
      rgb(971pt)=(0.960838,0.928555,0.0903916); 
      rgb(972pt)=(0.961151,0.930378,0.0891919); 
      rgb(973pt)=(0.961509,0.932216,0.0879922); 
      rgb(974pt)=(0.961866,0.934054,0.0867925); 
      rgb(975pt)=(0.962223,0.935892,0.0855928); 
      rgb(976pt)=(0.96262,0.937768,0.0843802); 
      rgb(977pt)=(0.963053,0.939683,0.083155); 
      rgb(978pt)=(0.963487,0.941597,0.0819297); 
      rgb(979pt)=(0.963921,0.943512,0.0807045); 
      rgb(980pt)=(0.964415,0.945426,0.0794643); 
      rgb(981pt)=(0.964951,0.947341,0.0782135); 
      rgb(982pt)=(0.965487,0.949255,0.0769628); 
      rgb(983pt)=(0.966023,0.951169,0.075712); 
      rgb(984pt)=(0.966594,0.953118,0.0744441); 
      rgb(985pt)=(0.967181,0.955083,0.0731679); 
      rgb(986pt)=(0.967768,0.957049,0.0718916); 
      rgb(987pt)=(0.968355,0.959014,0.0706153); 
      rgb(988pt)=(0.96898,0.960999,0.0693006); 
      rgb(989pt)=(0.969619,0.96299,0.0679733); 
      rgb(990pt)=(0.970257,0.964981,0.0666459); 
      rgb(991pt)=(0.970895,0.966972,0.0653186); 
      rgb(992pt)=(0.971554,0.968984,0.0639486); 
      rgb(993pt)=(0.972218,0.971001,0.0625703); 
      rgb(994pt)=(0.972882,0.973017,0.0611919);  
      rgb(995pt)=(0.973545,0.975034,0.0598135); 
      rgb(996pt)=(0.974232,0.97705,0.058318); 
      rgb(997pt)=(0.974922,0.979067,0.056812); 
      rgb(998pt)=(0.975611,0.981083,0.055306); 
      rgb(999pt)=(0.9763,0.9831,0.0538) 
    }
}

\usetikzlibrary{arrows,shapes,positioning}
\usetikzlibrary{decorations.markings}

\newcounter{tikzsubfigcounter}[figure]
\renewcommand{\thetikzsubfigcounter}{\thesection.\the\numexpr\value{figure}+1\relax\alph{tikzsubfigcounter}}
\newcommand*{\tikztitle}[1]{ %
  \refstepcounter{tikzsubfigcounter}
  (\alph{tikzsubfigcounter})\space\space #1
}
\newcounter{tikzsubfigcounterinvisible}[figure]
\renewcommand{\thetikzsubfigcounterinvisible}{\thesection.\the\numexpr\value{figure}+1\relax\alph{tikzsubfigcounterinvisible}}

\newcommand*{\inlinecode}[1]{\texttt{#1}}

\begin{document}

\begin{abstract}
  In this contribution we investigate in mathematical modeling and efficient simulation of biological
  cells
  with a particular emphasis on effective modeling of structural properties that originate from
  active forces
  generated from polymerization and depolymerization of cytoskeletal components.
  In detail, we propose a nonlinear continuum approach to model microtubule-based forces
  which have recently been established as central components of cell mechanics during early fruit fly
  wing development.
  The model is discretized in space using the finite-element method. Although the
  individual equations are decoupled by a semi-implicit time discretization, the discrete model is
  still computationally demanding.
  In addition, the parameters needed for the effective model equations are not easily available and
  have to be estimated or
  determined by repeatedly solving the model and fitting the results to measurements.
  This drastically increases the computational cost. Reduced basis methods have been used
  successfully to speed up
  such repeated solves, often by several orders of magnitude.  However, for the complex nonlinear
  models regarded
  here, the application of these model order reduction methods is not always straight-forward and
  comes
  with its own set of challenges. In particular, subspace construction using the Proper
  Orthogonal Decomposition (POD) becomes prohibitively expensive for reasonably fine grids.
  We thus propose to combine the Hierarchical Approximate POD, which is a general, easy-to-implement
  approach to compute
  an approximate POD, with an Empirical Interpolation Method to efficiently generate a fast to
  evaluate reduced order model.
  Numerical experiments are given to demonstrate the applicability and efficiency of the proposed
  modeling and simulation approach.
\end{abstract}

\maketitle

\noindent


\def\tikzpath{Images/}

\section{Introduction}
Cells are often called the ``smallest unit of life'', as all living creatures are composed of these
small entities. Notwithstanding their size,
cells are incredibly complex biological machines. This is reflected in the fact that even after
many years of research we are nowhere near having a
full understanding of all the processes that run even in a single cell. The picture gets even more
complicated when regarding the tissue level,
where cell-cell and
cell-extracellular matrix interactions have to taken into account.
Coarse-grained models try to remove some of the complexity by looking at the cell from a
macroscopic view point, averaging over the microscopic
details to focus on the essential features. Although these models are much simpler than full
microscopic models they are often still very demanding
computationally. In addition, the parameters needed for the model equations are often not available
and have to be estimated or
determined by repeatedly solving the model and comparing the results to measurements, which
drastically increases the computational cost.
In this contribution we focus on effective mathematical modeling of cytoskeletal proteins as the
main source of forces that determine the shape of biological cells and tissues.
While actin-based forces are thoroughly investigated, forces originating from
other cytoskeletal components came into focus only recently.
More specifically, we propose to use a nonlinear continuum approach \cite{Kruseetal05,VERDIER2009790,Marth2015,FAN2015559} to model
microtubule-based forces
which have recently been established as central components of cell mechanics during early fruit fly
wing
development \cite{Singh2018}.

While there is a broad discussion on continuous versus discrete mathematical modeling approaches
for cell
metamorphosis \cite{BRODLAND201562,Ingber2014,MR3587800,Banavaretal2021}, we concentrate here on the development of efficient numerical
discretization and model order reduction
methods for the parameterized system of coupled non-linear partial differential equations resulting
from this continuum approach.
To this end, the mathematical model is discretized in space using the finite-element method (FEM).
Although the
individual equations are decoupled by a semi-implicit time discretization, the discrete model is
still very demanding
computationally.
This is particularly a problem since computational studies to identify effective parameters of the
model require many evaluations of
the finite element system with different parameter settings and thus involve a large amount of time
and experimental effort.
Reduced basis methods have been used successfully to speed up such repeated solves, often
by several orders of magnitude.
For an overview on recent developments in model order reduction,  in general, and the reduced basis
method, in particular,
we refer to the monographs and collections
\cite{RB_Rozza,Quarteroni2016, MR3672144, MR3701994}.

For time-dependent problems, the POD-Greedy method \cite{Pod_Greedy} has been established as a
standard method for
reduced basis constructions. However, for an efficient implementation of the POD-Greedy method,
rigorous and cheap to evaluate a posteriori error estimates have to be available.
As this is not the case for the non-linear cell model at hand, we propose to speed-up a subspace
construction
based on the Proper Orthogonal Decomposition (POD) by employing a suitable variant of the
Hierarchical Approximate POD (HAPOD) \cite{Himpe2018}. This procedure enormously reduces the
so called
offline-complexity of the model order reduction approach. In order to also achieve a good
online-performance
of the method for our non-linear system, a variant of the empirical interpolation method (EIM)
\cite{eim,Chaturantabut2010,eim_operator}
is used to efficiently compute the reduced model.

The rest of this paper is organized as follows. In \cref{sec:biological_background} we will shortly review the
biological background. In
\cref{sec:cell_model} we will specify the mathematical cell model and its discretization. In
\cref{sec:modelreduction} we will present and discuss
the model reduction approach. Finally, we will use the reduced model for some numerical studies in
\cref{sec:numericalexperiments}.

\section{Biological background}\label{sec:biological_background}
A cell's shape and structural properties are largely controlled by the cytoskeleton, an
interconnected network of filaments and associated
proteins~\cite{Fletcher2010}. The filamentous proteins are often divided into three classes:
intermediate filaments, microtubules, and actin
filaments. Intermediate filaments are the most durable and long-lived of these three and their main
function is to withstand tensile mechanical
stress~\cite{Fletcher2010,Alberts2014}. Actin filaments and microtubules are more dynamic. They are build
from small monomers that can assemble into
large rod-like polymers. However, they can also rapidly disassemble again. While this applies to
both actin and
microtubules, the latter are particularly known for their
\emph{dynamic
  instability}~\cite{Alberts2014}, constantly assembling and then disassembling
again if not stabilized by special cap proteins.
Unlike intermediate filaments, actin filaments and microtubules are polar, i.e., they are
asymmetric on molecular level. Specialized \emph{motor
  proteins} use this polarity to actively
transport cell components in a directed manner along these
filaments. Motor proteins can also cross-link
filaments and exert active forces by sliding filaments against each other. Most prominently, the
contraction of muscle cells is achieved by a large
number of parallel actin filaments that are pulled together by myosin
motors~\cite{Alberts2014,Squire2019}. It is also long known that actin and myosin form the
contractile ring during cytokinesis (the process that divides a cell in two) and are involved in
cell motility,
mechanosensing and mechanotransduction (e.g.~\cite{Rappaport1996,Robinson2000,Pelham2002,Pollard2009,Meitinger2016,Korn1978,Ladoux2016,Pandya2017}).
Moreover, they form the actomyosin cortex, a thin actin-myosin network attached to the cell
membrane that plays an essential role in shaping the cell
and controlling mechanical properties of the cell surface~\cite{Salbreux2012}

In contrast, microtubules were traditionally not viewed in the role of a force generator, partly
because early studies underestimated the
amount of force a microtubule can exert \emph{in vivo} by polymerization at its plus end by
about an order of magnitude~\cite{Matis2020}. Only
recently microtubule-based forces gained broader attention~\cite{Bouchet2017,Singh2018,DAngelo2019,Matis2020}. We will here
focus on a
recent study showing that microtubules play a major role in the mechanics of tissue morphogenesis
during fruit fly (drosophila) wing
development~\cite{Singh2018}.

\section{Modelling}\label{sec:cell_model}
\subsection{Notation}
Vectors are written in bold, lower-case notation and matrices in bold, upper-case notation.

Let $\MA = {(\Va_1, \ldots, \Va_d)}^T \in \R^{d \times d}$ be a matrix with rows $\Va_i$ and let
$\Vu\colon \R^d \to \R^d$.
Then the divergence of $\MA$ is defined as
\begin{equation}
  \Vdiverg \MA \coloneqq
  \begin{pmatrix}
    \diverg \Va_1 \\
    \vdots        \\
    \diverg \Va_d
  \end{pmatrix}
\end{equation}
and the gradient of $\Vu$ is defined as
\begin{equation}
  \Vgradient \Vu :=
  \begin{pmatrix}
    \frac{\partial u_1}{\partial x_1} & \cdots & \frac{\partial u_1}{\partial x_d} \\
    \vdots                            & \vdots & \vdots                            \\
    \frac{\partial u_d}{\partial x_1} & \cdots & \frac{\partial u_d}{\partial x_d} \\
  \end{pmatrix}
\end{equation}
We define
\begin{equation}
  (\Vu \cdot \nabla) \Vu = \sum_{i=1}^d u_i \frac{\partial}{\partial x_i} \Vu.
\end{equation}
The Laplace operator is defined as
\begin{equation}
  \Vlaplace \Vu = \sum \limits_{i=1}^{d} \frac{\partial^2}{\partial x_i^2} \Vu = \Vdiverg \Vgradient \Vu.
\end{equation}

Let $(\cdot, \cdot)$ denote the $L^2$ inner product, i.e.
\begin{equation}
  (u,v) := \int \limits_\Omega u(\Vx) v(\Vx) \mathrm{d}\Vx.
\end{equation}
For vector-valued functions we define
\begin{equation}
  (\Vu,\Vv) := \int \limits _\Omega \Vu(\Vx) \cdot \Vv(\Vx)
  \mathrm{d}\Vx.
\end{equation}
For matrix-valued functions $\MU = {(\Vu_1, \ldots, \Vu_d)}^T$, $\MV = {(\Vv_1, \ldots, \Vv_d)}^T$ we define
\begin{equation}
  (\MU,\MV) := \sum_{i=1}^d (\Vu_i, \Vv_i).
\end{equation}

\subsection{Derivation of model equations}
We will model the cell as an active polar gel surrounded by a membrane that separates it from the
surrounding extracellular
fluid~\cite{Kruse2005,Tjhung2012,Marth2015}. See \cite{Marth2015} and references therein for a detailed
derivation. Here, we
only shortly recapitulate the important features of the model.
The model uses a diffuse interface description of the cell, i.e., the cell is modelled by a phase
field parameter \(\pfield\) that takes on the value
\(1\) inside the cell and \(-1\) outside with a smooth transition in the membrane region. The width
of the transition region can be controlled by the
model parameter
\(\pfieldepsilon\).
The average orientation of the microtubules inside the cell is tracked by the vector-valued
orientation field \(\ofield\). For the fluid, we track the velocity \(\stokesvelocity\) and
pressure \(\stokespressure\). For simplicity, we will
assume equal
density for the cytoplasm and extracellular fluid which is justified as both mainly consist of
water~\cite{Alberts2014}.

The model equations are obtained by first stating the free energy of the system and then
assuming that the system evolves according to a gradient descent of this free
energy~\cite{DeGroot2013,Marth2015}. The free energy
\begin{equation}\label{eq:freeenergy}
  \freeenergy(\pfield, \ofield, \stokesvelocity) = \kineticenergy(\stokesvelocity) + \surfaceenergy(\pfield) + \filamentenergy(\pfield,\ofield)
\end{equation}
is composed of the kinetic energy of the fluid, the membrane (or surface) energy \(\surfaceenergy\)
and the energy of the filament network
\(\filamentenergy\).
In the following, we will state the energies and equations in non-dimensionalized form only,
see~\cite{Marth2015} for details on the
non-dimensionalization.
The kinetic energy then simply is
\begin{equation}\label{eq:kineticenergy}
  \kineticenergy(\stokesvelocity) = \frac{\reynoldsnumber}{2} \bulkintegral{\stokesvelocity^2}
\end{equation}
where \(\reynoldsnumber\) is the well-known Reynolds number giving the ratio of inertial to viscous
forces within the fluid.
The surface energy reads
\begin{equation}\label{eq:surfaceenergy}
  \surfaceenergy(\pfield) = \frac{1}{\capillarynumber}\bulkintegral{\frac{\pfieldepsilon}{2} \abs{\gradient \pfield}^2 +
    \frac{1}{\pfieldepsilon}\doublewellpotential(\pfield)}
  + \frac{1}{\bendingcapillarynumber}\bulkintegral{\frac{1}{2\pfieldepsilon}{\left( \pfieldepsilon \laplace \pfield - \frac{1}{\pfieldepsilon}
      \doublewellpotentialprime(\pfield)\right)}^2}.
\end{equation}
The first integral in~\eqref{eq:surfaceenergy} corresponds to a classic
Cahn-Hilliard~\cite{Cahn1958,Cahn1959} model
where \(\doublewellpotential(\pfield) = \frac{1}{4}{\left( \pfield^2 - 1 \right)}^2\) is a double well potential with two
minima at \(\pm 1\) giving
the two pure ``phases'' (interior and exterior of the cell) a lower free energy than the interface
region where \(\pfield \in (-1,1)\).
The second integral describes membrane bending energy by a Helfrich-type
model~\cite{Helfrich1973,Hausser2013}. Accordingly, the capillary
number \(\capillarynumber\) and the bending capillary number \(\bendingcapillarynumber\) describe
the ratio of viscous drag forces to surface tension
forces and forces resisting bending, respectively.

The energy of the filament network is modelled as
\begin{equation}\label{eq:filamentenergy}
  \filamentenergy(\pfield, \ofield) = \frac{1}{\polaritynumber}\bulkintegral{\frac{1}{2}\doubledotproduct{\Vgradient\ofield}{\Vgradient\ofield} +
    \frac{\ofielddoublewellparam}{4}\abs{\ofield}^2\left( -2\pfield + \abs{\ofield}^2 \right)}.
\end{equation}
Here, the first term penalizes distortion of the filaments and promotes alignment. The polarity
number \(\polaritynumber\) again characterizes the
ratio between viscous forces and forces due to distortion of the filaments. The second term
controlled by \(\ofielddoublewellparam\) (weakly)
restricts \(\ofield\) to the interior of the cell and (weakly) enforces \(\abs{\ofield} =
1\) inside the cell: If \(\pfield < 0\), the term has a
single minimum at zero, otherwise it is again a double well potential with minima at \(\pm 1\).

The phase field is advected with the flow and dissipates the energy by a gradient flow on the
zero-average subspace of
\(\Hminusone\)~\cite{Fife2000,Lee2014}, giving the phase field equations in
the domain \(\domain\)
\begin{subequations}\label{eq:pfield_strong}
  \begin{align}
    \dt \pfield + \diverg(\stokesvelocity \pfield) & = \pfieldmobility \laplace \pfieldphinat    \label{eq:pfield_strong1}                     \\
    \pfieldphinat                                  & = \variationalderiv{\freeenergy}{\pfield}                                                 \\
    \pfieldmu                                      & = \pfieldepsilon\laplace\pfield - \frac1\pfieldepsilon \doublewellpotentialprime(\pfield)
  \end{align}
  with
  \begin{align}
    \begin{split}
      \variationalderiv{\freeenergy}{\pfield} & = \variationalderiv{\filamentenergy}{\pfield} + \variationalderiv{\surfaceenergy}{\pfield} \\
      & =-\frac{\ofielddoublewellparam}{2\polaritynumber} \abs{\ofield}^2
      + \frac{1}{\bendingcapillarynumber} \left( \laplace \pfieldmu
      - \frac{1}{\pfieldepsilon^2}\doublewellpotentialtwoprime(\pfield)\pfieldmu \right) - \frac{1}{\capillarynumber}\pfieldmu.
    \end{split}
  \end{align}
  and initial and boundary conditions
  \begin{align}
    \pfield(0,\spatialvar)                      & = \pfieldinitial(\spatialvar)
                                                & \text { for all } \spatialvar \in \domain \label{eq:pfieldinitial}
    \\
    \normalderiv{\pfield}(\timevar,\spatialvar) & = \normalderiv{\pfieldphinat}(\timevar,\spatialvar)  =
    \normalderiv{\pfieldmu}(\timevar,\spatialvar)  = 0
                                                & \text{ for all } \timevar \in [0, \tf], \spatialvar \in \domainboundary
    \label{eq:pfield_boundary_conditions}
  \end{align}
  Here, \(\pfieldphinat\) and \(\pfieldmu\) are helper variables introduced to write the sixth-order
  equation as a system of second-order equations.
  The phase field mobility coefficient \(\pfieldmobility\) is a model parameter that regulates the
  rate of entropy dissipation.
\end{subequations}

The orientation field equations are lend from liquid crystal theory:
\begin{subequations}\label{eq:ofield_strong}
  \begin{align}
    \dt \ofield & + \left(\Vgradient \ofield\right) \stokesvelocity
    + \left(\vorticitytensor(\stokesvelocity) - \shapefactor \deformationtensor(\stokesvelocity)\right) \ofield
    = -\frac{1}{\fracrotationaldynamicviscosity} \ofieldPnat                                                                           \\
    \ofieldPnat & = \variationalderiv{\freeenergy}{\ofield} = \frac{1}{\polaritynumber} \left( -\ofielddoublewellparam \pfield \ofield
    + \ofielddoublewellparam (\ofield \cdot \ofield) \ofield
    - \Vlaplace \ofield \right)
  \end{align}
  with vorticity tensor
  \begin{align}
    \vorticitytensor(\stokesvelocity) & = \frac{1}{2}(\Vgradient \stokesvelocity^T - \Vgradient \stokesvelocity),
  \end{align}
  deformation tensor
  \begin{align}
    \deformationtensor(\stokesvelocity) & = \frac{1}{2}(\Vgradient \stokesvelocity^T + \Vgradient \stokesvelocity)
  \end{align}
  and initial and boundary conditions
  \begin{align}
    \ofield(0,\spatialvar)                         & = \ofieldinitial(\spatialvar)                                            & \text { for all }
    \spatialvar \in \domain \label{eq:ofieldinitial}
    \\
    \normalderivvec{\ofield}(\timevar,\spatialvar) & = \Vzero
                                                   & \text{ for all } \timevar \in [0, \tf], \spatialvar \in \domainboundary.
    \label{eq:ofield_boundary_conditions}
  \end{align}
\end{subequations}
The filaments are advected with the flow and energy is dissipated by a \(\Ltwo\) gradient flow with
scaling factor
\(\frac{1}{\fracrotationaldynamicviscosity}\). The shape factor \(\shapefactor\) influences
alignment of the filaments with the flow and is positive for the rod-like
microtubules~\cite{Marchetti2013,Marth2015}.

Finally, evolution of the flow is modelled by the Navier-Stokes equations
\begin{subequations}\label{eq:navier_stokes_strong}
  \begin{align}
    \reynoldsnumber\left(\dt \stokesvelocity + (\Vgradient \stokesvelocity)\stokesvelocity\right) + \gradient \stokespressure
                            & = \Vdiverg\left(\viscousstress + \activestress + \distortionstress + \ericksenstress\right) \\
    \diverg \stokesvelocity & = 0
  \end{align}
  where
  \begin{align}
    \viscousstress(\stokesvelocity)         & = \deformationtensor(\stokesvelocity) = \frac{1}{2} \left( {\Vgradient
      \stokesvelocity}^T +
    \Vgradient \stokesvelocity \right),\label{eq:viscousstress}
    \\
    \activestress(\pfield,\ofield)          & = \frac{1}{\activeforcenumber} \pfieldcellindicator \ofield \otimes
    \ofield, \ \
    \pfieldcellindicator =
    \frac12(\pfield + 1)\label{eq:activestress}
    \\
    \distortionstress(\ofield, \ofieldPnat) & = \frac12(\ofieldPnat \otimes \ofield - \ofield \otimes \ofieldPnat)
    + \frac\shapefactor{2} (\ofieldPnat \otimes \ofield + \ofield \otimes \ofieldPnat),\label{eq:distortionstress}
    \\
    \begin{split}\label{eq:ericksenstress}
      \ericksenstress(\pfield, \pfieldphinat, \ofield, \ofieldPnat)
      &=
      \frac{1}{\bendingcapillarynumber}
      (\gradient\pfield\otimes\gradient\pfieldmu -
      \pfieldmu \Vgradient\gradient\pfield)
      - \frac{1}{\capillarynumber}
      (\epsilon\gradient\pfield\otimes\gradient\pfield) \\
      & - \frac{1}{\polaritynumber} \left(\Vgradient\ofield^T\cdot\Vgradient\ofield\right)
    \end{split}
  \end{align}
  are the stresses stemming from viscosity of the fluid, active extension or contraction of the
  filaments, distortion of the filaments and additional
  stress from the energy-minimizing behavior of both filaments and membrane,
  respectively~\cite{Marth2015}. The active force number
  \(\activeforcenumber\) describes the ratio of viscous forces to the active forces and is positive
  for contractile and negative for extensile
  stress. The term \(\pfieldcellindicator\) serves as an indicator function for \(\pfield > 0\) and
  thus restricts the active stress to the
  cytoplasm.
  However, as we already weakly enforce \(\ofield = 0\) outside of the cell (see above), removing
  this factor probably would not change the dynamics
  of
  the model significantly.

  By redefining the pressure, we can simplify the divergence of the last stress term to obtain (see
  \cite[Section 6.2.4 and Remark 3.2]{Marth2016})
  \begin{equation}\label{eq:ericksenstressdiv}
    \Vdiverg(\ericksenstress(\pfield, \pfieldphinat, \ofield, \ofieldPnat)) = \pfieldphinat \gradient \pfield + {(\Vgradient \ofield)}^T
    \ofieldPnat,
  \end{equation}

  The equations are complemented with initial and boundary conditions
  \begin{align}
    \stokesvelocity(0, \spatialvar) = \stokesvelocityinitial(\spatialvar) & \text { for all } \spatialvar \in \domain              \label{eq:stokesvelocityinitial} \\
    \stokesvelocity(\timevar, \spatialvar) = \Vzero
                                                                          & \text{ for all } \timevar \in [0, \tf], \spatialvar \in
    \domainboundary.\label{eq:stokes_boundary_conditions}
  \end{align}
\end{subequations}
As cell motility generally occurs in a regime with very low Reynolds number, in the following,
we will restrict
ourselves to the Stokes approximation obtained by setting \(\reynoldsnumber = 0\)
in~\eqref{eq:navier_stokes_strong}.

The model~\eqref{eq:pfield_strong}-\eqref{eq:navier_stokes_strong} conserves the cell volume, i.e.,
\begin{equation}
  \dt \bulkintegral{\pfield} = 0,
\end{equation}
and is thermodynamically consistent, i.e., the free energy~\eqref{eq:freeenergy} decreases over
time
\begin{equation}
  \dt \freeenergy \leq 0.
\end{equation}

\subsection{Weak formulations}\label{sec:cellmodeldiscretization}
We will use a finite element discretization to numerically solve the model equations. To that end,
we first derive weak formulations for each of the
systems of equations~\eqref{eq:pfield_strong}-\eqref{eq:navier_stokes_strong}.

\subsubsection{Phase field}
Multiplying~\eqref{eq:pfield_strong} by test functions \((\firstpfieldtest, \secondpfieldtest,
\thirdpfieldtest) \in {(\Hone(\domain))}^3\)
and integrating gives
\begin{align*}
  \dt \bulkintegral{\pfield \firstpfieldtest}    & + \boundaryint{\pfield \firstpfieldtest \stokesvelocity \cdot \outernormal}
  - \bulkintegral{\pfield \stokesvelocity \cdot \gradient \firstpfieldtest}
  = \pfieldmobility \left( \boundaryint \firstpfieldtest \gradient \pfieldphinat \cdot \outernormal
  - \bulkintegral{\gradient \pfieldphinat \cdot \gradient \firstpfieldtest} \right)
  \\
  \bulkintegral{\pfieldphinat \secondpfieldtest} & = \bulkintegral{
    \left( -\frac{\ofielddoublewellparam}{2 \polaritynumber} \ofield \cdot \ofield
    - \left(\frac1{\capillarynumber}  + \frac1{\bendingcapillarynumber \cdot \pfieldepsilon^2}
    \doublewellpotentialtwoprime(\pfield)\right)\pfieldmu\right)
    \secondpfieldtest }
  \\
                                                 & + \frac1{\bendingcapillarynumber} \left( \boundaryintegral{\secondpfieldtest \gradient \pfieldmu
    \cdot
    \outernormal}
  - \bulkintegral{\gradient \pfieldmu \cdot \gradient \secondpfieldtest} \right)
  \\
  \bulkintegral{\pfieldmu \thirdpfieldtest}      & = \pfieldepsilon \left( \boundaryintegral{\thirdpfieldtest \gradient \pfield
    \cdot \outernormal}
  - \bulkintegral{\gradient \pfield \cdot \gradient \thirdpfieldtest} \right) - \bulkintegral{\frac{1}{\pfieldepsilon}
    \doublewellpotentialprime(\pfield)
    \thirdpfieldtest}.
\end{align*}
Using the boundary conditions~\eqref{eq:pfield_boundary_conditions} and~\eqref{eq:stokes_boundary_conditions}, the boundary
integrals vanish and we
arrive at the weak formulation
\begin{subequations}\label{eq:pfield_weak}
  \begin{align}
    \dt \bulkintegral{\pfield \firstpfieldtest} -
    \bulkintegral{(\pfield \stokesvelocity - \pfieldmobility \gradient \pfieldphinat) \cdot \gradient
    \firstpfieldtest}           & = 0 \label{eq:pfield_weak1} \\
    \begin{split}
      \bulkintegral{
        \pfieldphinat \secondpfieldtest  +
        \frac{1}{\bendingcapillarynumber} \gradient \pfieldmu \cdot \gradient \secondpfieldtest +
        \left(\frac1{\capillarynumber} + \frac1{\bendingcapillarynumber \cdot \pfieldepsilon^2}
        \doublewellpotentialtwoprime(\pfield)\right)\pfieldmu \secondpfieldtest
      }
      &= - \bulkintegral{\left(\frac{\ofielddoublewellparam}{2 \polaritynumber} \ofield \cdot \ofield \right)\secondpfieldtest}
    \end{split}
    \\
    \bulkintegral{\pfieldmu \thirdpfieldtest + \frac1\pfieldepsilon \doublewellpotentialprime(\pfield)
      \thirdpfieldtest + \pfieldepsilon
      \gradient \pfield \cdot
    \gradient \thirdpfieldtest} & = 0
  \end{align}
\end{subequations}
for \((\pfield, \pfieldphinat, \pfieldmu) \in {(\Hone(\domain))}^3\).

\subsubsection{Orientation field}
Multiplying~\eqref{eq:ofield_strong} by test functions \((\firstofieldtest, \secondofieldtest)
\in {(\Hone(\domain))}^{\dimension} \times
{(\Ltwo(\domain))}^{\dimension}\)
and integrating gives
\begin{align*}
  0 & = \dt \bulkintegral{\ofield \cdot \secondofieldtest}
  + \bulkintegral{\left(\left(\Vgradient \ofield\right) \stokesvelocity + \Big(\vorticitytensor(\stokesvelocity) - \shapefactor
    \deformationtensor(\stokesvelocity)\Big)\ofield
  + \frac1{\fracrotationaldynamicviscosity}\ofieldPnat \right) \cdot \secondofieldtest}, \\
  \begin{split}
    0 &= \bulkintegral{\left( \ofieldPnat +\frac{\ofielddoublewellparam}{\polaritynumber} \pfield \ofield -
      \frac{\ofielddoublewellparam}{\polaritynumber} (\ofield \cdot \ofield)\ofield \right) \cdot \firstofieldtest -
      \frac{1}{\polaritynumber}\doubledotproduct{\Vgradient \ofield}{\Vgradient \firstofieldtest}} \\
    &+\frac{1}{\polaritynumber} \boundaryintegral{(\Vgradient \ofield \cdot \outernormal) \cdot \firstofieldtest}
  \end{split}
\end{align*}
Using the boundary conditions~\eqref{eq:ofield_boundary_conditions} the last term vanishes and we arrive at the
weak formulation
for the orientation field variables \((\ofield, \ofieldPnat) \in {(\Hone(\domain))}^{\dimension}
\times
{(\Ltwo(\domain))}^{\dimension}\)
\begin{subequations}\label{eq:ofield_weak}
  \begin{align} \dt \bulkintegral{\ofield \cdot \secondofieldtest} +
    \bulkintegral{\left(\left(\Vgradient \ofield\right) \stokesvelocity + \Big(\vorticitytensor(\stokesvelocity) - \shapefactor
      \deformationtensor(\stokesvelocity)\Big)\ofield
    + \frac1{\fracrotationaldynamicviscosity}\ofieldPnat \right) \cdot \secondofieldtest} & = 0,
    \\
    \bulkintegral{\left( \ofieldPnat +\frac{\ofielddoublewellparam}{\polaritynumber} \pfield \ofield -
      \frac{\ofielddoublewellparam}{\polaritynumber} (\ofield \cdot \ofield)\ofield
      \right) \cdot \firstofieldtest - \frac{1}{\polaritynumber}\doubledotproduct{\Vgradient \ofield}{\Vgradient \firstofieldtest}}
                                                                                          & = 0
  \end{align}
\end{subequations}

\subsubsection{Stokes}
We use a standard weak formulation~\cite{Quarteroni2017} of the Stokes
equations~\eqref{eq:navier_stokes_strong}: Find \((\stokesvelocity,
\stokespressure) \in {\Honezero(\domain)}^{\dimension} \times \Ltwozero(\domain)\) such that
\begin{subequations}\label{eq:stokes_weak}
  \begin{align}
    \begin{split}
      \bulkintegral{\stokespressure \cdot \diverg{\firststokestest} + \frac{1}{2}\doubledotproduct{\Vgradient \stokesvelocity}{\Vgradient \firststokestest}}
      & = -\bulkintegral{\doubledotproduct{\left(\activestress(\pfield, \ofield) + \distortionstress(\ofield,
          \ofieldPnat) + \ericksenstress(\pfield, \pfieldphinat, \ofield, \ofieldPnat)\right)}{\Vgradient \firststokestest}} \label{eq:stokes_weak_a}
    \end{split} \\
    \bulkintegral{\secondstokestest \diverg{\stokesvelocity}} & = 0
  \end{align}
\end{subequations}
for all \((\firststokestest,\secondstokestest) \in {\Honezero(\domain)}^{\dimension} \times
\Ltwozero(\domain)\).
Note that
\begin{align*}
  \Vdiverg(\viscousstress)
   & = \Vdiverg\left(\frac{1}{2}\left(\Vgradient \stokesvelocity + {(\Vgradient \stokesvelocity)}^T\right)\right)
  = \frac{1}{2}\left(\Vlaplace \stokesvelocity + \Vdiverg \left({(\Vgradient \stokesvelocity)}^T\right)\right)                                  \\
   & = \frac{1}{2}\left(\Vlaplace \stokesvelocity + \gradient\left(\diverg \stokesvelocity\right)\right) = \frac{1}{2}\Vlaplace \stokesvelocity \\
\end{align*}
due to the incompressibility constraint on the fluid, which after partial integration gives the
second term on the left-hand side
of~\eqref{eq:stokes_weak_a}. Moreover, we used the fact that
\begin{align*}
  \bulkintegral{\Vdiverg(\MA)\cdot\firststokestest}
   & = \bulkintegral{\sum_{\dimindex=0}^{\dimension-1} \diverg(\Va_\dimindex) v_\dimindex}                   \\
   & = \boundaryintegral{\sum_{\dimindex=0}^{\dimension-1} v_{\dimindex} \Va_{\dimindex} \cdot \outernormal}
  - \bulkintegral{\sum_{\dimindex=0}^{\dimension-1} \Va_{\dimindex} \cdot \gradient v_{\dimindex}}           \\
   & = - \bulkintegral{\doubledotproduct{\MA}{\Vgradient\firststokestest}}
\end{align*}
for any matrix \(\MA\) (here, the boundary integrals vanish since \(\firststokestest \in
\Honezero(\domain)\)).

For the Ericksen stress, instead of performing the partial integration, we will usually use the
reformulation \eqref{eq:ericksenstressdiv}
and thus directly use the term
\begin{equation}
  \bulkintegral{\Vdiverg(\ericksenstress(\pfield, \pfieldphinat, \ofield, \ofieldPnat)) \cdot \firststokestest}
\end{equation}
on the right-hand side of \eqref{eq:stokes_weak}.

\subsection{Finite element discretization}\label{sec:finiteelementdiscretization}
For the space discretization, we use a standard Galerkin finite element discretization. For
simplicity,
we use the same approximation order \(\pfieldpolynomialorder\) for the phase field variables,
the orientation field variables
and the Stokes pressure. We use Taylor-Hood elements to discretize the Stokes equations, i.e., the
approximation order for the Stokes velocity is chosen to be \(\pfieldpolynomialorder+1\).

Let \(\fespace \subset \Hone(\domain) \subset \Ltwo(\domain)\) be the continuous
finite element
approximation spaces with polynomial order \(\pfieldpolynomialorder\) and grid width
\(\gridwidth\), and let \(\fespacezero \subset \fespace \cap
\Honezero \) be the space of functions in \(\fespace\) that are zero on the Dirichlet boundary.
Let \(\pfieldfebasisfunc_1, \ldots, \pfieldfebasisfunc_{\pfieldfedimension}\) be a basis of \(\fespace\) and
\(\ofieldfebasisfunc_1, \ldots, \ofieldfebasisfunc_{\ofieldfedimension}\) be a basis of
\({(\fespace)}^{\dimension}\).
Moreover, let \(\stokesvelocityfebasisfunc_1, \ldots, \stokesvelocityfebasisfunc_{\stokesvelocityfedimension}\) be a basis of
\({(\fespacezero[\stokespressurepolynomialorder+1])}^{\dimension}\).
We then write the finite element approximations of our variables as
\(\pfieldfe = \sum_{\feindex=1}^{\pfieldfedimension} \pfieldfecoeff_\feindex \pfieldfebasisfunc_\feindex\),
\(\pfieldphinatfe = \sum_{\feindex=1}^{\pfieldfedimension} \pfieldphinatfecoeff_\feindex
\pfieldfebasisfunc_\feindex\),
\(\pfieldmufe = \sum_{\feindex=1}^{\pfieldfedimension} \pfieldmufecoeff_\feindex \pfieldfebasisfunc_\feindex\),
\(\ofieldfe = \sum_{\feindex=1}^{\ofieldfedimension} \ofieldfecoeff_\feindex \ofieldfebasisfunc_\feindex\),
\(\ofieldPnatfe = \sum_{\feindex=1}^{\ofieldfedimension} \ofieldPnatfecoeff_\feindex
\ofieldfebasisfunc_\feindex\),
\(\stokesvelocityfe = \sum_{\feindex=1}^{\stokesvelocityfedimension} \stokesvelocityfecoeff_\feindex
\ofieldfebasisfunc_\feindex\)
and
\(\stokespressurefe = \sum_{\feindex=1}^{\pfieldfedimension} \stokespressurefecoeff_\feindex
\pfieldfebasisfunc_\feindex\)
with corresponding coefficient vectors \(\pfieldfecoeffvec, \pfieldphinatfecoeffvec,
\pfieldmufecoeffvec, \stokespressurefecoeffvec \in
\Reals[\pfieldfedimension]\),
\(\ofieldfecoeffvec, \ofieldPnatfecoeffvec \in \Reals[\ofieldfedimension]\) and
\(\stokesvelocityfecoeffvec \in \Reals[\stokesvelocityfedimension]\).

\subsubsection{Phase field}
Inserting these basis representations in~\eqref{eq:pfield_weak} gives the system
\begin{subequations}\label{eq:pfield_weak_algebraic_semidiscrete}
  \begin{align}
    \dt \pfieldmassmatrix \pfieldfecoeffvec - \pfieldBmatrix(\stokesvelocityfe) \pfieldfecoeffvec + \pfieldmobility \pfieldellipticmatrix
    \pfieldphinatfecoeffvec               & = \Vzero                                                                       \\
    \frac{1}{\bendingcapillarynumber \cdot \pfieldepsilon^2} \pfieldfenonlinf(\pfieldfe, \pfieldmufe) + \pfieldmassmatrix \pfieldphinatfecoeffvec
    + \left(\frac{1}{\capillarynumber} \pfieldmassmatrix + \frac{1}{\bendingcapillarynumber} \pfieldellipticmatrix \right) \pfieldmufecoeffvec
                                          & = -\frac{\ofielddoublewellparam}{2 \polaritynumber} \pfieldferhsvec(\ofieldfe) \\
    \frac{1}{\pfieldepsilon}\pfieldfenonling(\pfieldfe) + \pfieldepsilon \pfieldellipticmatrix \pfieldfecoeffvec +
    \pfieldmassmatrix \pfieldmufecoeffvec & = \Vzero
  \end{align}
\end{subequations}
with
\begin{align*}
  \pfieldmassmatrix                    & \in \Reals[\pfieldfedimension \times \pfieldfedimension] \text{ s.t.\ }
  \pfieldmassmatrixentry{\feindex}{\feindexalt}
  =\bulkintegral{\pfieldfebasisfunc_\feindex \pfieldfebasisfunc_\feindexalt}                                                    \\
  \pfieldBmatrix(\stokesvelocity)      & \in \Reals[\pfieldfedimension \times \pfieldfedimension] \text{ s.t.\ }
  \pfieldBmatrixentry{\feindex}{\feindexalt}
  =\bulkintegral{ (\stokesvelocity \cdot \gradient \pfieldfebasisfunc_\feindex) \pfieldfebasisfunc_\feindexalt}
  \\
  \pfieldellipticmatrix                & \in \Reals[\pfieldfedimension \times \pfieldfedimension] \text{ s.t.\ }
  \pfieldellipticmatrixentry{\feindex}{\feindexalt}
  =\bulkintegral{\gradient \pfieldfebasisfunc_\feindex \cdot \gradient \pfieldfebasisfunc_\feindexalt}
  \\
  \pfieldfenonlinf(\pfield, \pfieldmu) & \in \Reals[\pfieldfedimension]
  \text{ s.t.\ } \pfieldfenonlinfentry{\feindex}(\pfield, \pfieldmu)
  = \bulkintegral{\doublewellpotentialtwoprime(\pfield)\pfieldmu \pfieldfebasisfunc_\feindex}
  \\
  \pfieldfenonling(\pfield)            & \in \Reals[\pfieldfedimension] \text{ s.t.\ } \pfieldfenonlingentry{\feindex}(\pfield)
  = \bulkintegral{\doublewellpotentialprime(\pfield) \pfieldfebasisfunc_\feindex}
  \\
  \pfieldferhsvec(\ofield)             & \in \Reals[\pfieldfedimension] \text{ s.t.\ }
  \pfieldferhsvecentry{\feindex}(\ofield) = \bulkintegral{\left(\ofield \cdot \ofield
    \right) \pfieldfebasisfunc_\feindex}
\end{align*}

\subsubsection{Orientation field}
Inserting the basis representations in~\eqref{eq:ofield_weak} gives the system
\begin{subequations}\label{eq:ofield_weak_algebraic_semidiscrete}
  \begin{align}
    \dt \ofieldmassmatrix \ofieldfecoeffvec + \ofieldBmatrix(\stokesvelocity) \ofieldfecoeffvec +  \frac{1}{\fracrotationaldynamicviscosity}
    \ofieldmassmatrix \ofieldPnatfecoeffvec
     & = \Vzero \\
    \ofieldmassmatrix \ofieldPnatfecoeffvec + \frac{\ofielddoublewellparam}{\polaritynumber} \ofieldCmatrix(\pfieldfe) \ofieldfecoeffvec
    - \frac{\ofielddoublewellparam}{\polaritynumber} \ofieldfenonlinf(\ofieldfe)
    - \frac{1}{\polaritynumber} \ofieldellipticmatrix \ofieldfecoeffvec
     & = \Vzero
  \end{align}
\end{subequations}
with
\begin{align*}
  \ofieldmassmatrix               & \in \Reals[\ofieldfedimension \times \ofieldfedimension] \text{ s.t.\ }
  \ofieldmassmatrixentry{\feindex}{\feindexalt}
  =\bulkintegral{\ofieldfebasisfunc_\feindex \cdot \ofieldfebasisfunc_\feindexalt}                          \\
  \ofieldellipticmatrix           & \in \Reals[\ofieldfedimension \times \ofieldfedimension] \text{ s.t.\ }
  \ofieldellipticmatrixentry{\feindex}{\feindexalt}
  =\bulkintegral{\doubledotproduct{\Vgradient \ofieldfebasisfunc_\feindex}{\Vgradient \ofieldfebasisfunc_\feindexalt}}
  \\
  \ofieldBmatrix(\stokesvelocity) & \in \Reals[\ofieldfedimension \times \ofieldfedimension] \text{ s.t.\ }
  \ofieldBmatrixentry{\feindex}{\feindexalt}
  =\bulkintegral{ \left(\left(\Vgradient \ofieldfebasisfunc_\feindexalt\right) \stokesvelocity +
    \Big(\vorticitytensor(\stokesvelocity)
    - \shapefactor \deformationtensor(\stokesvelocity)\Big)\ofieldfebasisfunc_\feindexalt\right) \cdot
  \ofieldfebasisfunc_\feindex}                                                                              \\
  \ofieldCmatrix(\pfield)         & \in \Reals[\ofieldfedimension \times \ofieldfedimension] \text{ s.t.\ }
  \ofieldCmatrixentry{\feindex}{\feindexalt}(\pfield)
  =\bulkintegral{\pfield \, \ofieldfebasisfunc_\feindex \cdot \ofieldfebasisfunc_\feindexalt}               \\
  \ofieldfenonlinf(\ofield)       & \in \Reals[\ofieldfedimension]
  \text{ s.t.\ } \ofieldfenonlinfentry{\feindex}(\ofield)
  = \bulkintegral{\left(\ofield \cdot \ofield\right) \ofield \cdot
    \ofieldfebasisfunc_\feindex}
\end{align*}

\subsubsection{Stokes}
The discretized Stokes system reads
\begin{subequations}\label{eq:stokes_weak_algebraic_semidiscrete}
  \begin{align}
    \stokesAmatrix \stokesvelocityfecoeffvec + \stokesBmatrix^T \stokespressurefecoeffvec & = \stokesferhsvec\left(\pfieldfe, \pfieldphinatfe,
    \ofieldfe, \ofieldPnatfe\right)                                                                                                            \\
    \stokesBmatrix \stokesvelocityfecoeffvec                                              & = \Vzero
  \end{align}
\end{subequations}
where
\begin{align*}
  \stokesAmatrix & \in \Reals[\stokesvelocityfedimension \times
    \stokesvelocityfedimension] \text{ s.t.\ }
  \stokesAmatrixentry{\feindex}{\feindexalt}
  = \frac{1}{2} \bulkintegral{\doubledotproduct{\Vgradient \stokesvelocityfebasisfunc_{\feindex}}{\Vgradient
      \stokesvelocityfebasisfunc_{\feindexalt}}}
  \\
  \stokesBmatrix & \in \Reals[\pfieldfedimension \times \stokesvelocityfedimension]
  \text{ s.t.\ } \stokesBmatrixentry{\feindex}{\feindexalt} =
  \bulkintegral{\pfieldfebasisfunc_{\feindex} \diverg \stokesvelocityfebasisfunc_{\feindexalt}}
  \\
  \begin{split}
    \stokesferhsvec\left(\pfield, \pfieldphinat, \ofield, \ofieldPnat\right) & \in \Reals[\stokesvelocityfedimension] \text{ s.t.\ }
    \stokesferhsvecentry{\feindex} = \bulkintegral{\Vdiverg\left(\ericksenstress(\pfield,\pfieldphinat,\ofield,\ofieldPnat)\right) \cdot
      \stokesvelocityfebasisfunc_{\feindex}} \\
    &\hspace{2.5cm}- \bulkintegral{\doubledotproduct{\left(\activestress(\pfield,\ofield) + \distortionstress(\ofield, \ofieldPnat)\right)}{\Vgradient
        \stokesvelocityfebasisfunc_{\feindex}}}
  \end{split}
\end{align*}

\subsection{Time discretization}\label{sec:timediscretization}
Exactly as~\cite{Marth2015}, we use a semi-implicit splitting scheme to decouple the phase
field, orientation field and flow equations. We introduce
a partition \(0 = \timevar_0 < \timevar_1 < \cdots < \timevar_{\numtimesteps} = \tf	\) of the
time interval \([0, \tf]\). For simplicity, we will
use a fixed time step \(\timestep\), i.e, \(\timevar^{\timeindex} = \timeindex \timestep\). We
will indicate the variables at the discrete times by
superscripts, e.g., \(\pfield^{\timeindex}(\spatialvar) \approx
\pfield(\timevar^{\timeindex}, \spatialvar)\) is the approximation of the phase field \(\pfield\) at time
\(\timevar^{\timeindex}\).
In each time step, we first solve the phase field equations~\eqref{eq:pfield_weak_algebraic_semidiscrete} where we
discretize the phase field
variables implicitly and the other variables explicitly.
\begin{subequations}\label{eq:pfield_discrete}
  \begin{align}
    \left(\pfieldmassmatrix - \timestep \pfieldBmatrix(\stokesvelocityfe^{\timeindex})\right) \pfieldfecoeffvec^{\timeindex+1} + \timestep
    \pfieldmobility \pfieldellipticmatrix
    \pfieldphinatfecoeffvectime{\timeindex+1}            & = \pfieldmassmatrix \pfieldfecoeffvec^{\timeindex}
    \\
    \frac{1}{\bendingcapillarynumber \cdot \pfieldepsilon^2} \pfieldfenonlinf(\pfieldfe^{\timeindex+1}, \pfieldmufe^{\timeindex+1}) +
    \pfieldmassmatrix
    \pfieldphinatfecoeffvectime{\timeindex+1}
    + \left(\frac{1}{\capillarynumber} \pfieldmassmatrix + \frac{1}{\bendingcapillarynumber} \pfieldellipticmatrix \right)
    \pfieldmufecoeffvec^{\timeindex+1}
                                                         & = -\frac{\ofielddoublewellparam}{2 \polaritynumber}\pfieldferhsvec(\ofieldfe^{\timeindex})
    \\
    \frac{1}{\pfieldepsilon}\pfieldfenonling(\pfieldfe^{\timeindex+1}) + \pfieldepsilon \pfieldellipticmatrix \pfieldfecoeffvec^{\timeindex+1} +
    \pfieldmassmatrix \pfieldmufecoeffvec^{\timeindex+1} & = \Vzero
  \end{align}
\end{subequations}
Second, we solve the orientation field equations~\eqref{eq:ofield_weak_algebraic_semidiscrete}, already using the phase
field variables at the new
time.
\begin{subequations}\label{eq:ofield_discrete}
  \begin{align}
    \ofieldmassmatrix \ofieldfecoeffvec^{\timeindex+1} + \timestep \ofieldBmatrix(\stokesvelocity^{\timeindex}) \ofieldfecoeffvec^{\timeindex+1} +
    \timestep \frac{1}{\fracrotationaldynamicviscosity}
    \ofieldmassmatrix \ofieldPnatfecoeffvectime{\timeindex+1}
     & = \ofieldmassmatrix \ofieldfecoeffvec^{\timeindex} \\
    \ofieldmassmatrix \ofieldPnatfecoeffvectime{\timeindex+1} + \frac{\ofielddoublewellparam}{\polaritynumber}
    \ofieldCmatrix(\pfieldfe^{\timeindex+1}) \ofieldfecoeffvec^{\timeindex+1}
    - \frac{\ofielddoublewellparam}{\polaritynumber} \ofieldfenonlinf(\ofieldfe^{\timeindex+1})
    - \frac{1}{\polaritynumber} \ofieldellipticmatrix \ofieldfecoeffvec^{\timeindex+1}
     & = \Vzero
  \end{align}
\end{subequations}
Finally, we solve the Stokes equations~\eqref{eq:stokes_weak_algebraic_semidiscrete} with the new phase field and
orientation field variables.
\begin{subequations}\label{eq:stokes_discrete}
  \begin{align}
    \stokesAmatrix \stokesvelocityfecoeffvec^{\timeindex+1} + \stokesBmatrix^T \stokespressurefecoeffvec^{\timeindex+1} & =
    \stokesferhsvec\left(\pfieldfe^{\timeindex+1}, \pfieldphinatfetime{\timeindex+1}, \ofieldfe^{\timeindex+1},
    \ofieldPnatfetime{\timeindex+1}\right)
    \\
    \stokesBmatrix \stokesvelocityfecoeffvec^{\timeindex+1}                                                             & = \Vzero
  \end{align}
\end{subequations}

\subsection{Solvers and preconditioners}\label{sec:solvers}
For the nonlinear systems~\eqref{eq:pfield_discrete} and~\eqref{eq:ofield_discrete}, we use a
backtracking Newton scheme using the values from the
last time step as an initial guess. The Jacobians needed for the Newton scheme could be inverted
using direct linear solvers~\cite{Marth2015}.
However, the Jacobians are large and sparse, and the variables from the last Newton iteration (or
the last time step) make good initial guesses. In
this setting, iterative linear solvers are usually both faster and less memory-intensive than
direct solvers. To choose suitable iterative solvers,
we will thus a closer look at the Jacobians of the nonlinear systems.

\subsubsection{Phase field}\label{sec:pfieldsolvers}
The Jacobian of the phase field system~\eqref{eq:pfield_discrete} has the form
\begin{align}
  \label{eq:phasefieldjacobian}
  \pfieldjacobianmatrix(\pfield, \pfieldmu)
  = \begin{pmatrix}
      \pfieldmassmatrix - \timestep \pfieldBmatrix(\stokesvelocity)                                                           &
      \timestep \pfieldmobility \pfieldellipticmatrix                                                                         &
      \Mzero                                                                                                                    \\
      \frac{1}{\bendingcapillarynumber \cdot \pfieldepsilon^2}\jacobianop_{\pfieldfecoeffvec}\pfieldfenonlinf(\pfield, \pfieldmu)
                                                                                                                              &
      \pfieldmassmatrix                                                                                                       &
      \frac{1}{\capillarynumber} \pfieldmassmatrix + \frac{1}{\bendingcapillarynumber} \pfieldellipticmatrix
      + \frac{1}{\bendingcapillarynumber \cdot \pfieldepsilon^2}\jacobianop_{\pfieldmufecoeffvec}\pfieldfenonlinf(\pfield)
      \\
      \pfieldepsilon \pfieldellipticmatrix + \frac{1}{\pfieldepsilon}\jacobianop_{\pfieldfecoeffvec}\pfieldfenonling(\pfield) &
      \Mzero                                                                                                                  &
      \pfieldmassmatrix
    \end{pmatrix}
\end{align}
where
\begin{subequations}\label{eq:pfield_nonlinear_jacobians}
  \begin{align}
    {(\jacobianop_{\pfieldfecoeffvec}\pfieldfenonlinf)}_{\feindex\feindexalt}
     & = \bulkintegral{\doublewellpotentialthreeprime(\pfield)\pfieldmu \pfieldfebasisfunc_{\feindex} \pfieldfebasisfunc_{\feindexalt}}
    = \bulkintegral{6\pfield\pfieldmu \pfieldfebasisfunc_{\feindex} \pfieldfebasisfunc_{\feindexalt}},                                  \\
    {(\jacobianop_{\pfieldmufecoeffvec}\pfieldfenonlinf)}_{\feindex\feindexalt}
    = {(\jacobianop_{\pfieldfecoeffvec}\pfieldfenonling)}_{\feindex\feindexalt}
     & = \bulkintegral{\doublewellpotentialtwoprime(\pfield) \pfieldfebasisfunc_{\feindex} \pfieldfebasisfunc_{\feindexalt}}
    = \bulkintegral{\left(3\pfield^2-1\right) \pfieldfebasisfunc_{\feindex} \pfieldfebasisfunc_{\feindexalt}}.
  \end{align}
\end{subequations}
Obviously, \(\pfieldjacobianmatrix\) is non-symmetric. The mass matrix \(\pfieldmassmatrix\) and
the stiffness matrix \(\pfieldellipticmatrix\) can
be computed once and for all before starting the
timestepping (as long as the grid does not change). The matrix \(\pfieldBmatrix\) depends on the
fluid velocity and
thus has to recomputed for each new time step (before starting the Newton iterations). The
Jacobians of the nonlinear terms \(\pfieldfenonlinf\) and
\(\pfieldfenonling\) have to be reassembled in each Newton iteration.

The convergence rate of iterative solvers usually depends on the condition number of the matrix.
Thus, preconditioners can greatly speed up (or merely enable) the solution process. To select a
preconditioner for the phase field
jacobian~\eqref{eq:phasefieldjacobian},
we first have a look at the scaling of the matrix blocks (cmp.~\cite{Boyanova2011}) which is
approximately
\begin{align*}
  \begin{pmatrix}
    \gridwidth^2 - \timestep \abs{\stokesvelocity} \gridwidth                                        &
    \timestep \pfieldmobility                                                                        &
    0                                                                                                  \\
    \frac{\abs{\pfield}\abs{\pfieldmu} \gridwidth^2}{\bendingcapillarynumber \cdot \pfieldepsilon^2} &
    \gridwidth^2                                                                                     &
    \frac{\gridwidth^2}{\capillarynumber} + \frac{1}{\bendingcapillarynumber} + \frac{\abs{3\pfield^2 - 1}\gridwidth^2}{\bendingcapillarynumber \cdot
    \pfieldepsilon^2}                                                                                  \\
    \pfieldepsilon + \frac{\abs{3\pfield^2 - 1}\gridwidth^2}{\pfieldepsilon}                         &
    0                                                                                                &
    \gridwidth^2
  \end{pmatrix}
\end{align*}
If we assume that \(\timestep\) is small (e.g., \(\landauO(\gridwidth^2)\)) and further assume
that the fluid velocity is of order \(1\), we see that
the \(\timestep \pfieldBmatrix(\stokesvelocity)\) term scales with \(\gridwidth^3\) and can thus be
neglected (for \(\gridwidth\) small enough)
compared to the mass matrix \(\pfieldmassmatrix\) which is of order \(\gridwidth^2\).

To guarantee a sufficient resolution of the interface region, the phase field parameter
\(\pfieldepsilon\) (or the grid width \(\gridwidth\)) is
usually chosen such that \(\gridwidth \leq c \pfieldepsilon\) where \(c\) is a small
integer\cite{Boyanova2011}. We can thus assume that these two
parameters scale similarly. Moreover, \(\pfield\) should take on values between
\([-1, 1]\) and is in particular bounded. If we further assume that
the capillary numbers \(\bendingcapillarynumber\) and \(\capillarynumber\) are of order \(1\), we
see that it is not clear that we can neglect the
Jacobians of the nonlinear terms as they seem to scale similar to the other terms in their
respective block. However, outside the interface region,
we have \(\pfield^2 \approx 1\) and \(\pfieldmu \approx 0\) and thus
\begin{subequations}
  \begin{align}
    \jacobianop_{\pfieldfecoeffvec}\pfieldfenonlinf \approx \Mzero \\
    \jacobianop_{\pfieldmufecoeffvec}\pfieldfenonlinf =
    \jacobianop_{\pfieldfecoeffvec}\pfieldfenonling
    \approx 2\pfieldmassmatrix
  \end{align}
\end{subequations}
are reasonable approximations. Taken together, we obtain the preconditioner
\begin{align}
  \label{eq:phasefieldjacobianpreconditioner}
  \pfieldjacobianlinearpreconditioner
  = \begin{pmatrix}
      \pfieldmassmatrix                                                                 & \timestep \pfieldmobility \pfieldellipticmatrix & \Mzero
      \\
      \Mzero                                                                            & \pfieldmassmatrix                               &
      \left(\frac{1}{\capillarynumber}
      + \frac{2}{\bendingcapillarynumber\cdot \pfieldepsilon^2}\right) \pfieldmassmatrix + \frac{1}{\bendingcapillarynumber} \pfieldellipticmatrix
      \\
      \pfieldepsilon \pfieldellipticmatrix + \frac{2}{\pfieldepsilon} \pfieldmassmatrix & \Mzero                                          &
      \pfieldmassmatrix
    \end{pmatrix}.
\end{align}
By the arguments above it is not totally clear that this preconditioner sufficiently approximates
the matrix entries corresponding to the interface
region as the Jacobians of the nonlinear parts (which we omitted) might contribute significantly
here. However, our numerical tests indicate
that~\eqref{eq:phasefieldjacobianpreconditioner} is indeed a good preconditioner for~\eqref{eq:phasefieldjacobian} (see~\cref{sec:numericalexperiments}).

\subsubsection{Orientation field}\label{sec:ofieldsolvers}
In each step of the Newton scheme, we have to invert the Jacobian
\begin{align}
  \label{eq:ofieldjacobian}
  \ofieldjacobianmatrix(\ofield)
  = \begin{pmatrix}
      \ofieldmassmatrix + \timestep \ofieldBmatrix(\stokesvelocity) & \frac{\timestep}{\fracrotationaldynamicviscosity} \ofieldmassmatrix \\
      \frac{\ofielddoublewellparam}{\polaritynumber} \ofieldCmatrix(\pfield)
      - \frac{1}{\polaritynumber} \ofieldellipticmatrix
      - \frac{\ofielddoublewellparam}{\polaritynumber}\jacobianop_{\ofieldfecoeffvec}\ofieldfenonlinf(\ofield)
                                                                    & \ofieldmassmatrix
    \end{pmatrix}
\end{align}
where the Jacobian of the nonlinear term has the entries
\begin{align}
  \begin{split}
    \label{eq:ofield_nonlinear_jacobians}
    {(\jacobianop_{\ofieldfecoeffvec}\ofieldfenonlinf)}_{\feindex\feindexalt}
    &= \bulkintegral{\partialderiv{(\ofield \cdot \ofield)}{\ofieldfecoeff_\feindexalt}
      (\ofield \cdot \ofieldfebasisfunc_\feindex) + (\ofield \cdot \ofield) \partialderiv{(\ofield \cdot
        \ofieldfebasisfunc_\feindex)}{\ofieldfecoeff_\feindexalt} } \\
    &= \bulkintegral{2(\ofield \cdot \ofieldfebasisfunc_\feindexalt) (\ofield \cdot
      \ofieldfebasisfunc_\feindex)
      + (\ofield \cdot \ofield) (\ofieldfebasisfunc_\feindexalt \cdot
      \ofieldfebasisfunc_\feindex)}
  \end{split}
\end{align}

We can invert~\eqref{eq:ofieldjacobian} by a Schur complement-type solver. Introducing
\begin{align}
  \ofieldGmatrix(\pfield, \ofield) =
  \frac{\ofielddoublewellparam}{\polaritynumber} \ofieldCmatrix(\pfield)
  - \frac{1}{\polaritynumber} \ofieldellipticmatrix
  - \frac{\ofielddoublewellparam}{\polaritynumber} \jacobianop_{\ofieldfecoeffvec}\ofieldfenonlinf(\ofield)
\end{align}
we factorize~\eqref{eq:ofieldjacobian} as
\begin{align}
  \begin{pmatrix}
    \ofieldmassmatrix + \timestep \ofieldBmatrix(\stokesvelocity) & \frac{\timestep}{\fracrotationaldynamicviscosity} \ofieldmassmatrix \\
    \ofieldGmatrix(\pfield, \ofield)
                                                                  & \ofieldmassmatrix
  \end{pmatrix}
  =
  \begin{pmatrix}
    \eyematrix & \frac{\timestep}{\fracrotationaldynamicviscosity} \eyematrix \\
    \Mzero     & \eyematrix
  \end{pmatrix}
  \begin{pmatrix}
    \ofieldmassmatrix + \timestep \ofieldBmatrix(\stokesvelocity) - \frac{\timestep}{\fracrotationaldynamicviscosity} \ofieldGmatrix(\pfield,
    \ofield)
           & \Mzero            \\
    \Mzero & \ofieldmassmatrix
  \end{pmatrix}
  \begin{pmatrix}
    \eyematrix                                              & \Mzero     \\
    \ofieldmassmatrix^{-1} \ofieldGmatrix(\pfield, \ofield) & \eyematrix
  \end{pmatrix}.
\end{align}
The inverse of \(\ofieldjacobianmatrix\) thus is
\begin{align}
  \ofieldjacobianmatrix^{-1}
  =
  \begin{pmatrix}
    \eyematrix                                               & \Mzero     \\
    -\ofieldmassmatrix^{-1} \ofieldGmatrix(\pfield, \ofield) & \eyematrix
  \end{pmatrix}
  \begin{pmatrix}
    {\left(\ofieldmassmatrix + \timestep \ofieldBmatrix(\stokesvelocity) - \frac{\timestep}{\fracrotationaldynamicviscosity}
      \ofieldGmatrix(\pfield,
    \ofield)\right)}^{-1} & \Mzero                 \\
    \Mzero                & \ofieldmassmatrix^{-1}
  \end{pmatrix}
  \begin{pmatrix}
    \eyematrix & -\frac{\timestep}{\fracrotationaldynamicviscosity} \eyematrix \\
    \Mzero     & \eyematrix
  \end{pmatrix}
\end{align}

Note that the Schur complement matrix does not contain any inverse matrices. As a consequence, it
is a sparse matrix (as \(\ofieldmassmatrix,
\ofieldBmatrix, \ofieldGmatrix\) are sparse) and can easily be assembled explicitly which enables
the use of direct solvers. However, we will
see in \cref{sec:numericalexperiments} that iterative solvers are also a good option since they converge
quickly even when used without preconditioner.

\section{Model order reduction}\label{sec:modelreduction}
The cell model we are considering depends on a multitude of parameters
\begin{equation*}
  \param := (\pfieldepsilon, \pfieldmobility, \ofielddoublewellparam,
  \fracrotationaldynamicviscosity, \shapefactor, \bendingcapillarynumber,
  \capillarynumber, \polaritynumber, \activeforcenumber) \in \mathbb{R}^9
\end{equation*}
we might want to vary.
However, each new combination of parameters will require a full new
simulation of the model, which will make extensive studies of the
influence of the parameters or even the optimization of the parameters
w.r.t.~some quantities of interest prohibitively expensive.
To mitigate this issue, we will in this section derive a low-rank
reduced order model using reduced-basis techniques, which
can be used as a high-quality surrogate for the original discrete
model \eqref{eq:pfield_discrete}--\eqref{eq:stokes_discrete}.

\subsection{Definition of the reduced order model}\label{sec:rom}
To simplify notation, we first stack the coefficient vectors such that there is only one vector for
the phase field, orientation field and Stokes
coefficients, respectively.
\begin{align*}
  \pfieldvec^{\timeindex} = \begin{pmatrix}
                              \pfieldfecoeffvec^{\timeindex} \\ \pfieldphinatfecoeffvectime{\timeindex} \\ \pfieldmufecoeffvec^{\timeindex}
                            \end{pmatrix} \in \mathbb{R}^{\pfieldvecdimension},\ \
  \ofieldvec^{\timeindex} = \begin{pmatrix}
                              \ofieldfecoeffvec^{\timeindex} \\ \ofieldPnatfecoeffvectime{\timeindex}
                            \end{pmatrix} \in \mathbb{R}^{\ofieldvecdimension},\ \
  \stokesvec^{\timeindex} = \begin{pmatrix}
                              \stokesvelocityfecoeffvec^{\timeindex} \\ \stokespressurefecoeffvec^{\timeindex}
                            \end{pmatrix} \in \mathbb{R}^{\stokesvecdimension}
\end{align*}
Here, the length of the vectors is $\pfieldvecdimension = 3\pfieldfedimension$, $\ofieldvecdimension = 2\ofieldfedimension$
and $\stokesvecdimension = \stokesvelocityfedimension + \pfieldfedimension$ (compare \cref{sec:finiteelementdiscretization}).
We then rewrite equations
\eqref{eq:pfield_discrete}--\eqref{eq:stokes_discrete}
from \cref{sec:timediscretization} in the following condensed residual form:
\begin{subequations}\label{eq:fom}
  \begin{align}
    \pfieldresidual(\pfieldvec^{\timeindex+1}(\param), \pfieldvec^{\timeindex}(\param),
    \ofieldvec^{\timeindex}(\param), \stokesvec^{\timeindex}(\param); \param)
                           & = \Vzero,                                \\
    \ofieldresidual(\ofieldvec^{\timeindex+1}(\param), \ofieldvec^{\timeindex}(\param),
    \pfieldvec^{\timeindex+1}(\param), \stokesvec^{\timeindex}(\param); \param)
                           & = \Vzero,                                \\
    \stokesresidual(\stokesvec^{\timeindex+1}(\param), \pfieldvec^{\timeindex+1}(\param),
    \ofieldvec^{\timeindex+1}(\param); \param)
                           & = \Vzero,                                \\
    \pfieldvec^{0}(\param) = \begin{pmatrix}
                               \pfieldfecoeffvec^{0} \\ \Vzero \\ \Vzero
                             \end{pmatrix},\,
    \ofieldvec^{0}(\param) =
    \begin{pmatrix}
      \ofieldfecoeffvec^{0} \\ \Vzero
    \end{pmatrix},\,
    \stokesvec^{0}(\param) & = \begin{pmatrix}
                                 \stokesvelocityfecoeffvec^{0} \\ \Vzero
                               \end{pmatrix}
  \end{align}
\end{subequations}
where the initial values \(\pfieldfecoeffvec^{0}, \ofieldfecoeffvec^{0},
\stokesvelocityfecoeffvec^{0}\) for the coefficient vectors are obtained from the initial values
\eqref{eq:pfieldinitial}, \eqref{eq:ofieldinitial}, \eqref{eq:stokesvelocityinitial}
by interpolation.

Assume that \(\param \in \paramspace \subseteq \mathbb{R}^9\) for some parameter
domain of
interest \(\paramspace\) we want to explore, and further assume that we are given low-rank matrices
\(\Bredp \in \mathbb{R}^{\pfieldvecdimension\times\redpfieldvecdimension},
\Bredo \in \mathbb{R}^{\ofieldvecdimension\times\redofieldvecdimension},
\Breds \in \mathbb{R}^{\stokesvecdimension\times\redstokesvecdimension}\) with orthonormal
columns such that each
\(\pfieldvec^\timeindex(\param), \ofieldvec^\timeindex(\param),
\stokesvec^\timeindex(\param)\)
can be well-approximated within the column-spans of these matrices for all \(0 \leq k \leq K\) and
\(\param\in\paramspace\). The construction of these basis matrices will be discussed in
\cref{sec:hapoddef}.

Our aim is to find vectors
\(\redpfieldvec^\timeindex(\param) \in \mathbb{R}^{\redpfieldvecdimension}\),
\(\redofieldvec^\timeindex(\param) \in \mathbb{R}^{\redofieldvecdimension}\) and
\(\redstokesvec^\timeindex(\param) \in \mathbb{R}^{\redstokesvecdimension}\)
such that
\begin{equation}\label{eq:red_state_approx}
  \pfieldvec^\timeindex(\param) \approx \Bredp \cdot \redpfieldvec^\timeindex(\param), \quad
  \ofieldvec^\timeindex(\param) \approx \Bredo \cdot \redofieldvec^\timeindex(\param), \quad
  \stokesvec^\timeindex(\param) \approx \Breds \cdot \redstokesvec^\timeindex(\param).
\end{equation}

To determine the reduced states
\(\redpfieldvec^\timeindex(\param)\), \(\redofieldvec^\timeindex(\param)\),
\(\redstokesvec^\timeindex(\param)\)
we use the original system equations \eqref{eq:fom}.
However, substituting \eqref{eq:red_state_approx} into \eqref{eq:fom} will lead to an
over-determined system.
To close the system and ensure its stability, we choose a residual-minimization approach and define
the reduced states as the solutions of

\begin{subequations}\label{eq:rom}
  \begin{align}
    \redpfieldvec^{\timeindex+1}
     & = \argmin_{\redpfieldvec\in\mathbb{R}^{\dim\Vredp}}
    \norm{\pfieldresidual(\Bredp\cdot\redpfieldvec, \Bredp\cdot\redpfieldvec^{\timeindex}(\param),
    \Bredo\cdot\redofieldvec^{\timeindex}(\param), \Breds\cdot\redstokesvec^{\timeindex}(\param); \param)},   \\
    \redofieldvec^{\timeindex+1}
     & = \argmin_{\redofieldvec\in\mathbb{R}^{\dim\Vredo}}
    \norm{\ofieldresidual(\Bredo\cdot\redofieldvec, \Bredo\cdot\redofieldvec^{\timeindex}(\param),
    \Bredp\cdot\redpfieldvec^{\timeindex+1}(\param), \Breds\cdot\redstokesvec^{\timeindex}(\param); \param)}, \\
    \redstokesvec^{\timeindex+1}
     & = \argmin_{\redstokesvec\in\mathbb{R}^{\dim\Vreds}}
    \norm{\stokesresidual(\Breds\cdot\redstokesvec, \Bredp\cdot\redpfieldvec^{\timeindex+1}(\param),
      \Bredo\cdot\redofieldvec^{\timeindex+1}(\param); \param)}. \label{eq:rom_stokes}
  \end{align}
\end{subequations}

These systems of equations can be solved iteratively using the Gauss-Newton method. However, while
the model
is now formulated on low-dimensional spaces, the computation of the residuals and Jacobians still
requires
high-dimensional computations, which will severely limit the computational speedup gained by this
formulation.
Thus, we apply an additional reduction step by employing the discrete empirical interpolation
method
(DEIM) \cite{Chaturantabut2010} to the residuals.
To this end, assume we are given collateral bases $\Credp \in \mathbb{R}^{\pfieldvecdimension \times \colpfieldvecdimension}$,
$\Credo \in \mathbb{R}^{\ofieldvecdimension \times \colofieldvecdimension}$ and $\Creds \in \mathbb{R}^{\stokesvecdimension \times \colstokesvecdimension}$
of dimensions $\colpfieldvecdimension, \colofieldvecdimension, \colstokesvecdimension$
and matrices
$\pfielddofmatrix =
  [\mathbf{e}_{\pfielddof{1}},\ldots,\mathbf{e}_{\pfielddof{\colpfieldvecdimension}}]
  \in \mathbb{R}^{\colpfieldvecdimension \times \colpfieldvecdimension}$,
$\ofielddofmatrix =
  [\mathbf{e}_{\ofielddof{1}},\ldots,\mathbf{e}_{\ofielddof{\colofieldvecdimension}}]
  \in \mathbb{R}^{\colofieldvecdimension \times \colofieldvecdimension}$ and
$\stokesdofmatrix =
  [\mathbf{e}_{\stokesdof{1}},\ldots,\mathbf{e}_{\stokesdof{\colstokesvecdimension}}]
  \in \mathbb{R}^{\colstokesvecdimension \times \colstokesvecdimension}$,
where
$1 \leq \pfielddof{i} \leq \pfieldvecdimension$ ($1 \leq i \leq \colpfieldvecdimension$),
$1 \leq \ofielddof{i} \leq \ofieldvecdimension$ ($1 \leq i \leq \colofieldvecdimension$) and
$1 \leq \stokesdof{i} \leq \stokesvecdimension$ ($1 \leq i \leq \colstokesvecdimension$)
are pairwise disjoint sets of indices and (by abuse of notation) $\mathbf{e}_i$ denotes the
$i$-th
canonical basis vector in the Euclidean space of appropriate dimension.
Thus, for instance, left-multiplying $\pfieldresidual$ with $\pfielddofmatrix^T$ extracts
the
$\colpfieldvecdimension$ DOFs of $\pfieldresidual$ with indices
$\pfielddof{1}, \ldots \pfielddof{\colpfieldvecdimension}$, and
$\Credp (\pfielddofmatrix^T\Credp)^{-1} \pfielddofmatrix^T \pfieldresidual$
is the interpolation of $\pfieldresidual$ within the column span of $\Credp$
that agrees
with $\pfieldresidual$ at all $\colpfieldvecdimension$ interpolation DOFs.
The computation of the interpolation bases and DOFs is again covered in \cref{sec:hapod}.

Replacing all residuals with their respective empirical interpolants leads us to the following
hyper-reduced residual-minimization problem, that we will use as surrogate for the full-order
model:

\begin{subequations}\label{eq:deim_rom}
  \begin{align}
    \redpfieldvec^{\timeindex+1}
     & = \argmin_{\redpfieldvec\in\mathbb{R}^{\dim\Vredp}}
    \norm{\Credp (\pfielddofmatrix^T\Credp)^{-1} \pfielddofmatrix^T
      \pfieldresidual(\Bredp\cdot\redpfieldvec, \Bredp\cdot\redpfieldvec^{\timeindex}(\param),
    \Bredo\cdot\redofieldvec^{\timeindex}(\param), \Breds\cdot\redstokesvec^{\timeindex}(\param); \param)},   \\
    \redofieldvec^{\timeindex+1}
     & = \argmin_{\redofieldvec\in\mathbb{R}^{\dim\Vredo}}
    \norm{\Credo (\ofielddofmatrix^T\Credo)^{-1} \ofielddofmatrix^T
      \ofieldresidual(\Bredo\cdot\redofieldvec, \Bredo\cdot\redofieldvec^{\timeindex}(\param),
    \Bredp\cdot\redpfieldvec^{\timeindex+1}(\param), \Breds\cdot\redstokesvec^{\timeindex}(\param); \param)}, \\
    \redstokesvec^{\timeindex+1}
     & = \argmin_{\redstokesvec\in\mathbb{R}^{\dim\Vreds}}
    \norm{\Creds (\stokesdofmatrix^T\Creds)^{-1} \stokesdofmatrix^T
      \stokesresidual(\Breds\cdot\redstokesvec, \Bredp\cdot\redpfieldvec^{\timeindex+1}(\param),
      \Bredo\cdot\redofieldvec^{\timeindex+1}(\param); \param)}. \label{eq:deim_rom_stokes}
  \end{align}
\end{subequations}

Thanks to the locality of all involved finite-element operators, the restricted residuals
$\pfielddofmatrix^T\pfieldresidual$,
$\ofielddofmatrix^T\ofieldresidual$,
$\stokesdofmatrix^T\stokesresidual$
and their Jacobians can be evaluated by only local computations in fixed-size neighborhoods of the
interpolation DOFs.
By storing the values of $\Bredp$, $\Bredo$, $\Breds$ at
the DOFs of these neighborhoods,
\eqref{eq:deim_rom} can be solved with an effort that scales in space and time only with the
reduced
dimensions $\redpfieldvecdimension, \redofieldvecdimension, \redstokesvecdimension, \colpfieldvecdimension,
  \colofieldvecdimension, \colstokesvecdimension$ and is independent of the dimensions of the finite-element
spaces.

\begin{remark}
  While the Stokes residual is affine linear in $\stokesvec$, it nonlinearly depends
  on $\pfieldvec$ and $\ofieldvec$.
  Hence, some form of hyper-reduction is also required for solving \eqref{eq:rom_stokes}.
  For ease of implementation, we chose here the same monolithic interpolation approach
  as for $\pfieldresidual$ and $\ofieldresidual$.
  The Gauss-Newton algorithm will converge after one iteration when applied to
  \eqref{eq:deim_rom_stokes}.
\end{remark}

\subsection{Distributed construction of approximation spaces}\label{sec:hapod}
The approximation quality of~\eqref{eq:deim_rom} is strongly dependent on the choice of
the reduced bases \(\Bredp\), \(\Bredo\), \(\Breds\), \(\Credp\), \(\Credo\), \(\Creds\).
In the following, we will outline the construction of these reduced bases using the hierarchical
approximate proper orthogonal decomposition (HAPOD) \cite{Himpe2018}.
To that end, we shortly recall the proper orthogonal decomposition (POD).

\subsubsection{POD}
The proper orthogonal decomposition (POD) builds
low-rank approximation spaces for the column vectors of a given
matrix of solution snapshots \(\snapshotmatrix \in \mathbb{R}^{\fedimension \times \snapshotcount}\)
from the left-singular vectors corresponding to the dominant singular values of
$\snapshotmatrix$.
We consider here the weighted version, where the singular value decomposition (SVD) is computed by
interpreting
$\snapshotmatrix$ as a linear operator on the finite-dimensional Hilbert space
$(\mathbb{R}^\fedimension, \innerproductmatrix)$,
where the symmetric positive matrix $\innerproductmatrix \in \mathbb{R}^{\fedimension \times \fedimension}$
induces the inner product and norm
\begin{align}\label{eq:innerproductdef}
  \Winnerproduct{\fevec_1}{\fevec_2} \coloneqq \fevec_1^T \innerproductmatrix \fevec_2 \quand
  \Wnorm{\fevec_1} \coloneqq \sqrt{\Winnerproduct{\fevec_1}{\fevec_1}}.
\end{align}
Using the so-called \emph{method of snapshots}~\cite{Sirovich1987} to determine the SVD of
$\snapshotmatrix$ via an eigendecomposition
of the associated weighted Gramian $\snapshotmatrix^T \cdot \innerproductmatrix \cdot \snapshotmatrix \in \mathbb{R}^{\snapshotcount \times \snapshotcount}$ yields the following algorithm:

\begin{algorithm}[H]
  \caption{$\POD(\snapshotmatrix, \innerproductmatrix, \varepsilon)$: weighted POD via method of snapshots}\label{alg:POD}
  \SetKwInOut{Input}{input}\SetKwInOut{Output}{output}
  \Input{snapshot matrix $\snapshotmatrix \in \mathbb{R}^{\fedimension \times \snapshotcount}$,
    inner product matrix $\innerproductmatrix \in \mathbb{R}^{\fedimension \times \fedimension}$,
    error tolerance $\varepsilon > 0$}
  \Output{POD modes $\MU \in \mathbb{R}^{\fedimension \times \reddimension}$,
    singular values $\boldsymbol{\sigma}$}
  $\MG := \snapshotmatrix^T \cdot \innerproductmatrix \cdot \snapshotmatrix$\tcp*{weighted $\snapshotcount \times \snapshotcount$ Gramian matrix}
  $\Vlambda, \MV := eigh(\MG)$\tcp*{$\Vlambda_i$, $\MV_{:,i}$ = $i$-th eigen-value/-vector ordered by magnitude}
  $\reddimension := \min \{M \in \mathbb{N}^{\geq 0} \vert \sum_{i=M+1}^\snapshotcount
    \Vlambda_i^2 \leq \varepsilon\}$\;
  $\Vlambda_\varepsilon, \MV_\varepsilon := \Vlambda(1:M), \MV(:, 1:M)$\tcp*{truncate according to given tolerance}
  $\boldsymbol{\sigma} := \sqrt{\Vlambda_\varepsilon}$\tcp*{singular values}
  $\MU := \MS \cdot \MV_\varepsilon \cdot
    diag(\boldsymbol{\sigma})^{-1}$\tcp*{left-signular vectors}
\end{algorithm}

By definition of the SVD, the POD modes $\MU$ are
$\innerproductmatrix$-orthonormal, i.e.
$\MU^T \innerproductmatrix \MU = \MI$.
In particular, the $\innerproductmatrix$-orthonormal projection $\orthprojection[\operatorname{colspan} \MU]$ onto the
approximation space
spanned by the POD modes is given by the formula
\begin{equation}\label{eq:Worthprojection}
  \orthprojection[\operatorname{colspan} \MU](\fevec) = \orthprojectionmatrix_{\MU} \cdot \fevec \coloneqq \MU \MU^T
  \innerproductmatrix \fevec.
\end{equation}

The main reason for the importance of the POD is the fact that it produces best approximating
spaces
for the snapshot data in the \(\ell^2\)-sense:

\begin{theorem}[Schmidt-Eckhard-Young-Mirsky]\label{thm:pod_error}
  Let $\MU, \boldsymbol{\sigma} := \POD(\snapshotmatrix, \innerproductmatrix, \varepsilon)$ be the
  POD modes and singular values
  for a given snapshot matrix \(\snapshotmatrix\), inner product matrix $\innerproductmatrix$ and
  truncation tolerance
  $\varepsilon$.
  Then $\operatorname{colspan} \MU$ is an \(\ell^2\)-best approximating space for \(\snapshotmatrix\) in
  the sense that
  \begin{equation}\label{eq:pod_error}
    \sum_{i = 1}^\snapshotcount \Wnorm{\snapshotmatrix_{:,i} - \orthprojection[\operatorname{colspan} \MU](\snapshotmatrix_{:,i})}^2
    = \min_{\substack{\subspace \subseteq \mathbb{R}^\fedimension \text{lin.\ subsp.}\\ \dim \subspace
        = \reddimension}}\ \sum_{i = 1}^\snapshotcount
    \Wnorm{\snapshotmatrix_{:,i} - \orthprojection[\subspace](\snapshotmatrix_{:,i})}^2 =
    \sum_{i=\reddimension+1}^{\snapshotcount} \Vlambda_i \leq \varepsilon^2,
  \end{equation}
  where $\Vlambda$ denotes the vector of untruncated squared singular values as in the
  definition of \cref{alg:POD}.
  The dimension of the POD space $\reddimension$ is the minimal dimension for which the
  $\ell^2$ error is not larger than
  the given tolerance $\varepsilon$.
\end{theorem}

\begin{remark}
  We note that we could have equally well used any other approach for the computation of the SVD
  in \cref{alg:POD},
  e.g. QR decomposition or a randomized method.
  We choose the method of snapshots here, as it is particularly easy to implement efficiently
  for small $\snapshotcount$, which is the case when it is used as part of the HAPOD algorithms
  introduced in the following section.
\end{remark}

\subsubsection{HAPOD}\label{sec:hapoddef}
Directly computing the POD for large snapshot sets can be computationally demanding.
This is, in particular, the case, when the snapshot data grows so large that it not longer fits
into memory.
The Hierarchical Approximate POD (HAPOD) introduced in~\cite{Himpe2018} is a flexible
approach
to partition the problem of computing the POD of a large dataset into smaller subproblems, which,
depending on the structure of the computational setup, can be either solved incrementally or in
parallel.
Each of these subproblems consists in computing a POD of a small subset of the snapshot data in
combination
with POD modes from previous subproblems.
The resulting local POD modes are then scaled by their singular values and sent to the next node in
the computation.

To formalize this procedure, we consider rooted trees where to each node of the tree a local POD is
associated.
A rooted tree is a connected acyclic graph $\T$ with nodes $\mathcal{N}_\T$
of which $\rho_\T \in \mathcal{N}_\T$
is designated as the root of the tree. For each node $\alpha \in \mathcal{N}_\T$, we denote by
$\childmap(\alpha) \subset \mathcal{N}_\T$ the children of $\alpha$, i.e., all nodes connected to
$\alpha$ which
cannot be reached from $\rho_\T$ without first visiting $\alpha$.
By $\mathcal{L}_\T\coloneqq \{\alpha \in \mathcal{N}_\T \ |\
  \children{\alpha} = \emptyset \}$ we denote the
leafs of $\T$, i.e., the nodes which do not have any children.

To define the HAPOD, let a snapshot matrix $\snapshotmatrix$ and an inner product matrix
$\innerproductmatrix$
be given.
For a given tree \(\T\), assign to each leaf $\gamma \in \mathcal{L}_\T$ a set of column indices
$D(\gamma) \subseteq \{1, \ldots, s\}$
s.t. $D(\gamma_1) \cap D(\gamma_2) = \emptyset$ for different leafs $\gamma_1 \neq
  \gamma_2$ and s.t.
$\bigcup_{\gamma \in \mathcal{L}_\T} D(\gamma) = \{1, \ldots, s\}$.
By $\snapshotmatrix_{D(\gamma)}$ we then denote the $\fedimension \times |D(\gamma)|$ matrix given by
\begin{equation*}
  \snapshotmatrix_{D(\gamma)} := \bigl[\snapshotmatrix_{:, {D(\gamma)_1}} \cdots \snapshotmatrix_{:, {D(\gamma)_{
              |D(\gamma)|}}}\bigr],
\end{equation*}
where ${D(\gamma)_1}, \ldots, {D(\gamma)_{|D(\gamma)|}}$ is an arbitrary enumeration of $D(\gamma)$.
Thus all snapshot vectors of $\snapshotmatrix$ are contained in exactly one local snapshot
matrix
$\snapshotmatrix_{D(\gamma)}$. Further choosing local POD tolerances $\varepsilon(\alpha) \in \mathbb{R}^{\geq 0}$
for each node $\alpha \in \mathcal{N}_\T$, we can recursively define:

\begin{algorithm}[H]
  \caption{$\HAPOD[\snapshotmatrix, \innerproductmatrix, \T, D, \varepsilon](\alpha)$: local HAPOD at node $\alpha \in
      \mathcal{N}_\T$}\label{alg:HAPOD}
  \SetKwInOut{Input}{input}\SetKwInOut{Output}{output}
  \Input{snapshot matrix $\snapshotmatrix \in \mathbb{R}^{\fedimension \times \snapshotcount}$,
    inner product matrix $\innerproductmatrix \in \mathbb{R}^{\fedimension \times \fedimension}$,
    snapshot distribution $D$,
    error tolerances $\varepsilon$,
    node $\alpha$}
  \Output{scaled HAPOD modes $\MU_\alpha \in \mathbb{R}^{\fedimension \times \reddimension}$}
  \eIf{$\alpha \in \mathcal{L}_\T$}{
    $\snapshotmatrix_\alpha := \snapshotmatrix_{D(\alpha)}$\;
  }{
    $\snapshotmatrix_\alpha := \bigl[
        \HAPOD[\snapshotmatrix, \innerproductmatrix, \T, D, \varepsilon](\children{\alpha}_1)\,\cdots\,
        \HAPOD[\snapshotmatrix, \innerproductmatrix, \T, D, \varepsilon](\children{\alpha}_{|\children{\alpha}|})\bigr]$\;
  }
  \eIf{$\varepsilon(\alpha) = 0$}{
    $\MU_\alpha, \boldsymbol{\sigma}_\alpha := \snapshotmatrix_\alpha,
      [1\, \cdots\, 1]$\tcp*{No POD when tolerance is $0$.}
  }{
    $\MU_\alpha, \boldsymbol{\sigma}_\alpha := \POD(\snapshotmatrix_\alpha, \innerproductmatrix, \varepsilon(\alpha))$\;
  }
  \If{$\alpha \neq \rho_\T$}{
    $\MU_\alpha := \MU_\alpha \cdot \operatorname{diag}{\boldsymbol{\sigma}_\alpha}$\tcp*{Scale POD modes when not at root.}
  }
\end{algorithm}

In this definition, $\children{\alpha}_i$ again denotes an arbitrary enumeration of
$\children{\alpha}$.
We note that the special case for $\varepsilon(\alpha) = 0$ is included in the algorithm to allow
incremental
POD computations where new snapshot data directly enters a POD with the current modes, without
first applying
a POD only to the new data.
The hierarchical approximate POD is now simply given by the local POD modes at the root node:
\begin{equation*}
  \HAPOD[\snapshotmatrix, \innerproductmatrix, \T, D, \varepsilon]:=
  \HAPOD[\snapshotmatrix, \innerproductmatrix, \T, D, \varepsilon](\rho_\T).
\end{equation*}
If the local POD tolerances are chosen as \cite{Himpe2018}
\begin{equation}
  \varepsilon(\alpha) = \begin{cases}
    \sqrt{\abs{\Snap}}\cdot \omega \cdot \varepsilon^\ast & \text{ if } \alpha = \rho_\T \\
    \sqrt{\abs{\AccSnap_\alpha}}\cdot {(L_\T - 1)}^{-1/2}\cdot \sqrt{1 - \omega^2}\cdot \varepsilon^\ast &\text{ else,}
  \end{cases}
\end{equation}
then the mean \(\ell^2\) projection error is bounded by \(\varepsilon^{\ast}\), i.e.\
\begin{equation}
  {\left(\frac{1}{\abs{\snapshotset}} \sum_{\snapshot \in \snapshotset} \Wnorm{\snapshot - \orthprojection[\operatorname{colspan} \MU](\snapshot)}^2\right)}^{\frac12} \leq \varepsilon^{\ast}.
\end{equation}
Here, the parameter \(0 \leq \omega \leq 1\) can be chosen arbitrarily and allows to find the best compromise between
efficiency of the HAPOD and the optimality of the resulting approximation space.

\subsubsection{Application of HAPOD-DEIM to our model}
We now want to use the HAPOD to construct the reduced bases
\(\Bredp\), \(\Bredo\), \(\Breds\), \(\Credp\), \(\Credo\), \(\Creds\)
for our model.

To that end, we first choose
symmetric positive definite weight matrices \(\innerproductmatrixpfield\),
\(\innerproductmatrixofield\) and \(\innerproductmatrixstokes\)
which, by \eqref{eq:innerproductdef}, define inner products
such that the coefficient spaces \(\hapodhilbertspacepfield = (\R^{\pfieldvecdimension},
\Winnerproductpfield{\cdot}{\cdot})\),
\(\hapodhilbertspaceofield = (\R^{\ofieldvecdimension}, \Winnerproductofield{\cdot}{\cdot})\)
and \(\hapodhilbertspacestokes = (\R^{\stokesvecdimension}, \Winnerproductstokes{\cdot}{\cdot})\)
are Hilbert spaces.

If we choose the weight matrices as the unit matrices, we obtain the standard Euclidean inner
products. On the other hand, if we use the mass matrices \(\pfieldmassmatrix, \ofieldmassmatrix,
\stokesmassmatrix\) of the finite element
discretization's basis functions, the resulting inner product is the \(\Lp{2}\)
inner product of the finite element approximation functions.

Now, to compute the reduced bases, we would like to solve the high-dimensional
problem~\eqref{eq:fom} for all parameters in a training parameter set
\begin{equation*}
  \trainparams = \Set{\param_\paramindex \given \paramindex\in\{0,\ldots,\ntrain-1\}},
\end{equation*}
which is assumed to
be representative for the set of parameters, and then perform a POD for each of the snapshot
multisets obtained by collecting the computed vectors and residuals:
\begin{subequations}\label{eq:snapshotsets}
  \begin{align}
    \snapshotsetany & = \Set{\anyvec^\timeindex(\param_\paramindex) \in \hapodhilbertspace_{\anyvec} \given
    \timeindex\in\{0,\ldots,\numtimesteps\},\ \paramindex\in\{0,\ldots,\ntrain-1\}}                                                                                                                                                                                      \\
    \residualsetany & = \Set{\anyresidual^{\timeindex,\newtonindex}(\param_\paramindex) \in \hapodhilbertspace_{\anyvec} \given \timeindex\in\{0,\ldots,\numtimesteps-1\}, \paramindex\in\{0,\ldots,\ntrain-1\}, \newtonindex\in\{0,\ldots,\newtonitsany^\timeindex-1\}}
  \end{align}
\end{subequations}
for \(\anyvec \in \Set{\pfieldvec, \ofieldvec, \stokesvec}\)
where \(\newtonitspfield^\timeindex, \newtonitsofield^\timeindex, \newtonitsstokes^\timeindex\) are
the number of iterations in the Newton scheme for the respective field in the \(\timeindex\)-th
time step and
\begin{align*}
  \pfieldresidual^{\timeindex,\newtonindex}(\param_\paramindex) & = \pfieldresidual(\pfieldvec^{\timeindex,\newtonindex}(\param_\paramindex), \pfieldvec^{\timeindex}(\param_\paramindex),
  \ofieldvec^{\timeindex}(\param_\paramindex), \stokesvec^{\timeindex}(\param_\paramindex); \param_\paramindex)                                                                              \\
  \ofieldresidual^{\timeindex,\newtonindex}(\param_\paramindex) & = \ofieldresidual(\ofieldvec^{\timeindex,\newtonindex}(\param_\paramindex), \pfieldvec^{\timeindex+1}(\param_\paramindex),
  \ofieldvec^{\timeindex}(\param_\paramindex), \stokesvec^{\timeindex}(\param_\paramindex); \param_\paramindex)                                                                              \\
  \stokesresidual^{\timeindex,\newtonindex}(\param_\paramindex) & = \stokesresidual(\stokesvec^{\timeindex,\newtonindex}(\param_\paramindex), \pfieldvec^{\timeindex+1}(\param_\paramindex),
  \ofieldvec^{\timeindex+1}(\param_\paramindex); \param_\paramindex)
\end{align*}
with \(\anyvec^{\timeindex,\newtonindex}\) the \(\newtonindex\)-th iterate in the Newton scheme for field
\(\anyvec\) in time step \(\timeindex\).
Note that, since the Stokes residual is affine linear in \(\stokesvec\), we have
\(\newtonitsstokes^\timeindex = 1\).

We could then perform a POD for each of the snapshot sets to obtain our reduced basis. However,
since we have relatively small time steps \(\timestep\) and since we often have to use hundreds of
training parameters
to get a good approximation of the parameter space, the snapshot sets \eqref{eq:snapshotsets} will
be very large. As a consequence, directly computing the PODs for the snapshot sets, for example via the method of snapshots (see \cref{alg:POD}), will
take a very long time or even be practically impossible. We will thus apply the HAPOD to obtain an
approximate POD decomposition with much less computing effort. Instead of computing the whole
snapshot sets \eqref{eq:snapshotsets} at a
time, which might exceed the available memory and require input/output operations to a mass storage
device, for the HAPOD we choose a chunk size \(\chunksize\) and only compute \(\chunksize\) time
steps at a time, obtaining the much smaller multisets
\begin{subequations}
  \begin{align}
    \snapshotsetanyi{\chunkindex}(\param_\paramindex) & = \Set{\anyvec^\timeindex(\param_\paramindex) \in \hapodhilbertspace_{\anyvec} \given
    \timeindex\in\{\chunkindex\cdot\chunksize,\ldots,\min\Big((\chunkindex+1)\cdot\chunksize, \numtimesteps\Big)\}}                                               \\
    \residualsetanyi{\chunkindex}(\param_\paramindex) & = \Set{\anyresidual^{\timeindex,\newtonindex}(\param_\paramindex) \in \hapodhilbertspace_{\anyvec} \given
      \timeindex\in\{\chunkindex\cdot\chunksize,\ldots,\min\Big((\chunkindex+1)\cdot\chunksize, \numtimesteps\Big)\}}
  \end{align}
\end{subequations}
for \(\paramindex\in\{0,\ldots,\ntrain-1\}\) and \(\anyvec \in \Set{\pfieldvec, \ofieldvec, \stokesvec}\). We then
already compute a local POD for each of these multisets separately, then scale the resulting
modes by their corresponding singular values, and collect these scaled modes to perform another POD
on the main rank of the compute node. Only then we compute the next chunk of time steps and again
perform a local POD and a
POD with the collected modes. This time, in the second POD with the collected modes, we also
include the scaled modes from the last timesteps. We proceed like this until
we computed all chunks of time steps. If we use several compute nodes, we then perform a final POD
with the collected output modes from all compute nodes.
The whole algorithm can be found in \cref{alg:HAPODCellModel}.

\begin{algorithm}
  \caption{HAPOD applied to cell model}\label{alg:HAPODCellModel}
  \DontPrintSemicolon
  \SetKw{KwIn}{in}
  \SetKw{KwFalse}{false}
  \SetKw{KwTrue}{true}
  \SetKw{KwNone}{None}
  \SetKwFunction{Fntol}{get\_local\_tolerance}
  \SetKwFunction{FnMain}{main}
  \SetKwProg{Fn}{function}{:}{}
  \Fn{\Fntol{\(\varepsilon^\ast\), \(\hapodtreedepth\), \texttt{num\_snaps\_in\_leafs}, \texttt{is\_root\_node}}}{
    \uIf{\texttt{is\_root\_node}}{
      \KwRet \(\varepsilon^\ast \cdot \omega \cdot \sqrt{\texttt{num\_snaps\_in\_leafs}}\)\;
    }
    \Else{
      \KwRet \(\varepsilon^\ast \cdot \sqrt{(1-\omega^2)\frac{\texttt{num\_snaps\_in\_leafs}}{\hapodtreedepth-1}}\)\;
    }
  }
  \;
  \Fn{\FnMain{\(\cdots\)}}{
    \(\param \gets\) \texttt{distribute\_parameters}(\(\cdots\))\tcp*[l]{For simplicity, assume one parameter \(\param\) per rank}
    \(\numtimesteps = \left\lceil\frac{\tf}{\timestep}\right\rceil\)\tcp*[l]{number of time steps}
    \(\numchunks = \left\lceil \frac{\numtimesteps+1}{\chunksize}\right\rceil\)\tcp*[l]{\(\numtimesteps+1\) to include initial values}
    \texttt{input\_for\_node\_level\_POD} \(\gets\) [[\,], [\,], [\,], [\,], [\,], [\,]] \tcp*[l]{List of \(6\) empty lists}
    \texttt{num\_snaps\_in\_leafs} \(\gets\) [0, 0, 0, 0, 0, 0]\;
    \For{\(\chunkindex \gets 0\) \KwTo \(\numchunks-1\)}{
      Compute next \(\chunksize\) time steps according to
      \eqref{eq:pfield_discrete}--\eqref{eq:stokes_discrete} (\(\chunksize-1\) for first,  \(\chunksizelast\) for last chunk)\;
      \For{\(i\), \(\snapshotset_{\chunkindex}\) \KwIn enumerate(\(\snapshotsetpfieldi{\chunkindex}(\param)\), \(\snapshotsetofieldi{\chunkindex}(\param)\), \(\snapshotsetstokesi{\chunkindex}(\param)\), \(\residualsetpfieldi{\chunkindex}(\param)\),
        \(\residualsetofieldi{\chunkindex}(\param)\), \(\residualsetstokesi{\chunkindex}(\param)\)) }{
        \tcp*[l]{Perform a POD for the local snapshots on each MPI rank}
        \(\varepsilon \gets\) \Fntol{\(\varepsilon^\ast\), \(\hapodtreedepth\), \(\abs{\snapshotset_{\chunkindex}}\), \KwFalse}\;
        \texttt{modes}, \texttt{singular\_values} \(\gets\) \(\POD(\snapshotset_{\chunkindex}, \varepsilon)\)\;
        \texttt{modes} \(\multeq\) \texttt{singular\_values}\tcp*[l]{Scale vectors in \texttt{modes} by their corresponding singular value}
        \tcp*[l]{Now collect the modes on the main MPI rank on each compute node.}
        \tcp*[l]{Instead of using MPI rank 0 for all \(i\), we could also use}
        \tcp*[l]{rank \(i\) to perform the PODs on node level in parallel for all \(6\) fields.}
        \texttt{collected\_modes} \(\gets\) \tcp*[l]{Gathered vectors from all \texttt{modes} variables}
        \texttt{num\_snaps} \(\gets\)\tcp*[l]{Sum of \(\abs{\snapshotset_{\chunkindex}}\) over all MPI ranks on node}
        \If{\texttt{node\_rank} = 0}{
          \(\texttt{num\_snaps\_in\_leafs}[i] \pluseq\) \texttt{num\_snaps}\;
          \(\varepsilon =\) \Fntol{\(\varepsilon^\ast\), \(\hapodtreedepth\), \texttt{num\_snaps\_in\_leafs}[i], \KwFalse}\;
          \texttt{input\_for\_node\_level\_POD}[\(i\)].append(\texttt{collected\_modes})\;
          \texttt{modes}, \texttt{singular\_values} \(\gets \POD(\texttt{input\_for\_node\_level\_POD}[i], \varepsilon)\)\;
          \texttt{input\_for\_node\_level\_POD}[\(i\)] = \texttt{modes} * \texttt{singular\_values}
        }
      }
    }
    \tcp*[l]{Now perform the final POD using the collected inputs from all compute nodes.}
    \tcp*[l]{For simplicity, we perform a single POD here. In our actual implementation,}
    \tcp*[l]{we again apply a HAPOD over a binary tree of compute nodes.}
    \texttt{bases} \(\gets\) [\KwNone, \KwNone, \KwNone, \KwNone, \KwNone, \KwNone]\;
    \For{\(i \gets 0\) \KwTo 5}{
      \texttt{collected\_inputs} \(\gets\) \tcp*[l]{Gathered vectors from all \texttt{input\_for\_node\_level\_POD}[\(i\)] variables}
      \texttt{num\_snaps} \(\gets\)\tcp*[l]{Sum of \texttt{num\_snaps\_in\_leafs}[i] over all compute nodes}
      \If{\texttt{world\_rank} = 0}{
        \(\varepsilon =\) \Fntol{\(\varepsilon^\ast\), \(\hapodtreedepth\), \texttt{num\_snaps}, \KwTrue}\;
        \texttt{bases}[i], \_ \(\gets \POD(\texttt{collected\_inputs}, \varepsilon)\)\;
      }
    }
  }
\end{algorithm}%

\section{Numerical experiments}\label{sec:numericalexperiments}
\subsection{Numerical solvers}
Before turning to the model order reduction, we first validate the choice of our numerical solvers
for the high-dimensional problem (see
\cref{sec:solvers}).  As a test case, we use an initially circular cell in the domain \(\domain
= [0,30]^2\) with
microtubules aligned across the x-direction, i.e.
\begin{equation*}
  \pfieldinitial(\spatialvar) = \tanh(\frac{\signeddistance(\spatialvar)}{\sqrt{2}\pfieldepsilon}),\quad
  \ofieldinitial(\spatialvar) = \begin{pmatrix}
    1 \\ 0
  \end{pmatrix}
  \cdot \frac{\pfieldinitial(\spatialvar) + 1}{2},\quad
  \stokesvelocityinitial(\spatialvar) = \Vzero\quad
\end{equation*}
where
\(\signeddistance\) is the signed distance function to the cell membrane (\(\signeddistance(\spatialvar)
= 5 - \norm{{(15, 15)}^T - \spatialvar}\)).
To test the solvers, we compute the first 50 time steps of the first-order
(\(\pfieldpolynomialorder = 1\)) finite element
discretization from \cref{sec:finiteelementdiscretization}
for different grid sizes and different
choices of the linear solvers. For all three systems, we use the UMFPACK~\cite{Davis2004}
solver to directly invert the system matrix.
For the phase field and orientation field jacobians, we also test an iterative GMRES solver which
is taken from dune-istl~\cite{Blatt2007} (\inlinecode{RestartedGMResSolver}) and
used without preconditioner for the
orientation field system. However, the GMRES solver does not converge in the phase field case. We
thus use an incomplete LU
factorization of the matrix~\eqref{eq:phasefieldjacobianpreconditioner} as a
preconditioner (using \inlinecode{IncompleLUT} from the Eigen package~\cite{Guennebaud2010} with a
fill factor of 80).
For the orientation field and Stokes matrices, we also consider a Schur-complement approach (see
\cref{sec:ofieldsolvers} for the orientation field
and, e.g.,~\cite{Benzi2005} for the Stokes system) coupled with an iterative solver. Here we
use GMRES for the non-symmetric Schur complement
of the
orientation field jacobian and CG (dune-istl's \inlinecode{CGSolver}) preconditioned with the
pressure mass matrix~\cite{Peters2005} for the
symmetric positive-definite Stokes Schur complement. All mass matrices needed in the Schur solvers
or for the
preconditioners are inverted using the sparse Cholesky factorization (\inlinecode{SimplicialLDLT}) of
the Eigen package. Since both the mass matrices and
the
preconditioner matrix~\eqref{eq:phasefieldjacobianpreconditioner} are time-independent, the (incomplete) factorization
has to be performed only once before starting the time stepping. During the time-stepping, we can
then efficiently apply the inverse matrix using
this factorization.

The results can be found in \cref{tab:Solvers} and are computed for a single fixed parameter
value of
\(\bendingcapillarynumber=\capillarynumber=\polaritynumber=\activeforcenumber=1\). However, we also
tested several other
parameters to make sure the results do not change significantly with the choice of parameter. The
rows labelled \emph{Setup} contain the
time that is
needed initially for setting up the solver before starting the timestepping. Almost all of this
time is used for the initial factorization of the
matrices. For the phase field and orientation field equations, which are nonlinear, the direct
solvers cannot precompute any
factorization and thus are set up very fast. However, the factorization then has to be done prior
to each solve which is why the average time per
matrix inversion (\emph{Solve} rows of the table) is about an order of magnitude smaller
for the iterative solvers, in particular for the
finer grids. This holds
true even for the orientation field system where the (unpreconditioned) GMRES solver on average
needs over a hundred iterations to converge.
Still, a preconditioner might be needed when refining the grid even further, as the iteration count
increases with the grid size. Fortunately,
the Schur complement seems to be well-conditioned even for the finer grids and has a very low cost
per GMRES iteration as no matrix inversions are
needed for its computation. As a consequence, it is again about an order of magnitude faster than
the full-system GMRES solver. Memory-wise, the
direct solver is much more demanding and uses about thrice as much memory as the two iterative
solvers. If necessary, memory-usage of the iterative
solvers could probably further be reduced by using iterative solvers also for the inversion of the
mass matrices.

For the phase field matrices, preconditioning using the preconditioner~\eqref{eq:phasefieldjacobianpreconditioner}
seems
to work well in practice, though the iterations
slowly increase when refining the grid. However, even on the finest grid tested, only about 14
iterations are needed for convergence, while the
solver does not
converge at all without preconditioning. The increase in the number of iterations might also be
partly due to the fact that we use a fixed time step
of \(\timestep = 10^{-3}\) across all grid sizes such that the ratio of time
step \(\timestep\) to grid width \(\gridwidth\)
increases when refining the grid (remember that we assumed that \(\timestep\) is small compared to
\(\gridwidth\)
when deriving the preconditioner, see \cref{sec:pfieldsolvers}). Though significantly faster than
the
direct solver, the iterative solver is still
the limiting factor for the overall solution process (remember that we have to invert the phase
field jacobian once in each Newton iteration,
while only once per timestep for the linear Stokes system). However, note that all measurements
were done in serial computations. In principle,
parallelization could be used both for the preconditioner and the GMRES solver to reduce
computation times.

The Stokes system matrix is time-independent such that in this case also the direct solver has to
perform only a single matrix factorization, which
is why it is even a little faster than
the iterative Schur complement CG solver. However, the difference is less than a factor of two, and
the iterative solver uses significantly less
memory, which is why we are using the iterative solver for the numerical experiments in this paper.
Of course, there has been a
lot of research on linear solvers for the Stokes system such that more advanced (and probably much
faster)
solvers would be available (see e.g.,~\cite{Benzi2005,Gmeiner2015}). However, as long as the phase field
solver limits overall solution time, the simple
Schur-CG solver is sufficiently fast
for our purposes.

In the following, we will always use the preconditioned GMRES solver for the phase field system,
the unpreconditioned Schur-GMRES solver for the
orientation field
system and the Schur-CG solver preconditioned by the pressure mass matrix for the Stokes system.

\begin{table}[htbp]
  \centering
  \small
  \begin{tabular}{c*{13}{c}}
                              &            & \multicolumn{2}{c}{Phase field solver} & \multicolumn{3}{c}{Orientation field solver} & \multicolumn{2}{c}{Stokes solver} \\
    \cmidrule(r){3-4} \cmidrule(r){5-7} \cmidrule(r){8-9}
    Elements                  &            & Direct & GMRES                         & Direct & GMRES & Schur-GMRES                 & Direct & Schur-CG                 \\
    \midrule
                              & It.        & ---    & 7                             & ---    & 65    & 33                          & ---    & 21                       \\
                              & Setup [s]  & 0.1    & 1.0                           & 0.2    & 0.2   & 0.1                         & 0.5    & 0.1                      \\
                              & Solve [s]  & 0.6    & 0.1                           & 1.1    & 0.3   & 0.0                         & 0.0    & 0.1                      \\
    \multirow{-4}{*}{3,600}   & Mem.\ [GB] & 0.3    & 0.2                           & 0.5    & 0.3   & 0.2                         & 0.3    & 0.2                      \\
    \midrule
                              & It.        & ---    & 8.3                           & ---    & 83    & 26                          & ---    & 22                       \\
                              & Setup [s]  & 0.4    & 6.1                           & 1.1    & 1.1   & 0.6                         & 2.1    & 0.6                      \\
                              & Solve [s]  & 3.4    & 0.5                           & 7.5    & 1.9   & 0.1                         & 0.2    & 0.4                      \\
    \multirow{-4}{*}{14,400}  & Mem.\ [GB] & 1.1    & 0.8                           & 2.0    & 0.9   & 0.8                         & 1.0    & 0.8                      \\
    \midrule
                              & It.        & ---    & 10                            & ---    & 93    & 22                          & ---    & 22                       \\
                              & Setup [s]  & 0.9    & 20                            & 3.0    & 3.0   & 1.8                         & 6.2    & 1.9                      \\
                              & Solve [s]  & 11.3   & 1.4                           & 24.8   & 5.2   & 0.3                         & 0.7    & 1.1                      \\
    \multirow{-4}{*}{32,400}  & Mem.\ [GB] & 2.6    & 1.7                           & 4.9    & 2.1   & 1.7                         & 2.3    & 1.7                      \\
    \midrule
                              & It.        & ---    & 11                            & ---    & 100   & 20                          & ---    & 23                       \\
                              & Setup [s]  & 1.7    & 43                            & 6.6    & 7.1   & 5.1                         & 14.2   & 5.2                      \\
                              & Solve [s]  & 24.9   & 2.8                           & 51.4   & 9.7   & 0.5                         & 1.2    & 2.2                      \\
    \multirow{-4}{*}{57,600}  & Mem.\ [GB] & 4.3    & 3.0                           & 9.2    & 3.7   & 3.0                         & 4.4    & 3.0                      \\
    \midrule
                              & It.        & ---    & 11                            & ---    & 105   & 21                          & ---    & 23                       \\
                              & Setup [s]  & 2.7    & 73                            & 9.8    & 10.0  & 6.8                         & 27.8   & 6.9                      \\
                              & Solve [s]  & 58.4   & 4.5                           & 109.0  & 16.6  & 0.8                         & 2.3    & 3.5                      \\
    \multirow{-4}{*}{90,000}  & Mem.\ [GB] & 8.2    & 4.6                           & 15.1   & 5.7   & 4.6                         & 6.9    & 4.6                      \\
    \midrule
                              & It.        & ---    & 13                            & ---    & 111   & 23                          & ---    & 24                       \\
                              & Setup [s]  & 4.7    & 130                           & 20.9   & 26.5  & 17.9                        & 44.7   & 18.5                     \\
                              & Solve [s]  & 91.2   & 8.0                           & 185.6  & 29.4  & 1.4                         & 3.8    & 6.5                      \\
    \multirow{-4}{*}{129,600} & Mem.\ [GB] & 11.8   & 6.7                           & 22.5   & 8.3   & 6.7                         & 10.1   & 6.6                      \\
    \midrule
                              & It.        & ---    & 14                            & ---    & 114   & 26                          & ---    & 23                       \\
                              & Setup [s]  & 7.5    & 190                           & 38.9   & 32.9  & 28                          & 80.5   & 34.0                     \\
                              & Solve [s]  & 152.4  & 11.5                          & 314.4  & 35.2  & 2.2                         & 5.6    & 9.3                      \\
    \multirow{-4}{*}{176,400} & Mem.\ [GB] & 16.8   & 9.0                           & 31.1   & 11.3  & 9.0                         & 14.3   & 9.0                      \\
    \bottomrule
  \end{tabular}
  \caption{Performance of different solvers in the first 50 time steps of (\(\timestep =
    10^{-3}\)). The number of degrees of freedom in
    the finite
    element approximation
    (and thus the number of
    rows and columns in the system matrix) is about 6 times, 8 times and 4.5 times the number of grid elements, for the phase field,
    orientation field, and Stokes flow, respectively. The direct solver used for the full systems is UMFPACK~\cite{Davis2004}. The
    GMRES solver is
    taken from dune-istl~\cite{Blatt2007} (\inlinecode{RestartedGMResSolver}) and used without preconditioner for the orientation
    field systems (both full system matrix and Schur
    complement, see \cref{sec:ofieldsolvers}). For the phase field system, we use an incomplete LU factorization of the matrix~\eqref{eq:phasefieldjacobianpreconditioner} as a preconditioner
    (using \inlinecode{IncompleLUT} from the Eigen
    package~\cite{Guennebaud2010}). The Stokes schur solver uses a CG solver (dune-istl's
    \inlinecode{CGSolver}) preconditioned with the pressure mass
    matrix~\cite{Peters2005} to invert the Schur complement. All mass matrices needed in the Schur
    solvers or for the
    preconditioners are inverted using the sparse Cholesky factorization (\inlinecode{SimplicialLDLT}) of
    the Eigen package. All computations are serial (no
    parallelization).
    \emph{It.}: Average number of iterations for a single solve, \emph{Setup}: Time needed for initial setup of the solver (mostly factorization of
    matrices). \emph{Solve}: Average (wall) time
    needed
    for a single solve. \emph{Mem.}: Maximum whole-program memory usage (maximum resident
    set size). For the systems that are not
    mentioned in the
    column header, we always use the preconditioned GMRES solver for the phase field system, the Schur solver for the orientation field system and the CG Schur solver preconditioned by the
    pressure mass matrix for the Stokes system.
  }\label{tab:Solvers}
\end{table}

\subsection{Cell isolation experiment}\label{sec:cellisolationexperiment}
\begin{figure}[htbp]
  \centering
  \includegraphics[width=\textwidth]{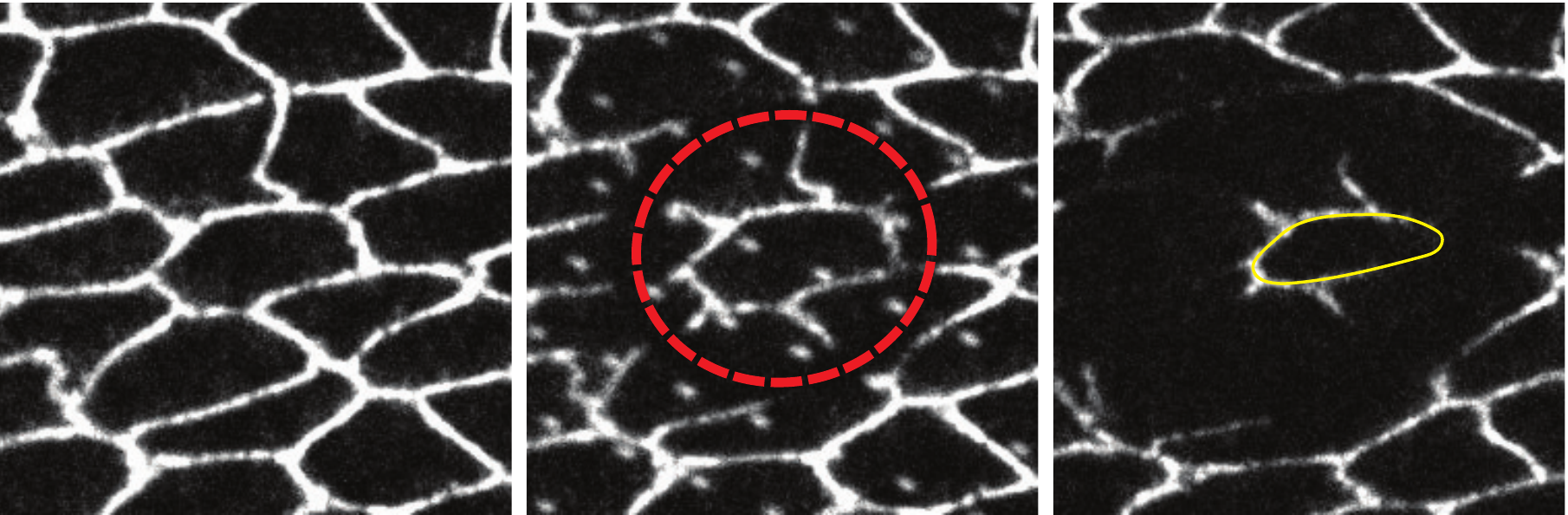}
  \caption{A drosophila wing epithelial cell is isolated from the surrounding tissue by laser
    ablation. The isolated cell mostly
    maintains its shape, the only change is a small isotropic decrease in cell size. This shows that at this stage (16h APF)
    of the wing development,
    cells are not shaped by their neighbors but mostly by cell-autonomous forces
    (see~\cite{Singh2018}).}\label{fig:isolate_cell_by_laser_ablation}
\end{figure}

As a first test for our numerical setup, we consider a cell isolation experiment as in
\cite{Singh2018}.
A drosophila wing epithelial cell is isolated from the surrounding tissue by laser ablation (see
\cref{fig:isolate_cell_by_laser_ablation}).
Intriguingly, the ablation has almost no effect on the cell shape except for a small isotropic
decrease in volume. The authors in \cite{Singh2018}
conclude that the epithelial cells are mostly force-autonomous during this developmental stage.
They further show evidence that there is
a balance between extensile forces generated by microtubules and the contractile forces of the
actomyosin cortex.

\begin{figure}[htbp]
  \centering
  \includegraphics[width=\textwidth]{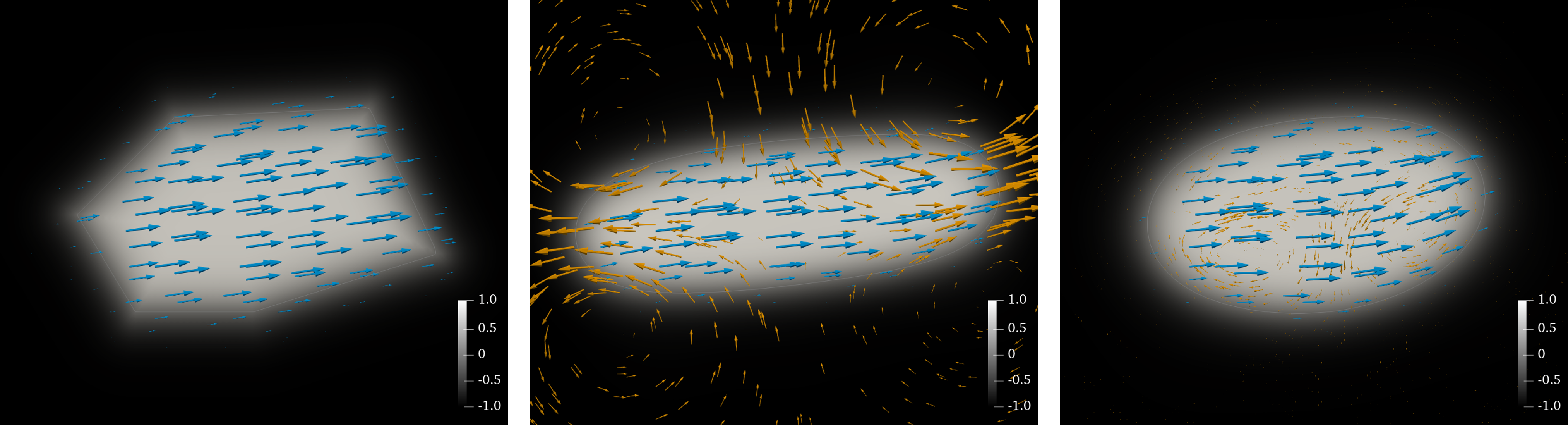}
  \caption{Initial values (left) and solution at time \(\timevar=5\) for \(\capillarynumber = 1\) (center) and
    \(\capillarynumber = 0.2\) (right) for the cell isolation test case. The other parameters are
    chosen as \(\bendingcapillarynumber=\polaritynumber=\activeforcenumber=1\). Numerical solutions
    were computed
    using the first-order (\(\pfieldpolynomialorder=1\)) finite element discretization from
    \cref{sec:finiteelementdiscretization} on a uniform simplicial grid with \(240\times240\) elements and a time step
    of \(\timestep=0.001\). The color bar
    gives the values of the phase field \(\pfield\). The orientation field \(\ofield\) and the fluid
    velocity \(\stokesvelocity\) are represented by
    blue and orange arrows, respectively.}\label{fig:initial_values_numex1}
\end{figure}

To reproduce this setting, we consider a polygonal cell with vertices
\({(10., 18.)}^T\),
\({(13., 13.)}^T\),
\({(19., 13.)}^T\),
\({(28.5, 16.)}^T\),
\({(25., 23.5)}^T\),
\({(15., 23.)}^T\) in
the domain \(\domain = {[0, 40]}^2\). Let again
\(\signeddistance(\spatialvar)\)
be the signed distance of \(\spatialvar\) to the
polygonal cell boundary such that \(\signeddistance(\spatialvar) \geq 0\) if \(\spatialvar\) is inside
the cell and \(\signeddistance(\spatialvar) < 0\) for
\(\spatialvar\) outside of the cell. We then choose the initial condition for the phase field as
\(\pfieldinitial = \tanh\left(\frac{\signeddistance(\spatialvar)}{\sqrt{2}\pfieldepsilon}\right)\). In addition, we
assume that the microtubules
are initially approximately aligned with the elongation of the cell, i.e., \(\ofieldinitial =
{(0.99, 0.14)}^T\), and that the fluid is at rest
(\(\stokesvelocityinitial = \Vzero\)). \cref{fig:initial_values_numex1} shows the initial values and
numerical solutions for
two different parameter choices. As can be seen in the figure, for \(\capillarynumber = 1\), the
cell becomes
elongated and there is still significant fluid flow visible at time \(\timevar = 5\). In contrast,
for \(\capillarynumber = 0.2\),
the cell is almost in steady state at time \(\timevar = 5\). Moreover, the cell approximately
maintains its shape except
that it becomes rounded.

\subsection{Model order reduction}
\subsubsection{HAPOD-DEIM approximation of single fields}
As described in \cref{sec:modelreduction}, we build reduced bases separately for phase field,
orientation field and the Stokes variables. However, since the equations are all coupled, errors in
one of these fields will
also influence the other fields, i.e., an approximation of the phase field variables will also
result in errors in orientation field and Stokes variables and vice versa. In the worst case, these
errors might even be amplified, e.g.,
a small error in the phase field might results in large errors in the orientation field. To test
whether this is the case, we first reduce only one field at a time and analyse the errors that are
introduced in all fields.

We use the cell isolation setup from the previous section. We only vary the two parameters
\(\capillarynumber\) and \(\polaritynumber\) by an order of
magnitude around \(1\) such that the parameter space is \(\paramspace =
[\frac{1}{\sqrt{10}},\sqrt{10}]^2\). We use \(32\) MPI ranks on a single compute
node and sample each parameter uniformly by \(8\) values (\(\ntrain = 64\)) such that each
MPI rank computes solutions for two parameters. To compute the errors, we randomly choose \(32\)
new parameters (\(1\) parameter per rank) and compute both full-order and reduced solutions for
these parameters.
As error measures, we use the absolute and relative mean \(\Lp{2}\) errors
\begin{equation}
  \LpMeanErrorAbs{2} = {\left(\frac{1}{\abs{\snapshotset}} \sum_{\snapshot \in \snapshotset} \Wnorm{\snapshot - \orthprojection[\operatorname{colspan} \MU](\snapshot)}^2\right)}^{\frac12}, \quad
  \LpMeanErrorRel{2} = {\left(\frac{1}{\abs{\snapshotset}} \sum_{\snapshot \in \snapshotset} \frac{\Wnorm{\snapshot - \orthprojection[\operatorname{colspan} \MU](\snapshot)}^2}{\Wnorm{\snapshot}^2}\right)}^{\frac12}.
\end{equation}

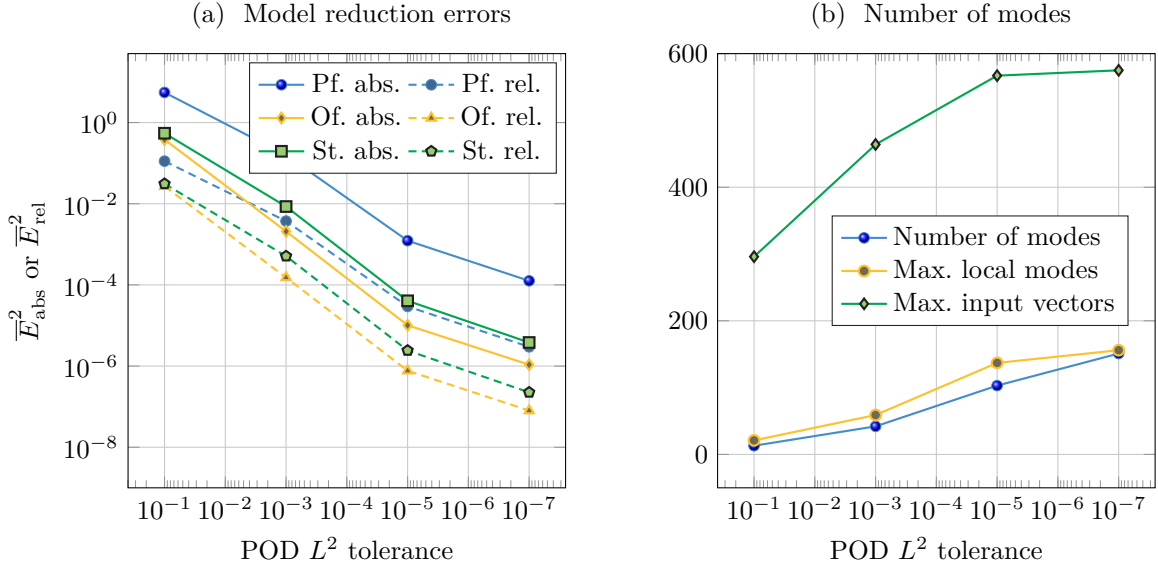
\begin{figure}
  \centering
  \externaltikz{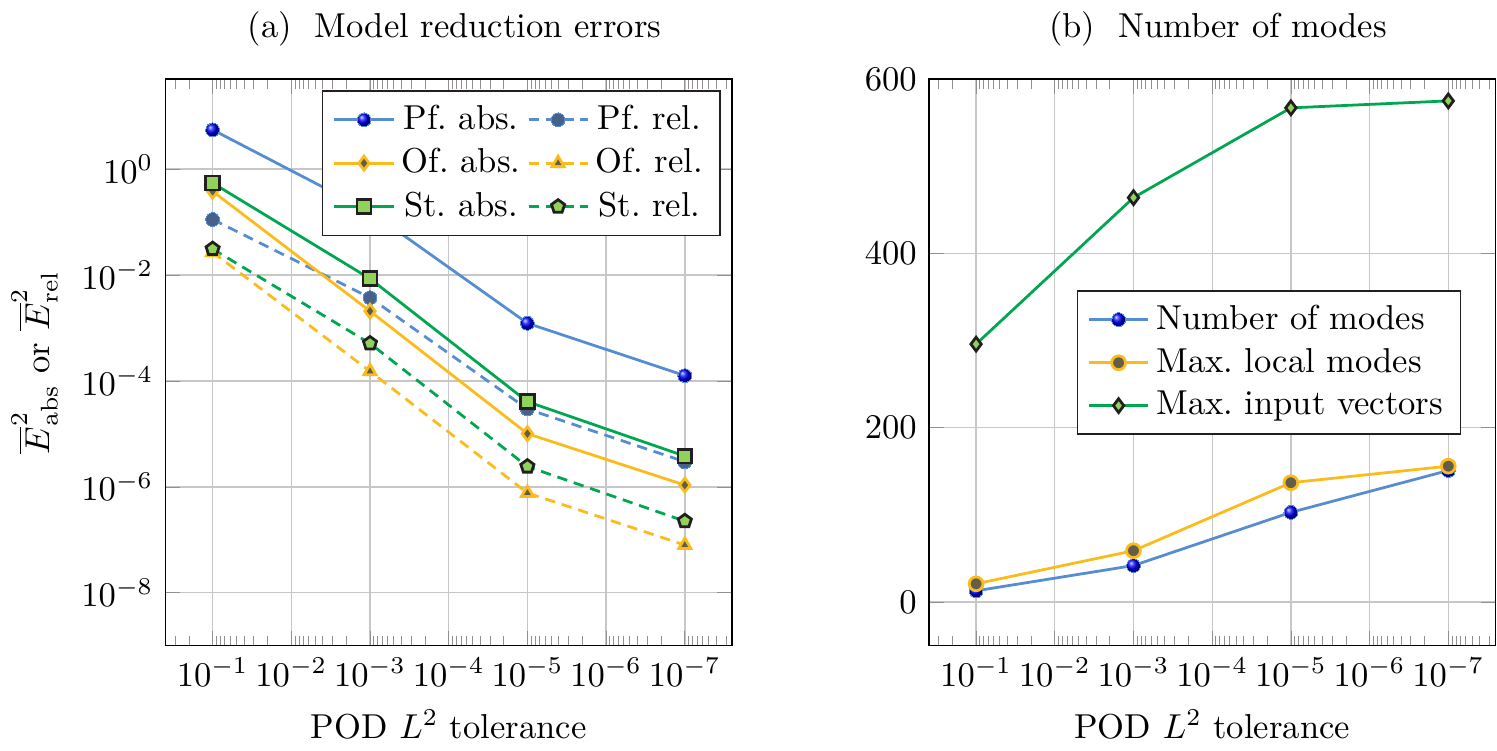}{\begin{tikzpicture}
  \def\pfad{Figures/SingleFieldPODs/}
  \begin{groupplot}[group style={group size=2 by 1, horizontal sep = 2cm,  vertical sep = 2cm},
      width=\figurewidth,
      height=\figureheight,
      scale only axis,
      grid = major,
      xmode = log,
      x dir = reverse,
      xlabel = {POD \(\Lp{2}\) tolerance},
    ]
    \nextgroupplot[
      title = \tikztitle{Model reduction errors},
      ylabel = {\(\LpMeanErrorAbs{2}\) or \(\LpMeanErrorRel{2}\)},
      ymode = log,
      cycle list name=color parula pairwise,
      ymin = 1e-9,
      ymax = 50,
      legend columns = 2,
    ]
    \addplot+ [thick] plot table[x=pod_L2tol,y=pfield_red_err_new] {\pfad pod_only_pfield_no_deim.txt};
    \addplot+ [thick] plot table[x=pod_L2tol,y=pfield_rel_red_err_new] {\pfad pod_only_pfield_no_deim.txt};
    \addplot+ [thick] plot table[x=pod_L2tol,y=ofield_red_err_new] {\pfad pod_only_pfield_no_deim.txt};
    \addplot+ [thick] plot table[x=pod_L2tol,y=ofield_rel_red_err_new] {\pfad pod_only_pfield_no_deim.txt};
    \addplot+ [thick] plot table[x=pod_L2tol,y=stokes_red_err_new] {\pfad pod_only_pfield_no_deim.txt};
    \addplot+ [thick] plot table[x=pod_L2tol,y=stokes_rel_red_err_new] {\pfad pod_only_pfield_no_deim.txt};
    \addlegendentry{Pf.\ abs.};
    \addlegendentry{Pf.\ rel.};
    \addlegendentry{Of.\ abs.};
    \addlegendentry{Of.\ rel.};
    \addlegendentry{St.\ abs.};
    \addlegendentry{St.\ rel.};

    \nextgroupplot[
      title = \tikztitle{Number of modes},
      ylabel = {},
      legend style={at={(0.6, 0.5)},xshift=-0cm,yshift=0cm,anchor=center,nodes=right},
      cycle list name=color parula,
      ymin = -50,
      ymax = 600,
    ]
    \addplot+ [thick] plot table[x=pod_L2tol,y=num_modes] {\pfad pod_only_pfield_no_deim.txt};
    \addplot+ [thick] plot table[x=pod_L2tol,y=max_local_modes] {\pfad pod_only_pfield_no_deim.txt};
    \addplot+ [thick] plot table[x=pod_L2tol,y=max_input_vecs] {\pfad pod_only_pfield_no_deim.txt};
    \addlegendentry{Number of modes};
    \addlegendentry{Max.\ local modes};
    \addlegendentry{Max.\ input vectors};
  \end{groupplot}
\end{tikzpicture}%
}
  \caption{Errors and number of modes used for HAPOD approximation of the phase field with different
    prescribed tolerances
    (cell isolation test case, \(80\times80\) rectangular grid, \(\tf = 0.2\), \(\timestep = 0.001\), \(\omega = 0.95\)).
    (a) Mean \(\Lp{2}\) model reduction error in the different fields.
    (b) Number of modes obtained from the HAPOD.}\label{fig:SingleFieldPODsPfield}
\end{figure}

\begin{figure}
  \centering
  \externaltikz{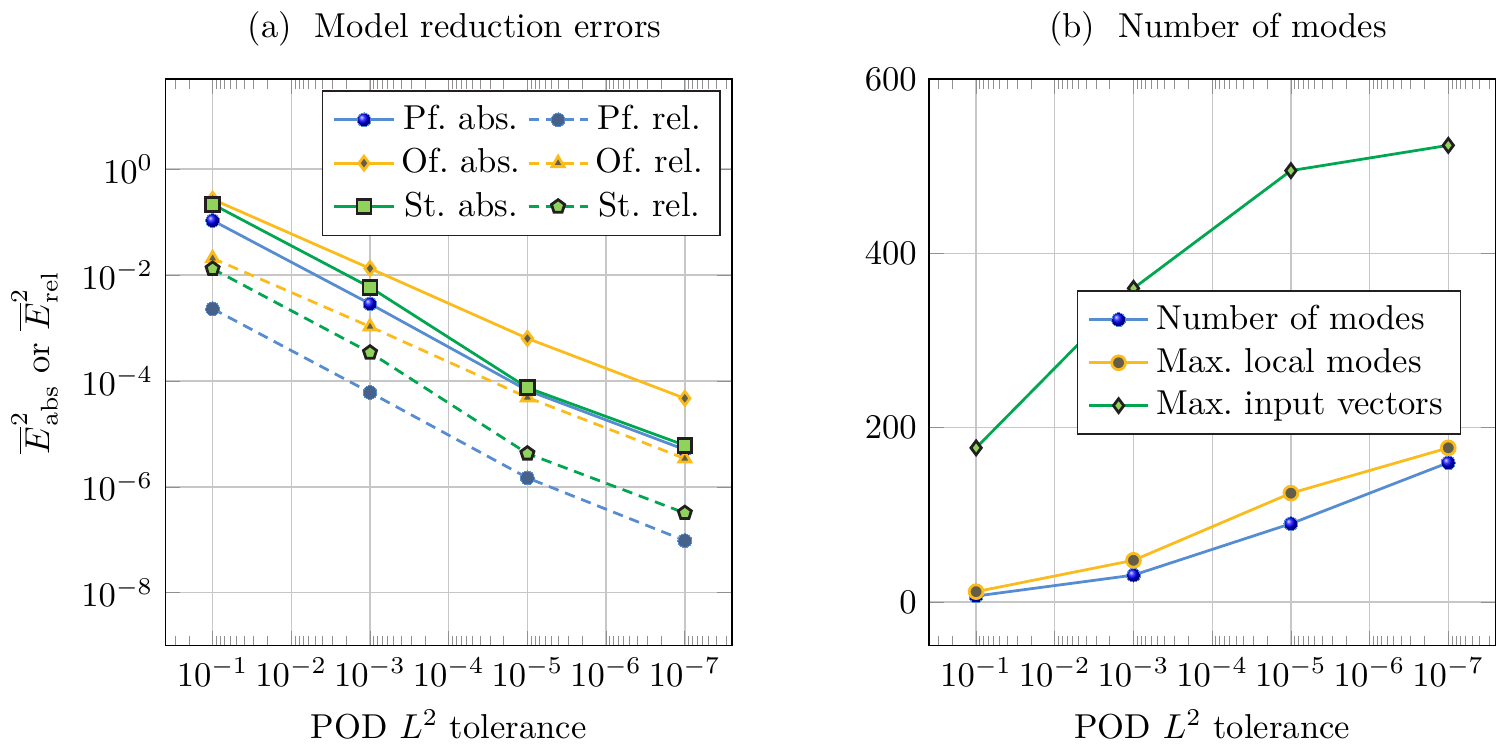}{\begin{tikzpicture}
  \def\pfad{Figures/SingleFieldPODs/}
  \begin{groupplot}[group style={group size=2 by 1, horizontal sep = 2cm,  vertical sep = 2cm},
      width=\figurewidth,
      height=\figureheight,
      scale only axis,
      grid = major,
      xmode = log,
      x dir = reverse,
      xlabel = {POD \(\Lp{2}\) tolerance},
    ]
    \nextgroupplot[
      title = \tikztitle{Model reduction errors},
      ylabel = {\(\LpMeanErrorAbs{2}\) or \(\LpMeanErrorRel{2}\)},
      ymin = 1e-9,
      ymax = 50,
      legend columns = 2,
      ymode = log,
      cycle list name=color parula pairwise,
    ]
    \addplot+ [thick] plot table[x=pod_L2tol,y=pfield_red_err_new] {\pfad pod_only_ofield_no_deim.txt};
    \addplot+ [thick] plot table[x=pod_L2tol,y=pfield_rel_red_err_new] {\pfad pod_only_ofield_no_deim.txt};
    \addplot+ [thick] plot table[x=pod_L2tol,y=ofield_red_err_new] {\pfad pod_only_ofield_no_deim.txt};
    \addplot+ [thick] plot table[x=pod_L2tol,y=ofield_rel_red_err_new] {\pfad pod_only_ofield_no_deim.txt};
    \addplot+ [thick] plot table[x=pod_L2tol,y=stokes_red_err_new] {\pfad pod_only_ofield_no_deim.txt};
    \addplot+ [thick] plot table[x=pod_L2tol,y=stokes_rel_red_err_new] {\pfad pod_only_ofield_no_deim.txt};
    \addlegendentry{Pf.\ abs.};
    \addlegendentry{Pf.\ rel.};
    \addlegendentry{Of.\ abs.};
    \addlegendentry{Of.\ rel.};
    \addlegendentry{St.\ abs.};
    \addlegendentry{St.\ rel.};

    \nextgroupplot[
      title = \tikztitle{Number of modes},
      ylabel = {},
      legend style={at={(0.6, 0.5)},xshift=-0cm,yshift=0cm,anchor=center,nodes=right},
      cycle list name=color parula,
      ymin = -50,
      ymax = 600,
    ]
    \addplot+ [thick] plot table[x=pod_L2tol,y=num_modes] {\pfad pod_only_ofield_no_deim.txt};
    \addplot+ [thick] plot table[x=pod_L2tol,y=max_local_modes] {\pfad pod_only_ofield_no_deim.txt};
    \addplot+ [thick] plot table[x=pod_L2tol,y=max_input_vecs] {\pfad pod_only_ofield_no_deim.txt};
    \addlegendentry{Number of modes};
    \addlegendentry{Max.\ local modes};
    \addlegendentry{Max.\ input vectors};
  \end{groupplot}
\end{tikzpicture}%
}
  \caption{Errors and number of modes used for HAPOD approximation of the orientation field with different
    prescribed tolerances
    (cell isolation test case, \(80\times80\) rectangular grid, \(\tf = 0.2\), \(\timestep = 0.001\), \(\omega = 0.95\)).
    (a) Mean \(\Lp{2}\) model reduction error in the different fields.
    (b) Number of modes obtained from the HAPOD.}\label{fig:SingleFieldPODsOfield}
\end{figure}

\begin{figure}
  \centering
  \externaltikz{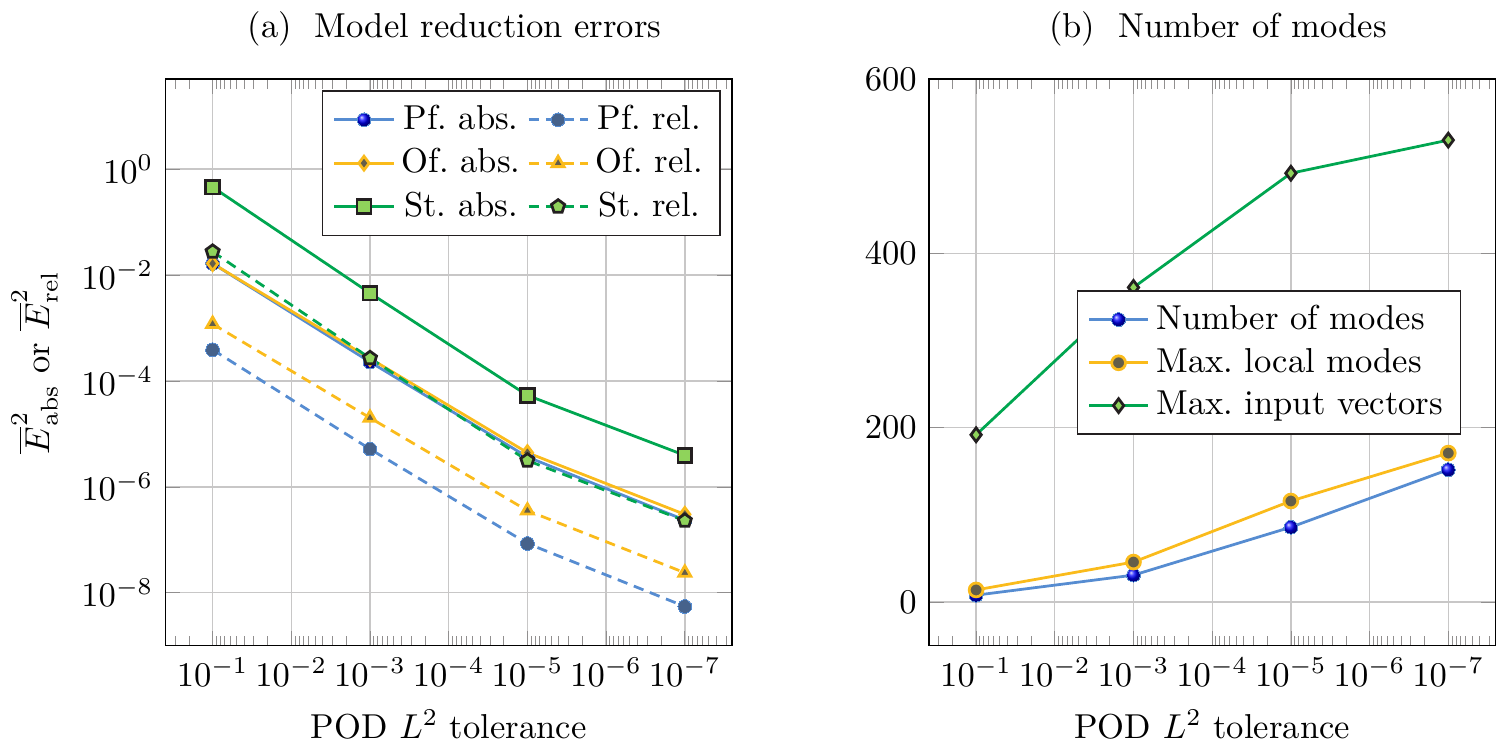}{\begin{tikzpicture}
  \def\pfad{Figures/SingleFieldPODs/}
  \begin{groupplot}[group style={group size=2 by 1, horizontal sep = 2cm,  vertical sep = 2cm},
      width=\figurewidth,
      height=\figureheight,
      scale only axis,
      grid = major,
      xmode = log,
      x dir = reverse,
      xlabel = {POD \(\Lp{2}\) tolerance},
    ]
    \nextgroupplot[
      title = \tikztitle{Model reduction errors},
      ylabel = {\(\LpMeanErrorAbs{2}\) or \(\LpMeanErrorRel{2}\)},
      ymin = 1e-9,
      ymax = 50,
      ymode = log,
      cycle list name=color parula pairwise,
      legend columns = 2,
    ]
    \addplot+ [thick] plot table[x=pod_L2tol,y=pfield_red_err_new] {\pfad pod_only_stokes_no_deim.txt};
    \addplot+ [thick] plot table[x=pod_L2tol,y=pfield_rel_red_err_new] {\pfad pod_only_stokes_no_deim.txt};
    \addplot+ [thick] plot table[x=pod_L2tol,y=ofield_red_err_new] {\pfad pod_only_stokes_no_deim.txt};
    \addplot+ [thick] plot table[x=pod_L2tol,y=ofield_rel_red_err_new] {\pfad pod_only_stokes_no_deim.txt};
    \addplot+ [thick] plot table[x=pod_L2tol,y=stokes_red_err_new] {\pfad pod_only_stokes_no_deim.txt};
    \addplot+ [thick] plot table[x=pod_L2tol,y=stokes_rel_red_err_new] {\pfad pod_only_stokes_no_deim.txt};
    \addlegendentry{Pf.\ abs.};
    \addlegendentry{Pf.\ rel.};
    \addlegendentry{Of.\ abs.};
    \addlegendentry{Of.\ rel.};
    \addlegendentry{St.\ abs.};
    \addlegendentry{St.\ rel.};

    \nextgroupplot[
      title = \tikztitle{Number of modes},
      legend style={at={(0.6, 0.5)},xshift=-0cm,yshift=0cm,anchor=center,nodes=right},
      cycle list name=color parula,
      ymin = -50,
      ymax = 600,
    ]
    \addplot+ [thick] plot table[x=pod_L2tol,y=num_modes] {\pfad pod_only_stokes_no_deim.txt};
    \addplot+ [thick] plot table[x=pod_L2tol,y=max_local_modes] {\pfad pod_only_stokes_no_deim.txt};
    \addplot+ [thick] plot table[x=pod_L2tol,y=max_input_vecs] {\pfad pod_only_stokes_no_deim.txt};
    \addlegendentry{Number of modes};
    \addlegendentry{Max.\ local modes};
    \addlegendentry{Max.\ input vectors};
  \end{groupplot}
\end{tikzpicture}%
}
  \caption{Errors and number of modes used for HAPOD approximation of the Stokes field with different
    prescribed tolerances
    (cell isolation test case, \(80\times80\) rectangular grid, \(\tf = 0.2\), \(\timestep = 0.001\), \(\omega = 0.95\)).
    (a) Mean \(\Lp{2}\) model reduction error in the different fields.
    (b) Number of modes obtained from the HAPOD.}\label{fig:SingleFieldPODsStokes}
\end{figure}

The results without DEIM can be found in \cref{fig:SingleFieldPODsPfield,fig:SingleFieldPODsOfield,fig:SingleFieldPODsStokes}. We see that both relative and
absolute error for all fields are
proportional to the prescribed tolerance, except for the very low tolerance
\(10^{-7}\) (probably due to numerical inaccuracies). In all cases, the
error for the approximated field is about an order of magnitude larger
than the induced error in the other fields. We can thus conclude that the errors are not amplified
between equations. Regarding the absolute values, we see that, at least for the smaller tolerances,
the error (in the Stokes variables) introduced by the Stokes approximation is about an order of
magnitude smaller than the errors introduced by the phase field and orientation field
approximations in their respective variables.
In the following, we thus might have to choose the tolerances for these fields about an order of
magnitude smaller than the Stokes tolerance.

The maximal number of local modes, the maximal number of input vectors to a local POD and the
number of final POD modes are similar for the three fields.

\begin{figure}
  \centering
  \externaltikz{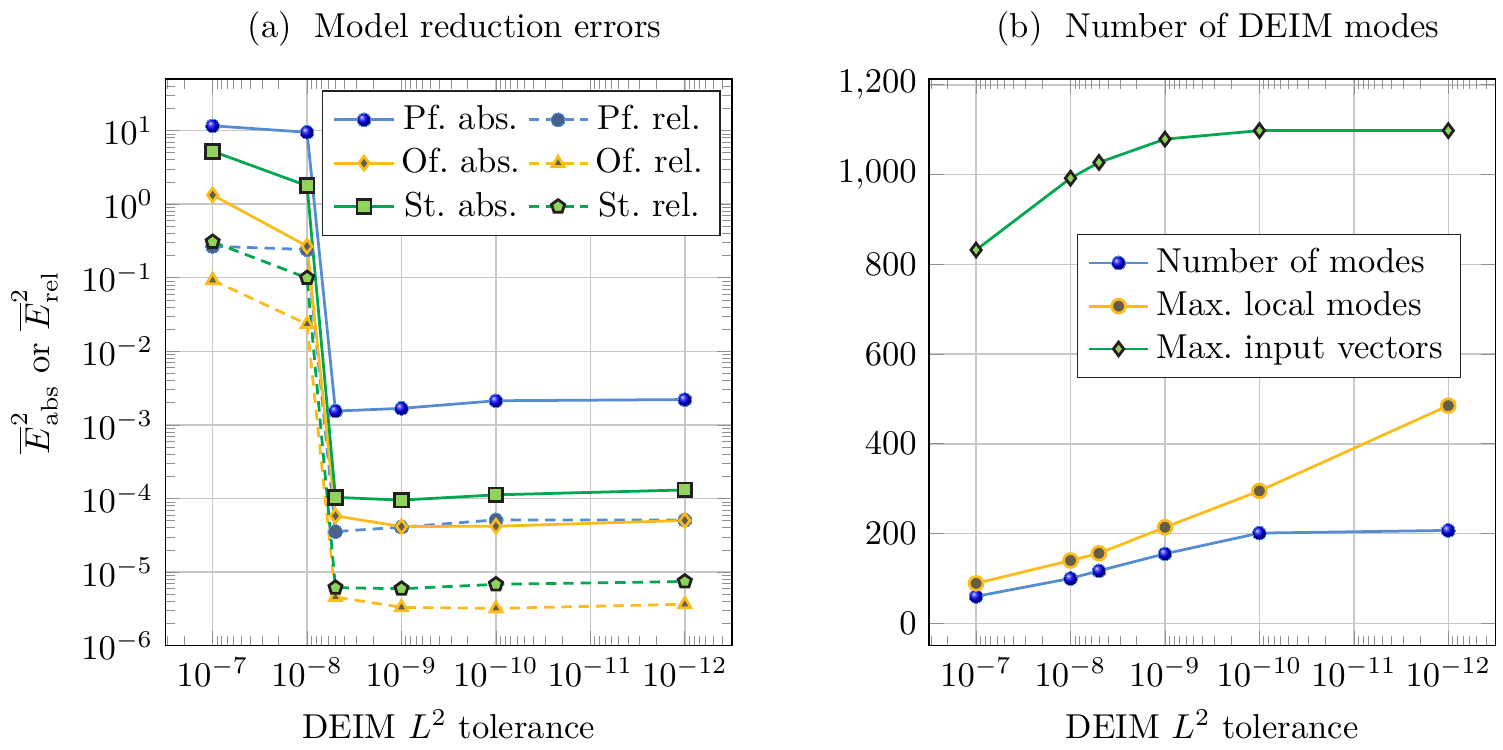}{\begin{tikzpicture}
  \def\pfad{Figures/SingleFieldPODs/}
  \begin{groupplot}[group style={group size=2 by 1, horizontal sep = 2cm,  vertical sep = 2cm},
      width=\figurewidth,
      height=\figureheight,
      scale only axis,
      grid = major,
      xmode = log,
      x dir = reverse,
      xlabel = {DEIM \(\Lp{2}\) tolerance},
    ]
    \nextgroupplot[
      title = \tikztitle{Model reduction errors},
      ylabel = {\(\LpMeanErrorAbs{2}\) or \(\LpMeanErrorRel{2}\)},
      ymode = log,
      cycle list name=color parula pairwise,
      ymin = 1e-6,
      ymax = 50,
      legend columns = 2,
    ]
    \addplot+ [thick] plot table[x=L2tol,y=pfield_red_err_new] {\pfad pod_only_pfield_with_deim.txt};
    \addplot+ [thick] plot table[x=L2tol,y=pfield_rel_red_err_new] {\pfad pod_only_pfield_with_deim.txt};
    \addplot+ [thick] plot table[x=L2tol,y=ofield_red_err_new] {\pfad pod_only_pfield_with_deim.txt};
    \addplot+ [thick] plot table[x=L2tol,y=ofield_rel_red_err_new] {\pfad pod_only_pfield_with_deim.txt};
    \addplot+ [thick] plot table[x=L2tol,y=stokes_red_err_new] {\pfad pod_only_pfield_with_deim.txt};
    \addplot+ [thick] plot table[x=L2tol,y=stokes_rel_red_err_new] {\pfad pod_only_pfield_with_deim.txt};
    \addlegendentry{Pf.\ abs.};
    \addlegendentry{Pf.\ rel.};
    \addlegendentry{Of.\ abs.};
    \addlegendentry{Of.\ rel.};
    \addlegendentry{St.\ abs.};
    \addlegendentry{St.\ rel.};

    \nextgroupplot[
      title = \tikztitle{Number of DEIM modes},
      ylabel = {},
      legend style={at={(0.6, 0.6)},xshift=-0cm,yshift=0cm,anchor=center,nodes=right},
      cycle list name=color parula,
      ymin = -50,
    ]
    \addplot+ [thick] plot table[x=L2tol,y=num_modes] {\pfad pod_only_pfield_with_deim.txt};
    \addplot+ [thick] plot table[x=L2tol,y=max_local_modes] {\pfad pod_only_pfield_with_deim.txt};
    \addplot+ [thick] plot table[x=L2tol,y=max_input_vecs] {\pfad pod_only_pfield_with_deim.txt};
    \addlegendentry{Number of modes};
    \addlegendentry{Max.\ local modes};
    \addlegendentry{Max.\ input vectors};
  \end{groupplot}
\end{tikzpicture}%
}
  \caption{Errors and number of DEIM modes used for HAPOD-DEIM approximation of the phase field using a fixed
    prescribed tolerance of \(10^{-5}\) for the POD and varying DEIM tolerances
    (cell isolation test case, \(80\times80\) rectangular grid, \(\tf = 0.2\), \(\timestep = 0.001\), \(\omega = 0.95\)).
    (a) Mean \(\Lp{2}\) model reduction error in the different fields.
    (b) Number of DEIM modes obtained from the HAPOD.}\label{fig:SingleFieldPODDEIMsPfield}
\end{figure}

\begin{figure}
  \centering
  \externaltikz{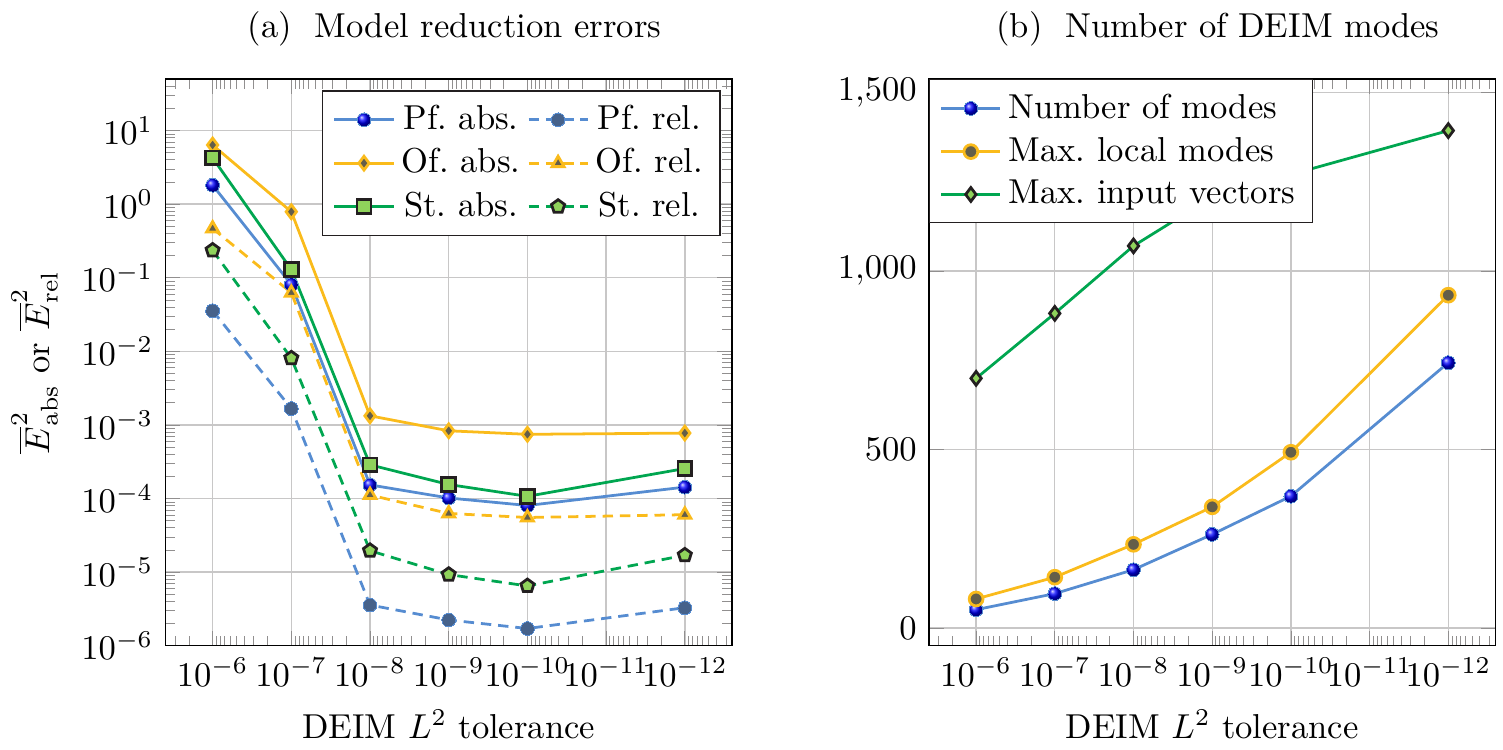}{\begin{tikzpicture}
  \def\pfad{Figures/SingleFieldPODs/}
  \begin{groupplot}[group style={group size=2 by 1, horizontal sep = 2cm,  vertical sep = 2cm},
      width=\figurewidth,
      height=\figureheight,
      scale only axis,
      grid = major,
      xmode = log,
      x dir = reverse,
      xlabel = {DEIM \(\Lp{2}\) tolerance},
    ]
    \nextgroupplot[
      title = \tikztitle{Model reduction errors},
      ylabel = {\(\LpMeanErrorAbs{2}\) or \(\LpMeanErrorRel{2}\)},
      ymode = log,
      cycle list name=color parula pairwise,
      ymin = 1e-6,
      ymax = 50,
      legend columns = 2,
    ]
    \addplot+ [thick] plot table[x=L2tol,y=pfield_red_err_new] {\pfad pod_only_ofield_with_deim.txt};
    \addplot+ [thick] plot table[x=L2tol,y=pfield_rel_red_err_new] {\pfad pod_only_ofield_with_deim.txt};
    \addplot+ [thick] plot table[x=L2tol,y=ofield_red_err_new] {\pfad pod_only_ofield_with_deim.txt};
    \addplot+ [thick] plot table[x=L2tol,y=ofield_rel_red_err_new] {\pfad pod_only_ofield_with_deim.txt};
    \addplot+ [thick] plot table[x=L2tol,y=stokes_red_err_new] {\pfad pod_only_ofield_with_deim.txt};
    \addplot+ [thick] plot table[x=L2tol,y=stokes_rel_red_err_new] {\pfad pod_only_ofield_with_deim.txt};
    \addlegendentry{Pf.\ abs.};
    \addlegendentry{Pf.\ rel.};
    \addlegendentry{Of.\ abs.};
    \addlegendentry{Of.\ rel.};
    \addlegendentry{St.\ abs.};
    \addlegendentry{St.\ rel.};

    \nextgroupplot[
      title = \tikztitle{Number of DEIM modes},
      ylabel = {},
      legend style={at={(0., 1.)},xshift=-0cm,yshift=0cm,anchor=north west,nodes=right},
      cycle list name=color parula,
      ymin = -50,
    ]
    \addplot+ [thick] plot table[x=L2tol,y=num_modes] {\pfad pod_only_ofield_with_deim.txt};
    \addplot+ [thick] plot table[x=L2tol,y=max_local_modes] {\pfad pod_only_ofield_with_deim.txt};
    \addplot+ [thick] plot table[x=L2tol,y=max_input_vecs] {\pfad pod_only_ofield_with_deim.txt};
    \addlegendentry{Number of modes};
    \addlegendentry{Max.\ local modes};
    \addlegendentry{Max.\ input vectors};
  \end{groupplot}
\end{tikzpicture}%
}
  \caption{Errors and number of DEIM modes used for HAPOD-DEIM approximation of the orientation field using a
    fixed prescribed tolerance of \(10^{-5}\) for the POD and varying DEIM
    tolerances
    (cell isolation test case, \(80\times80\) rectangular grid, \(\tf = 0.2\), \(\timestep = 0.001\), \(\omega = 0.95\)).
    (a) Mean \(\Lp{2}\) model reduction error in the different fields.
    (b) Number of DEIM modes obtained from the HAPOD.}\label{fig:SingleFieldPODDEIMsOfield}
\end{figure}

\begin{figure}
  \centering
  \externaltikz{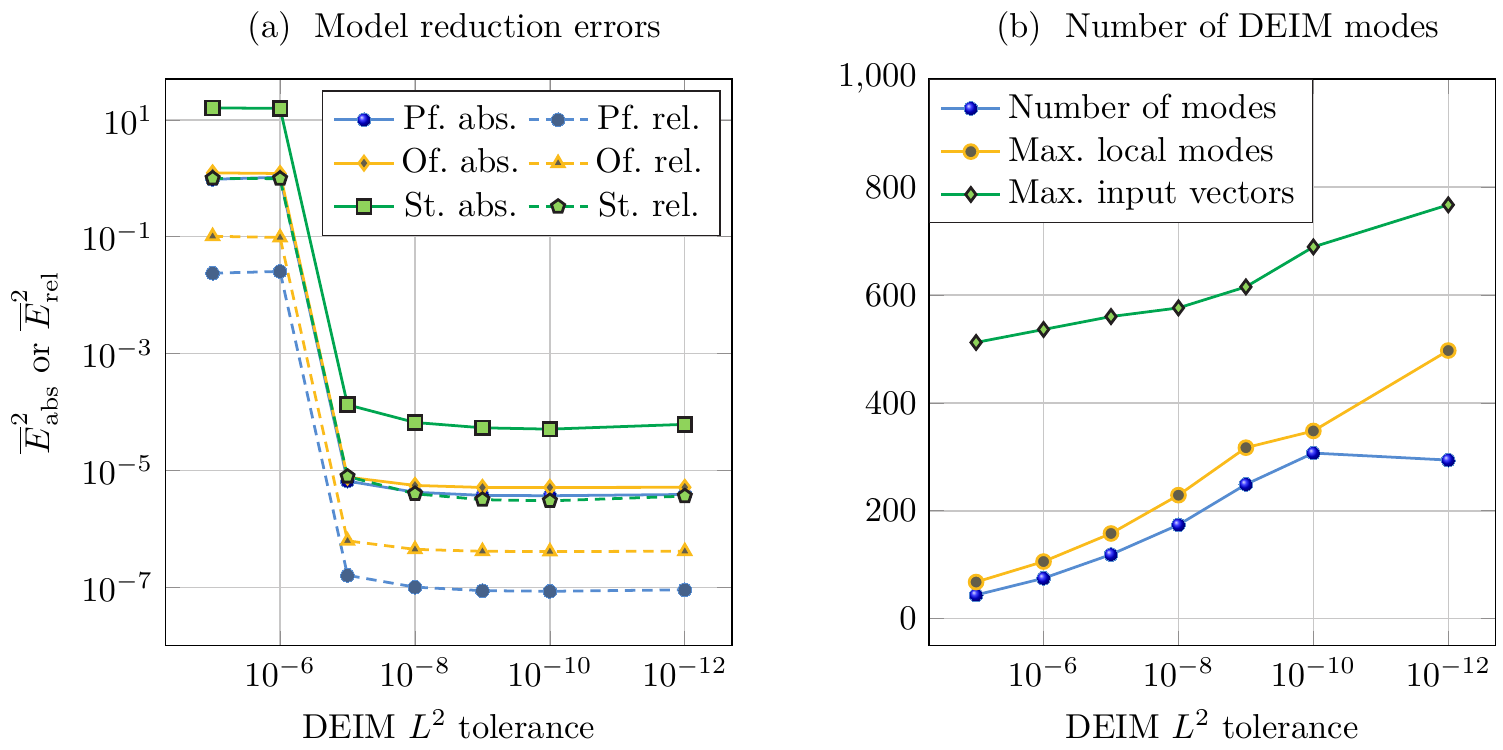}{\begin{tikzpicture}
  \def\pfad{Figures/SingleFieldPODs/}
  \begin{groupplot}[group style={group size=2 by 1, horizontal sep = 2cm,  vertical sep = 2cm},
      width=\figurewidth,
      height=\figureheight,
      scale only axis,
      grid = major,
      xmode = log,
      x dir = reverse,
      xlabel = {DEIM \(\Lp{2}\) tolerance},
    ]
    \nextgroupplot[
      title = \tikztitle{Model reduction errors},
      ylabel = {\(\LpMeanErrorAbs{2}\) or \(\LpMeanErrorRel{2}\)},
      ymode = log,
      cycle list name=color parula pairwise,
      ymin = 1e-8,
      ymax = 50,
      legend columns = 2,
    ]
    \addplot+ [thick] plot table[x=L2tol,y=pfield_red_err_new] {\pfad pod_only_stokes_with_deim.txt};
    \addplot+ [thick] plot table[x=L2tol,y=pfield_rel_red_err_new] {\pfad pod_only_stokes_with_deim.txt};
    \addplot+ [thick] plot table[x=L2tol,y=ofield_red_err_new] {\pfad pod_only_stokes_with_deim.txt};
    \addplot+ [thick] plot table[x=L2tol,y=ofield_rel_red_err_new] {\pfad pod_only_stokes_with_deim.txt};
    \addplot+ [thick] plot table[x=L2tol,y=stokes_red_err_new] {\pfad pod_only_stokes_with_deim.txt};
    \addplot+ [thick] plot table[x=L2tol,y=stokes_rel_red_err_new] {\pfad pod_only_stokes_with_deim.txt};
    \addlegendentry{Pf.\ abs.};
    \addlegendentry{Pf.\ rel.};
    \addlegendentry{Of.\ abs.};
    \addlegendentry{Of.\ rel.};
    \addlegendentry{St.\ abs.};
    \addlegendentry{St.\ rel.};

    \nextgroupplot[
      title = \tikztitle{Number of DEIM modes},
      ylabel = {},
      legend style={at={(0., 1.)},xshift=-0cm,yshift=0cm,anchor=north west,nodes=right},
      cycle list name=color parula,
      ymin = -50,
      ymax = 1000,
    ]
    \addplot+ [thick] plot table[x=L2tol,y=num_modes] {\pfad pod_only_stokes_with_deim.txt};
    \addplot+ [thick] plot table[x=L2tol,y=max_local_modes] {\pfad pod_only_stokes_with_deim.txt};
    \addplot+ [thick] plot table[x=L2tol,y=max_input_vecs] {\pfad pod_only_stokes_with_deim.txt};
    \addlegendentry{Number of modes};
    \addlegendentry{Max.\ local modes};
    \addlegendentry{Max.\ input vectors};
  \end{groupplot}
\end{tikzpicture}%
}
  \caption{Errors and number of DEIM modes used for HAPOD-DEIM approximation of the Stokes variables using a
    fixed prescribed tolerance of \(10^{-5}\) for the POD and varying DEIM
    tolerances
    (cell isolation test case, \(80\times80\) rectangular grid, \(\tf = 0.2\), \(\timestep = 0.001\), \(\omega = 0.95\)).
    (a) Mean \(\Lp{2}\) model reduction error in the different fields.
    (b) Number of DEIM modes obtained from the HAPOD.}\label{fig:SingleFieldPODDEIMsStokes}
\end{figure}
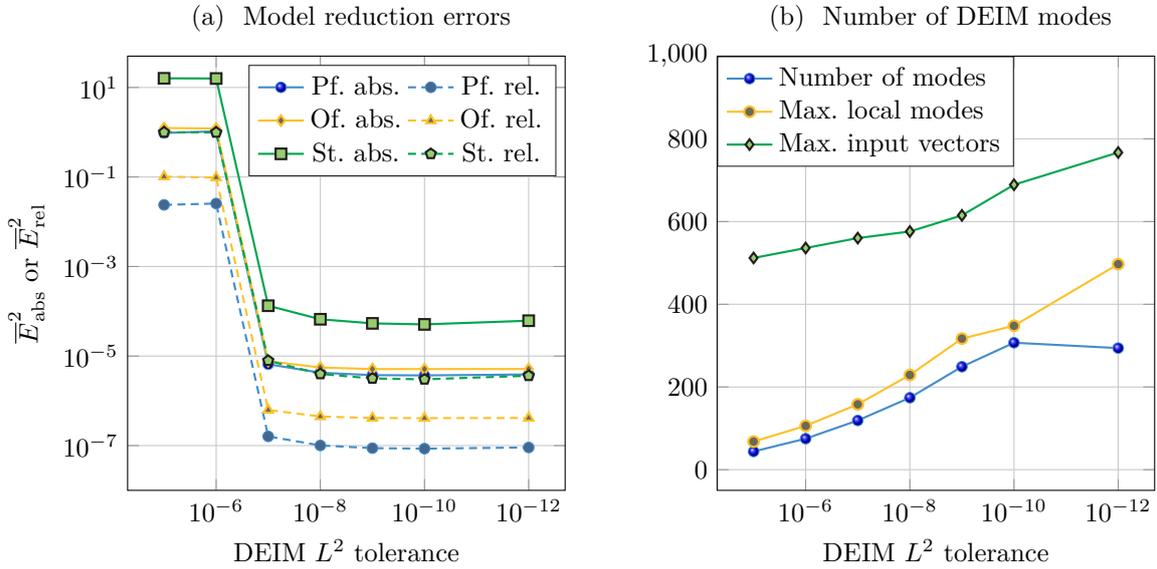

To analyse the effect of the DEIM hyperreduction, we then fix the tolerance for the snapshot POD to
\(10^{-5}\) and perform the same computations again, this time with added
DEIM hyperreduction for different prescribed tolerances.
As can be seen in \cref{fig:SingleFieldPODDEIMsPfield,fig:SingleFieldPODDEIMsOfield,fig:SingleFieldPODDEIMsStokes}, for DEIM tolerance small enough the errors are basically
independent of the DEIM tolerance and
similar to the errors without DEIM hyperreduction
(compare results for a tolerance of \(10^{-5}\) in \cref{fig:SingleFieldPODsPfield,fig:SingleFieldPODsOfield,fig:SingleFieldPODsStokes}).
However, if the DEIM is increased, there is a point where the DEIM interpolation becomes
insufficient and the errors jump by several orders of magnitude. For the phase field, for a POD
tolerance of \(10^{-5}\), this jump occurs at a
DEIM tolerance between \(5\cdot10^{-9}\) and \(10^{-8}\), for the
orientation field between \(\cdot10^{-8}\) and \(10^{-7}\) and
for the Stokes variables between
\(\cdot10^{-7}\) and \(10^{-6}\). These results indicate that the
optimal DEIM tolerance for a POD tolerance of \(10^{-5}\) is in the order of
\(10^{-8}\). For larger tolerances,
large errors can be observed, and further decreasing the tolerance is inefficient since it increase
the number of DEIM modes without improving the model reduction errors.

\section{Conclusion and Outlook}
We suggested to use a nonlinear continuum approach to model microtubule-based forces in planar
fruit fly wing tissues. We presented preconditioned iterative solvers for the discretized model
and outlined an approach to obtain reduced models that can be used, e.g, to make parameter studies
possible where the model has been solved repeatedly for many different parameters.
We proposed to use the hierarchical approximate POD (HAPOD) to efficiently compute the reduced
bases needed
for the model reduction and for an additional hyperreduction by the discrete empirical
interpolation method (DEIM).
Finally, we presented some first numerical studies on the choice of tolerances for the basis
construction showing that a sufficiently efficient approximation can be obtained with a reasonably
small number
of reduced basis functions for each variable.
By now, we only reduced one of the equations at a time. Since the other variables are still
high-dimensional, we cannot expect that the reduced model is significantly faster than the
full-order model in this setting.
To get an actual speed-up, we have to reduce all of the three equations. Preliminary results
indicate that a speed-up of about 100x can realistically be achieved.

\section*{Acknowledgement}
We acknowledge funding by the Deutsche Forschungsgemeinschaft under Germany’s Excellence Strategy
EXC 2044-390685587, Mathematics Münster: Dynamics – Geometry – Structure (T. Leibner, M. Ohlberger,
S. Rave) and  EXC 1003 FF-2015-07 (T. Leibner, M. Matis, M. Ohlberger).

\bibliographystyle{siamplain}
\bibliography{jabreflib}


\end{document}